\documentclass[11pt,twoside]{amsart}
\usepackage{amssymb}
\newenvironment{thlist}{
\renewcommand{\theenumi}{{\rm\roman{enumi}}}
\renewcommand{\labelenumi}{(\theenumi).}
\renewcommand{\labelenumii}{\theenumii)}
\begin{enumerate}}{\end{enumerate}}

\newenvironment{deflist}{
\renewcommand{\theenumi}{\roman{enumi}}
\renewcommand{\theenumii}{\arabic{enumii}}
\renewcommand{\labelenumi}{(\theenumi)}
\renewcommand{\labelenumii}{\theenumii)}
\begin{enumerate}}{\end{enumerate}}
\catcode`@=11

\let\latex@newcommand\newcommand
\let\latex@renewcommand\renewcommand

\def\if@undefined#1{\ifx#1\un@@@defined@@@}
\def\newcommand#1{\if@undefined{#1}
   \let\next@=\latex@newcommand \else
   \let\next@=\latex@renewcommand \fi
   \next@{#1}}

\catcode`@=12

\theoremstyle{definition}

    \newtheorem{thm}{Theorem}[section]
    \newtheorem{prop}[thm]{Proposition}
    \newtheorem{lemma}[thm]{Lemma}
    \newtheorem{defn}[thm]{Definition}
    \newtheorem{cor}[thm]{Corollary}

    \newtheorem*{thm*}{Theorem}
    \newtheorem*{prop*}{Proposition}
    \newtheorem*{lemma*}{Lemma}
    \newtheorem*{defn*}{Definition}
    \newtheorem*{cor*}{Corollary}
    \newtheorem{rems}[thm]{Remark}
    \newtheorem*{rems*}{Remark}
    \newtheorem*{proof*}{Proof}
    \newtheorem*{not*}{Notation}

\renewenvironment{proof}{{\it Proof.}}{$\Box$\par}

\newcommand{\ndef}{\newcommand}
\def\rndef{\renewcommand}

\ndef{\myaddress}[1]{\begin{center} \it\tiny #1 \end{center}}

\ndef{\clA}{{\mathcal A}} \ndef{\rmA}{{\mathrm A}} \ndef{\mbA}{{\mathbb A}} \ndef{\bfA}{{\mathbb A}} \ndef{\euA}{{\EuScript A}} \ndef{\frA}{{\mathfrak A}}
\ndef{\clB}{{\mathcal B}} \ndef{\rmB}{{\mathrm B}} \ndef{\mbB}{{\mathbb B}} \ndef{\bfB}{{\mathbb B}} \ndef{\euB}{{\EuScript B}} \ndef{\frB}{{\mathfrak B}}
\ndef{\clC}{{\mathcal C}} \ndef{\rmC}{{\mathrm C}} \ndef{\mbC}{{\mathbb C}} \ndef{\bfC}{{\mathbb C}} \ndef{\euC}{{\EuScript C}} \ndef{\frC}{{\mathfrak C}}
\ndef{\clD}{{\mathcal D}} \ndef{\rmD}{{\mathrm D}} \ndef{\mbD}{{\mathbb D}} \ndef{\bfD}{{\mathbb D}} \ndef{\euD}{{\EuScript D}} \ndef{\frD}{{\mathfrak D}}
\ndef{\clE}{{\mathcal E}} \ndef{\rmE}{{\mathrm E}} \ndef{\mbE}{{\mathbb E}} \ndef{\bfE}{{\mathbb E}} \ndef{\euE}{{\EuScript E}} \ndef{\frE}{{\mathfrak E}}
\ndef{\clF}{{\mathcal F}} \ndef{\rmF}{{\mathrm F}} \ndef{\mbF}{{\mathbb F}} \ndef{\bfF}{{\mathbb F}} \ndef{\euF}{{\EuScript F}} \ndef{\frF}{{\mathfrak F}}
\ndef{\clG}{{\mathcal G}} \ndef{\rmG}{{\mathrm G}} \ndef{\mbG}{{\mathbb G}} \ndef{\bfG}{{\mathbb G}} \ndef{\euG}{{\EuScript G}} \ndef{\frG}{{\mathfrak G}}
\ndef{\clH}{{\mathcal H}} \ndef{\rmH}{{\mathrm H}} \ndef{\mbH}{{\mathbb H}} \ndef{\bfH}{{\mathbb H}} \ndef{\euH}{{\EuScript H}} \ndef{\frH}{{\mathfrak H}}
\ndef{\clI}{{\mathcal I}} \ndef{\rmI}{{\mathrm I}} \ndef{\mbI}{{\mathbb I}} \ndef{\bfI}{{\mathbb I}} \ndef{\euI}{{\EuScript I}} \ndef{\frI}{{\mathfrak I}}
\ndef{\clJ}{{\mathcal J}} \ndef{\rmJ}{{\mathrm J}} \ndef{\mbJ}{{\mathbb J}} \ndef{\bfJ}{{\mathbb J}} \ndef{\euJ}{{\EuScript J}} \ndef{\frJ}{{\mathfrak J}}
\ndef{\clK}{{\mathcal K}} \ndef{\rmK}{{\mathrm K}} \ndef{\mbK}{{\mathbb K}} \ndef{\bfK}{{\mathbb K}} \ndef{\euK}{{\EuScript K}} \ndef{\frK}{{\mathfrak K}}
\ndef{\clL}{{\mathcal L}} \ndef{\rmL}{{\mathrm L}} \ndef{\mbL}{{\mathbb L}} \ndef{\bfL}{{\mathbb L}} \ndef{\euL}{{\EuScript L}} \ndef{\frL}{{\mathfrak L}}
\ndef{\clM}{{\mathcal M}} \ndef{\rmM}{{\mathrm M}} \ndef{\mbM}{{\mathbb M}} \ndef{\bfM}{{\mathbb M}} \ndef{\euM}{{\EuScript M}} \ndef{\frM}{{\mathfrak M}}
\ndef{\clN}{{\mathcal N}} \ndef{\rmN}{{\mathrm N}} \ndef{\mbN}{{\mathbb N}} \ndef{\bfN}{{\mathbb N}} \ndef{\euN}{{\EuScript N}} \ndef{\frN}{{\mathfrak N}}
\ndef{\clO}{{\mathcal O}} \ndef{\rmO}{{\mathrm O}} \ndef{\mbO}{{\mathbb O}} \ndef{\bfO}{{\mathbb O}} \ndef{\euO}{{\EuScript O}} \ndef{\frO}{{\mathfrak O}}
\ndef{\clP}{{\mathcal P}} \ndef{\rmP}{{\mathrm P}} \ndef{\mbP}{{\mathbb P}} \ndef{\bfP}{{\mathbb P}} \ndef{\euP}{{\EuScript P}} \ndef{\frP}{{\mathfrak P}}
\ndef{\clQ}{{\mathcal Q}} \ndef{\rmQ}{{\mathrm Q}} \ndef{\mbQ}{{\mathbb Q}} \ndef{\bfQ}{{\mathbb Q}} \ndef{\euQ}{{\EuScript Q}} \ndef{\frQ}{{\mathfrak Q}}
\ndef{\clR}{{\mathcal R}} \ndef{\rmR}{{\mathrm R}} \ndef{\mbR}{{\mathbb R}} \ndef{\bfR}{{\mathbb R}} \ndef{\euR}{{\EuScript R}} \ndef{\frR}{{\mathfrak R}}
\ndef{\clS}{{\mathcal S}} \ndef{\rmS}{{\mathrm S}} \ndef{\mbS}{{\mathbb S}} \ndef{\bfS}{{\mathbb S}} \ndef{\euS}{{\EuScript S}} \ndef{\frS}{{\mathfrak S}}
\ndef{\clT}{{\mathcal T}} \ndef{\rmT}{{\mathrm T}} \ndef{\mbT}{{\mathbb T}} \ndef{\bfT}{{\mathbb T}} \ndef{\euT}{{\EuScript T}} \ndef{\frT}{{\mathfrak T}}
\ndef{\clU}{{\mathcal U}} \ndef{\rmU}{{\mathrm U}} \ndef{\mbU}{{\mathbb U}} \ndef{\bfU}{{\mathbb U}} \ndef{\euU}{{\EuScript U}} \ndef{\frU}{{\mathfrak U}}
\ndef{\clV}{{\mathcal V}} \ndef{\rmV}{{\mathrm V}} \ndef{\mbV}{{\mathbb V}} \ndef{\bfV}{{\mathbb V}} \ndef{\euV}{{\EuScript V}} \ndef{\frV}{{\mathfrak V}}
\ndef{\clW}{{\mathcal W}} \ndef{\rmW}{{\mathrm W}} \ndef{\mbW}{{\mathbb W}} \ndef{\bfW}{{\mathbb W}} \ndef{\euW}{{\EuScript W}} \ndef{\frW}{{\mathfrak W}}
\ndef{\clX}{{\mathcal X}} \ndef{\rmX}{{\mathrm X}} \ndef{\mbX}{{\mathbb X}} \ndef{\bfX}{{\mathbb X}} \ndef{\euX}{{\EuScript X}} \ndef{\frX}{{\mathfrak X}}
\ndef{\clY}{{\mathcal Y}} \ndef{\rmY}{{\mathrm Y}} \ndef{\mbY}{{\mathbb Y}} \ndef{\bfY}{{\mathbb Y}} \ndef{\euY}{{\EuScript Y}} \ndef{\frY}{{\mathfrak Y}}
\ndef{\clZ}{{\mathcal Z}} \ndef{\rmZ}{{\mathrm Z}} \ndef{\mbZ}{{\mathbb Z}} \ndef{\bfZ}{{\mathbb Z}} \ndef{\euZ}{{\EuScript Z}} \ndef{\frZ}{{\mathfrak Z}}

\ndef{\tA}{{\widetilde A}} \ndef{\tcA}{{\widetilde\clA}} \ndef{\ttcA}{\widetilde{\tcA}} \ndef{\sfA}{{\textsf A}} \ndef{\ttA}{\widetilde{\tA}} \ndef{\dzA}{{A^sharp}}
\ndef{\tB}{{\widetilde B}} \ndef{\tcB}{{\widetilde\clB}} \ndef{\ttcB}{\widetilde{\tcB}} \ndef{\sfB}{{\textsf B}} \ndef{\ttB}{\widetilde{\tB}} \ndef{\dzB}{{B^sharp}}
\ndef{\tC}{{\widetilde C}} \ndef{\tcC}{{\widetilde\clC}} \ndef{\ttcC}{\widetilde{\tcC}} \ndef{\sfC}{{\textsf C}} \ndef{\ttC}{\widetilde{\tC}} \ndef{\dzC}{{C^sharp}}
\ndef{\tD}{{\widetilde D}} \ndef{\tcD}{{\widetilde\clD}} \ndef{\ttcD}{\widetilde{\tcD}} \ndef{\sfD}{{\textsf D}} \ndef{\ttD}{\widetilde{\tD}} \ndef{\dzD}{{D^sharp}}
\ndef{\tE}{{\widetilde E}} \ndef{\tcE}{{\widetilde\clE}} \ndef{\ttcE}{\widetilde{\tcE}} \ndef{\sfE}{{\textsf E}} \ndef{\ttE}{\widetilde{\tE}} \ndef{\dzE}{{E^sharp}}
\ndef{\tF}{{\widetilde F}} \ndef{\tcF}{{\widetilde\clF}} \ndef{\ttcF}{\widetilde{\tcF}} \ndef{\sfF}{{\textsf F}} \ndef{\ttF}{\widetilde{\tF}} \ndef{\dzF}{{F^sharp}}
\ndef{\tG}{{\widetilde G}} \ndef{\tcG}{{\widetilde\clG}} \ndef{\ttcG}{\widetilde{\tcG}} \ndef{\sfG}{{\textsf G}} \ndef{\ttG}{\widetilde{\tG}} \ndef{\dzG}{{G^sharp}}
\ndef{\tH}{{\widetilde H}} \ndef{\tcH}{{\widetilde\clH}} \ndef{\ttcH}{\widetilde{\tcH}} \ndef{\sfH}{{\textsf H}} \ndef{\ttH}{\widetilde{\tH}} \ndef{\dzH}{{H^sharp}}
\ndef{\tI}{{\widetilde I}} \ndef{\tcI}{{\widetilde\clI}} \ndef{\ttcI}{\widetilde{\tcI}} \ndef{\sfI}{{\textsf I}} \ndef{\ttI}{\widetilde{\tI}} \ndef{\dzI}{{I^sharp}}
\ndef{\tJ}{{\widetilde J}} \ndef{\tcJ}{{\widetilde\clJ}} \ndef{\ttcJ}{\widetilde{\tcJ}} \ndef{\sfJ}{{\textsf J}} \ndef{\ttJ}{\widetilde{\tJ}} \ndef{\dzJ}{{J^sharp}}
\ndef{\tK}{{\widetilde K}} \ndef{\tcK}{{\widetilde\clK}} \ndef{\ttcK}{\widetilde{\tcK}} \ndef{\sfK}{{\textsf K}} \ndef{\ttK}{\widetilde{\tK}} \ndef{\dzK}{{K^sharp}}
\ndef{\tL}{{\widetilde L}} \ndef{\tcL}{{\widetilde\clL}} \ndef{\ttcL}{\widetilde{\tcL}} \ndef{\sfL}{{\textsf L}} \ndef{\ttL}{\widetilde{\tL}} \ndef{\dzL}{{L^sharp}}
\ndef{\tM}{{\widetilde M}} \ndef{\tcM}{{\widetilde\clM}} \ndef{\ttcM}{\widetilde{\tcM}} \ndef{\sfM}{{\textsf M}} \ndef{\ttM}{\widetilde{\tM}} \ndef{\dzM}{{M^sharp}}
\ndef{\tN}{{\widetilde N}} \ndef{\tcN}{{\widetilde\clN}} \ndef{\ttcN}{\widetilde{\tcN}} \ndef{\sfN}{{\textsf N}} \ndef{\ttN}{\widetilde{\tN}} \ndef{\dzN}{{N^sharp}}
\ndef{\tO}{{\widetilde O}} \ndef{\tcO}{{\widetilde\clO}} \ndef{\ttcO}{\widetilde{\tcO}} \ndef{\sfO}{{\textsf O}} \ndef{\ttO}{\widetilde{\tO}} \ndef{\dzO}{{O^sharp}}
\ndef{\tP}{{\widetilde P}} \ndef{\tcP}{{\widetilde\clP}} \ndef{\ttcP}{\widetilde{\tcP}} \ndef{\sfP}{{\textsf P}} \ndef{\ttP}{\widetilde{\tP}} \ndef{\dzP}{{P^sharp}}
\ndef{\tQ}{{\widetilde Q}} \ndef{\tcQ}{{\widetilde\clQ}} \ndef{\ttcQ}{\widetilde{\tcQ}} \ndef{\sfQ}{{\textsf Q}} \ndef{\ttQ}{\widetilde{\tQ}} \ndef{\dzQ}{{Q^sharp}}
\ndef{\tR}{{\widetilde R}} \ndef{\tcR}{{\widetilde\clR}} \ndef{\ttcR}{\widetilde{\tcR}} \ndef{\sfR}{{\textsf R}} \ndef{\ttR}{\widetilde{\tR}} \ndef{\dzR}{{R^sharp}}
\ndef{\tS}{{\widetilde S}} \ndef{\tcS}{{\widetilde\clS}} \ndef{\ttcS}{\widetilde{\tcS}} \ndef{\sfS}{{\textsf S}} \ndef{\ttS}{\widetilde{\tS}} \ndef{\dzS}{{S^sharp}}
\ndef{\tT}{{\widetilde T}} \ndef{\tcT}{{\widetilde\clT}} \ndef{\ttcT}{\widetilde{\tcT}} \ndef{\sfT}{{\textsf T}} \ndef{\ttT}{\widetilde{\tT}} \ndef{\dzT}{{T^sharp}}
\ndef{\tU}{{\widetilde U}} \ndef{\tcU}{{\widetilde\clU}} \ndef{\ttcU}{\widetilde{\tcU}} \ndef{\sfU}{{\textsf U}} \ndef{\ttU}{\widetilde{\tU}} \ndef{\dzU}{{U^sharp}}
\ndef{\tV}{{\widetilde V}} \ndef{\tcV}{{\widetilde\clV}} \ndef{\ttcV}{\widetilde{\tcV}} \ndef{\sfV}{{\textsf V}} \ndef{\ttV}{\widetilde{\tV}} \ndef{\dzV}{{V^sharp}}
\ndef{\tW}{{\widetilde W}} \ndef{\tcW}{{\widetilde\clW}} \ndef{\ttcW}{\widetilde{\tcW}} \ndef{\sfW}{{\textsf W}} \ndef{\ttW}{\widetilde{\tW}} \ndef{\dzW}{{W^sharp}}
\ndef{\tX}{{\widetilde X}} \ndef{\tcX}{{\widetilde\clX}} \ndef{\ttcX}{\widetilde{\tcX}} \ndef{\sfX}{{\textsf X}} \ndef{\ttX}{\widetilde{\tX}} \ndef{\dzX}{{X^sharp}}
\ndef{\tY}{{\widetilde Y}} \ndef{\tcY}{{\widetilde\clY}} \ndef{\ttcY}{\widetilde{\tcY}} \ndef{\sfY}{{\textsf Y}} \ndef{\ttY}{\widetilde{\tY}} \ndef{\dzY}{{Y^sharp}}
\ndef{\tZ}{{\widetilde Z}} \ndef{\tcZ}{{\widetilde\clZ}} \ndef{\ttcZ}{\widetilde{\tcZ}} \ndef{\sfZ}{{\textsf Z}} \ndef{\ttZ}{\widetilde{\tZ}} \ndef{\dzZ}{{Z^sharp}}

\ndef{\bfc}{{\bf c}}

\let\geq\geqslant
\let\leq\leqslant

\ndef{\lims}[1]{\lim\limits_{#1}}
\ndef{\sums}[1]{\sum\limits_{#1}}
\ndef{\ints}[1]{\int\limits_{#1}}
\ndef{\sups}[1]{\sup\limits_{#1}}
\ndef{\liminfty}[1]{\lims{#1\to\infty}}
\ndef{\suminf}[1]{\sums{#1=1}^\infty}

\ndef{\limo}[1]{\omega\mbox{-}\!\!\!\lims{#1\to\infty}}
\ndef{\limL}[1]{\rmL\mbox{-}\!\!\!\lims{#1\to\infty}}
\ndef{\limLOne}[1]{\clL_1\mbox{-}\!\!\!\lims{#1}}
\ndef{\tildelimo}[1]{\tilde\omega\mbox{-}\!\!\!\lims{#1\to\infty}}

\ndef{\normE}[1]{\norm{#1}_E}               
\ndef{\Aut}{\operatorname{Aut}}
\ndef{\Ch}{\operatorname{ch}}        
\ndef{\End}{\operatorname{End}}
\ndef{\Hom}{\operatorname{Hom}}
\ndef{\Ker}{\operatorname{Ker}}
\ndef{\Log}{\operatorname{Log}}
\ndef{\OP}{\operatorname{OP}}
\ndef{\Op}{\operatorname{Op}}
\ndef{\Symb}{\operatorname{Symb}}
\ndef{\Tr}{\operatorname{Tr}}
\ndef{\Wres}{\operatorname{Wres}}
\ndef{\cl}{\operatorname{cl}}
\ndef{\com}{\operatorname{com}}
\ndef{\const}{\operatorname{const}}
\ndef{\conv}{\operatorname{conv}}
\rndef{\det}{\operatorname{det}}
\ndef{\detFK}{\operatorname{det_{FK}}}
\ndef{\diag}{\operatorname{diag}}
\ndef{\dist}{\operatorname{dist}}
\ndef{\dom}{\operatorname{dom}}
\ndef{\ec}{\operatorname{ec}}        
\ndef{\id}{1}
\ndef{\ind}{\operatorname{ind}}
\ndef{\mydeg}{\operatorname{deg}}
\ndef{\op}{\operatorname{op}}
\ndef{\rank}{\operatorname{rank}}
\ndef{\res}{\operatorname{res}}      
\ndef{\rng}{\operatorname{ran}}      
\ndef{\sflow}{\operatorname{sf}}     
\ndef{\sign}{\operatorname{sign}}
\ndef{\sing}{\operatorname{sing}}
\ndef{\supp}{\operatorname{supp}}
\ndef{\tr}{\operatorname{tr}}
\ndef{\vol}{\operatorname{vol}}      
\ndef{\wn}{\operatorname{wn}}        
\ndef{\wres}{\operatorname{wres}}    
\rndef{\Im}{\operatorname{Im}}
\rndef{\Re}{\operatorname{Re}}

\ndef{\rslv}[1]{R_z(#1)}
\ndef{\HH}{H}
\ndef{\tHH}{\tilde \HH}
\ndef{\VV}{V}
\ndef{\Rz}{R_z}
\ndef{\tRz}{\tR_z}
\ndef{\psif}[1]{#1^{[1]}} 
\ndef{\bndl}{\xi}                         
\ndef{\bndlA}{\eta}                      
\ndef{\GlueMap}{\varphi}                        
\ndef{\ChartMap}{h}                             

\ndef{\hilb}{\clH}                     
   \ndef{\hilbasargument}{(\hilb)} 
\ndef{\LpH}[1]{\clL^{#1}}                      
\ndef{\TrCl}{\LpH{1}}                          
\ndef{\TrClH}{\LpH{1}\hilbasargument}          
\ndef{\clBH}{\clB\hilbasargument}              
\ndef{\clCH}{\clC\hilbasargument}              
\ndef{\clKH}{\clK\hilbasargument}              
\ndef{\clFH}{\clF\hilbasargument}              
\ndef{\clUH}{\clU\hilbasargument}              
\ndef{\clCFH}{{\clC\clF}\hilbasargument}       
\ndef{\saBH}{\clB\hilbasargument_{sa}}         
\ndef{\saCH}{\clC\hilbasargument_{sa}}         
\ndef{\saFH}{\clF\hilbasargument_{sa}}         
\ndef{\saKH}{\clK\hilbasargument_{sa}}         
\ndef{\saCFH}{\clC\clF\hilbasargument_{sa}}    
\ndef{\clUFH}{\clU\clF\hilbasargument}         
\ndef{\Uinj}{\clU_{inj}\hilbasargument}        
\ndef{\UFinj}{\clU\clF_{inj}\hilbasargument}   

\ndef{\LpN}[1]{\clL^{#1}(\clN,\tau)}       
\ndef{\rLpN}[1]{L^{#1}(\clN,\tau)}       
\ndef{\clAND}{(\clA,\clN,\clD)}          
\ndef{\clBA}{{\clB(\clA)}}
\ndef{\clKN}{{\clK(\clN,\tau)}}          
\ndef{\clKtN}{{\clK_\tau(\tilde\clN)}}   
\ndef{\clFN}{{\clF_\tau(\clN)}}          
\ndef{\clPN}{\clP(\clN)}                 
\ndef{\clQN}{\clQ_\tau(\clN)}            
\ndef{\infPN}{{\clP_\tau^\infty(\clN)}}  
\ndef{\clOF}[2]{\clF_{#1\mbox{-}#2}(\clN,\tau)}         
\ndef{\oind}[2]{{\rm \tau\mbox{-}ind}_{#1\mbox{-}#2}}   
\ndef{\tind}{\ind_\tau}                  
\ndef{\DInd}{\ind_{\clD,\tau}}           
\ndef{\BF}{$\tau$-Fredholm }
\ndef{\affl}{\eta}
\ndef{\vNa}{von Neumann algebra}
\ndef{\nsf}{faithful semifinite normal }  
\ndef{\taubrs}[1]{\tau\brackets{#1}}

\ndef{\domd}{\bigcap\limits_{n\ge 0} \dom\;\delta^n}          
\ndef{\DiffOP}{{\rm \clD}}
\ndef{\ADA}{\clA \cup [\clD,\clA]}
\ndef{\DixIdeal}[1]{\LpH{#1,\infty}}               
\ndef{\dixideal}{\ell^{1,\infty}}                  
\ndef{\WDixIdeal}{\LpH{1,\mathrm w}}               
\ndef{\DixIdealPos}[1]{\DixIdeal{#1}_+}          
\ndef{\DixIdealN}[1]{\LpN{#1,\infty}}            
\ndef{\DixIdealNPar}[2]{\clL_{#1,\infty}^{#2}(\clN,\tau)}    
\ndef{\DixIdealNPos}[1]{\LpN{#1,\infty}_+}       
\ndef{\TrD}{\Tr_\omega}                                       
\ndef{\tauD}{{\tau_\omega}}                                   
\ndef{\ILog}{\frac 1{\log(1+t)}}
\ndef{\ILogN}{\frac 1{\log(1+N)}}
\ndef{\DixNorm}[1]{\norm{#1}_{(1,\infty)}}                        
\ndef{\DixInt}[1]{\ints 0^t \mu_s(#1)\,ds}
\ndef{\DixIntL}[1]{\ints 0^{\lambda_{1/t}(#1)}\mu_s(#1)\,ds}
    \ndef{\SmallIdeal}{\clL_{1, \mathrm w}}
    \ndef{\DixIntII}[1]{\ints 0^t \mu_s(#1)ds}
    \ndef{\DixIntf}[1]{f_t(#1)}
    \ndef{\DixIntg}[1]{g_t(#1)}

\ndef{\lpi}{\LpN{1,(\pi)}}
\ndef{\HaarMeasBohrs}{\nu}            
\ndef{\BrownsMeas}{\mu}               
\ndef{\BohrCont}[1]{\tilde{#1}}        
\ndef{\APMean}{{M}}                   
\ndef{\CDSS}{{\clA_B}}                
\ndef{\matr}{{\rm Mat}}
\ndef{\seque}[1]{\ensuremath{\{#1_j\}_{j=1}^\infty}}    
\ndef{\sequen}[2]{\ensuremath{\{#1_#2\}_{#2=1}^\infty}}    
\ndef{\Seque}[1]{\ensuremath{\left(#1_0,#1_1,#1_2,\dots\right)}}    
\ndef{\Cesaro}{H}                           
\ndef{\CesaroRPlus}{M}                      
\ndef{\Dilation}{D}                         
\ndef{\Shift}{T}                            
\ndef{\norm}[1]{\left\Vert#1\right\Vert}          
\ndef{\TrNorm}[1]{\norm{#1}_1}              
\ndef{\HSNorm}[1]{\norm{#1}_2}              
\ndef{\InftyNorm}[1]{\norm{#1}_\infty}      
\ndef{\abs}[1]{\left\lvert#1\right\rvert}   
\ndef{\set}[1]{\left\{#1\right\}}           
\ndef{\brackets}[1]{\left(#1\right)}
\ndef{\brs}[1]{\brackets{#1}}
\ndef{\scalprod}[2]{\la #1,#2\ra}
\ndef{\precprec}{\prec\!\!\!\prec}
\ndef{\qeq}{\stackrel?=}
\ndef{\spectrum}[1]{\sigma_{#1}} 
\ndef{\numrange}[1]{\mathrm{W}(#1)}
\rndef{\emptyset}{\varnothing}
\ndef{\sss}[1]{\subsubsection{}\label{#1}}
\rndef{\phi}{\varphi}
\ndef{\OpenUnitDisk}{D}
\ndef{\RHS}{RHS}
\ndef{\LHS}{LHS} 
\ndef{\ttt}{\Leftrightarrow}
\ndef{\then}{\Rightarrow}
\ndef{\tto}{\longrightarrow}
\ndef{\nno}{\nonumber\\}
\ndef{\newn}[1]{\index{#1} \emph{#1}}                                
\ndef{\la}{\langle}
\ndef{\ra}{\rangle}
\ndef{\dbar}{{\;\bar{\phantom{o}} \!\!\!\! d}}
\ndef{\stl}[1]{\stackrel{\vbox to 0pt{\vss\hbox{$\scriptstyle #1$}}}}
\ndef{\mathcomment}[1]{{\scriptstyle\text{(#1)}}\qquad}        
\ndef{\details}[1]{\smallskip\begin{center} {\bf Here:} #1\end{center}\medskip}

     \ndef{\npartial}{\slash\!\!\!\partial}
     \ndef{\Heis}{\operatorname{Heis}}
     \ndef{\Solv}{\operatorname{Solv}}
     \ndef{\Spin}{\operatorname{Spin}}
     \ndef{\SO}{\operatorname{SO}}
     \ndef{\Index}{\operatorname{index}}
             \ndef{\coker}{{\mbox coker}}
             \ndef{\p}{\partial}
             \ndef{\dd}{|\clD|}
             \ndef{\n}{\parallel}
     \ndef{\gf}[2]{\genfrac{}{}{0pt}{}{#1}{#2}}
     \ndef{\ta}{\widetilde{\alpha}}
     \ndef{\tb}{\widetilde{\beta}}
     \ndef{\txi}{\widetilde{\xi}}
     \ndef{\tk}{\widetilde{K}}
     \ndef{\CGh}{\widetilde{\CG}}
     \ndef{\boe}{{\bf e}}\ndef{\bt}{{\bf t}}
     \ndef{\vth}{\vartheta}
     \ndef{\db}{\overline{\partial}}
     \ndef{\hV}{\hat{V}}
     \ndef{\cag}{{\clA^\Gamma}}
     \ndef{\sind}{\sigma{\rm -ind}}

\begin{document}
\theoremstyle{definition}

\def\clA{{\mathcal A}}    \def\rmA{{\mathrm A}}    \def\mbA{{\mathbb A}}
\def\clB{{\mathcal B}}    \def\rmB{{\mathrm B}}    \def\mbB{{\mathbb B}}
\def\clC{{\mathcal C}}    \def\rmC{{\mathrm C}}    \def\mbC{{\mathbb C}}
\def\clD{{\mathcal D}}    \def\rmD{{\mathrm D}}    \def\mbD{{\mathbb D}}
\def\clE{{\mathcal E}}    \def\rmE{{\mathrm E}}    \def\mbE{{\mathbb E}}
\def\clF{{\mathcal F}}    \def\rmF{{\mathrm F}}    \def\mbF{{\mathbb F}}
\def\clG{{\mathcal G}}    \def\rmG{{\mathrm G}}    \def\mbG{{\mathbb G}}
\def\clH{{\mathcal H}}    \def\rmH{{\mathrm H}}    \def\mbH{{\mathbb H}}
\def\clI{{\mathcal I}}    \def\rmI{{\mathrm I}}    \def\mbI{{\mathbb I}}
\def\clJ{{\mathcal J}}    \def\rmJ{{\mathrm J}}    \def\mbJ{{\mathbb J}}
\def\clK{{\mathcal K}}    \def\rmK{{\mathrm K}}    \def\mbK{{\mathbb K}}
\def\clL{{\mathcal L}}    \def\rmL{{\mathrm L}}    \def\mbL{{\mathbb L}}
\def\clM{{\mathcal M}}    \def\rmM{{\mathrm M}}    \def\mbM{{\mathbb M}}
\def\clN{{\mathcal N}}    \def\rmN{{\mathrm N}}    \def\mbN{{\mathbb N}}
\def\clO{{\mathcal O}}    \def\rmO{{\mathrm O}}    \def\mbO{{\mathbb O}}
\def\clP{{\mathcal P}}    \def\rmP{{\mathrm P}}    \def\mbP{{\mathbb P}}
\def\clQ{{\mathcal Q}}    \def\rmQ{{\mathrm Q}}    \def\mbQ{{\mathbb Q}}
\def\clR{{\mathcal R}}    \def\rmR{{\mathrm R}}    \def\mbR{{\mathbb R}}
\def\clS{{\mathcal S}}    \def\rmS{{\mathrm S}}    \def\mbS{{\mathbb S}}
\def\clT{{\mathcal T}}    \def\rmT{{\mathrm T}}    \def\mbT{{\mathbb T}}
\def\clU{{\mathcal U}}    \def\rmU{{\mathrm U}}    \def\mbU{{\mathbb U}}
\def\clV{{\mathcal V}}    \def\rmV{{\mathrm V}}    \def\mbV{{\mathbb V}}
\def\clW{{\mathcal W}}    \def\rmW{{\mathrm W}}    \def\mbW{{\mathbb W}}
\def\clX{{\mathcal X}}    \def\rmX{{\mathrm X}}    \def\mbX{{\mathbb X}}
\def\clY{{\mathcal Y}}    \def\rmY{{\mathrm Y}}    \def\mbY{{\mathbb Y}}
\def\clZ{{\mathcal Z}}    \def\rmZ{{\mathrm Z}}    \def\mbZ{{\mathbb Z}}
\def\clNt{\widetilde{\clN}}
\def\nt{(\clN,\tau)}
\def\ent{E\nt}
\def\clMt{\widetilde{\clM}}
\def\LL{{\mathbb L}}
\def\R{{\mathbb R}}
\def\frB{{\mathfrak B}}

\def\lims#1{\lim\limits_{#1}}
\def\sums#1{\sum\limits_{#1}}
\def\ints#1{\int\limits_{#1}}
\def\sups#1{\sup\limits_{#1}}
\def\liminfty#1{\lim\limits_{#1\to\infty}}

\def\limo#1{\omega\mbox{-}\!\!\!\lim\limits_{#1\to\infty}} 
\def\limL#1{\rmL\mbox{-}\!\!\!\lim\limits_{#1\to\infty}} 
\def\tildelimo#1{\tilde\omega\mbox{-}\!\!\!\lim\limits_{#1\to\infty}}

\def\domd{\bigcap\limits_{n\ge 0} \dom\;\delta^n}          
\def\DiffOP{{\rm \clD}}
\def\ADA{\clA \cup [\clD,\clA]}
\def\DixIdeal#1{\clL^{(#1,\infty)}(\hilb)}               
\def\DixIdealPos#1{\clL_+^{(#1,\infty)}(\hilb)}          
\def\DixIdealN#1{\clL^{(#1,\infty)}(\clN,\tau)}            
\def\DixIdealNPar#1#2{\clL^{({#1},\infty)}_{#2}(\clN,\tau)}    
\def\DixIdealNPos#1{\clL^{(#1,\infty)}_+(\clN,\tau)}       
\def\tauD{{\tau_\omega}}                                   
\def\ILog{\frac 1{\log(1+t)}}
\def\ILogN{\frac 1{\log(1+N)}}
\def\DixNorm#1{||#1||_{(1,\infty)}}                        
\def\DixInt#1{\ints 0^t \mu_s(#1)\,ds}
\def\DixIntL#1{\ints 0^{\l_{1/t}(#1)}\mu_s(#1)\,ds}

\def\seque#1{\ifmmode \{#1_j\}_{j=1}^\infty \else $\{#1_j\}_{j=1}^\infty$
\fi}    
\def\sequen#1#2{\ifmmode \{#1_#2\}_{#2=1}^\infty \else
$\{#1_#2\}_{#2=1}^\infty$ \fi}    
\def\Seque#1{\ifmmode \left(#1_1,#1_2,#1_3,\dots\right) \else
$\left(#1_1,#1_2,#1_3,\dots\right)$ \fi}    
\def\Cesaro{H}                     
\def\CesaroRPlus{M}                     
\def\Dilation{D}                   
\def\Shift{T}                                 
\def\then{\Rightarrow}
\def\newn#1{\index{#1} \emph{#1}}
\def\la{\langle} \def\ra{\rangle}

\def\title#1{\begin{center}\bf\large #1\end{center}\vskip 0.3 in}
\def\author#1{\begin{center}#1\end{center}}

     \title{DIXMIER TRACES AND SOME APPLICATIONS \\ IN NONCOMMUTATIVE GEOMETRY}

\hbox{\hsize=0.5\hsize\vbox{\noindent
{\bf Alan L. Carey}\\Mathematical Sciences Institute\\
Australian National University\\
Canberra, ACT. 0200, AUSTRALIA\\
e-mail: acarey@maths.anu.edu.au\\
}\hfill
\vbox{\noindent
{\bf Fyodor A. Sukochev}\\
School of Informatics and Engineering\\
Flinders University\\
Bedford Park S.A 5042 AUSTRALIA\\
e-mail: sukochev@infoeng.flinders.edu.au\\
} }

\tableofcontents

\newpage\section{Introduction}

The Dixmier trace $\tau_\omega$ arose from the problem of whether
the algebra $B(\clH)$ of all bounded linear operators  on a
Hilbert space $\clH$ had a unique non-trivial trace. Dixmier
resolved this question in the negative in 1966 in a note in
Comptes Rendus~\cite{Dix66CR}. To construct the trace he used an
invariant mean $\omega$ on the solvable `ax+b' group. His trace
vanishes on the ideal of trace class operators and hence is
completely disjoint from the usual trace. It is also non-normal.

Applications of the Dixmier trace to classical geometry are
facilitated by a remarkable result \cite{Co6} relating the trace
to the Wodzicki~\cite{Wo84KT} residue for pseudo-differential
operators on a closed manifold. The latter was introduced
independently by Adler \cite{Ad}, Manin \cite{Man}  in the one
dimensional case and developed further by Wodzicki \cite{Wo84KT}
and Guillemin \cite{Gu} in higher dimensions. The intriguing and
most useful property of the Wodzicki residue apart from its
computability in examples is the fact that it makes sense for
pseudo-differential operators of arbitrary order.  Moreover it is
the unique trace on the pseudo-differential operators which
extends the Dixmier trace on operators of order $\leq$ minus the
dimension of the underlying manifold \cite{Co6}.

Both the Dixmier trace and the Wodzicki residue play important
roles in noncommutative geometry and its
applications~\cite{Co4,GVF}. In particular, the Dixmier trace
$\tau_\omega$ is an appropriate noncommutative analogue of
integration on a compact $n$-dimensional Riemannian spin manifold
$M,$ more exactly
\begin{equation}
   \tau_\omega(f(x)|D|^{-n})=\frac 1{(2\pi)^n n} 2^{[n/2]}
\Omega_{n-1} \int_M f(x) d\vol, \label{integration}
\end{equation}
where $f$ is
a smooth function on $M,$ $D$ is the Dirac operator on $M,$
$\Omega_{n-1}$ is the volume of $n-1$-dimensional unit sphere and
$d\vol$ is the Riemannian volume form. The usefulness of the
Dixmier trace is extended by the results of Connes \cite{Co4}
which relates it to residues of zeta functions.

Important applications which we do not have the space to include
are to the theory of gravitation, classical field theory and
particle physics. The former
  is well covered in the book \cite{GVF} and has its origins in the
  relationship of the Wodzicki residue of powers of the Dirac operator
  to the Einstein-Hilbert action
\cite{Co6}.
  Equally we do not try to cover other
  material in the books \cite{Co4} and \cite{GVF} which describe
  several interesting physical applications of the Dixmier trace.
We refer the reader to the survey \cite{Co7} and the extensive
literature on the application of noncommutative
  geometry to the `standard model' of particle physics where, starting with
  an elegant expression for a noncommutative action principle expressed
  in terms of the Dixmier trace,
Connes and Lott \cite{CoL} and Connes \cite{Co8} show how to
derive from it the Euclidean version
  of the action for the standard model.
Analogously, one may also give a noncommutative formulation of the
Hamiltonian version of classical field theory again using the
Dixmier trace \cite{Kal}. Some background to these developments is
provided by \cite{Kas}.

These applications to physics  can be traced back to the research
announcement \cite{Co12} where a series of foundational results
are described (some of these were later collected in \cite{Co4}).
After calculating the Dixmier trace for pseudodifferential
operators Connes  obtained the Yang-Mills and Polyakov actions
from an action functional involving the Dixmier trace.
  The theorem on a residue formula for the Hochschild class
of the Chern character, which we will describe in Section
17, was also announced. The next step appears in \cite{Co11} where
  Connes introduced the axioms of noncommutative spin
geometry for a noncommutative algebra $\clA$. The very first axiom
uses the Dixmier trace to introduce a noncommutative integration
theory on $\clA$ which is completely natural in view of
(\ref{integration}).
  This point of view resurfaces in our discussion of Lesch's
index theorem in Section 16 although we will not digress further
to introduce the other axioms of Connes except to mention in
Section 17 the role of the Hochschild class.

While these applications to physical theory form a motivational
background they are not the focus of this article. Our aim is to
give a unified and coherent account of some recent functional
analytic advances in the theory of Dixmier traces.
In addition to surveying these new results we also offer
in some cases new proofs.
Part of our motivation is to extend and clarify questions
raised by \cite{Co4}
  Chapter IV.
  Specifically we characterise the class of measurable
operators defined in \cite{Co4}, explain the role of the Cesaro
mean in Connes' version of Dixmier traces (called Connes-Dixmier
traces here) and give a complete analysis of zeta function
formulae for the Dixmier trace.  Most importantly, however, our
whole treatment is within the framework of 'semifinite spectral
triples' which we explain in Section 8. This notion arises when
one extends the theory in \cite{Co4}, which deals with subalgebras
of the bounded operators on a separable Hilbert space equipped
with the standard trace, to the situation where the von Neumann
algebra of bounded operators is replaced by a general semifinite
von Neumann algebra equipped with some faithful normal semifinite
trace (cf \cite{A}). As a result of the extra effort needed to handle this level
of generality it is possible to find improvements even in the
standard type $I$ theory of \cite{Co4}. A number of authors have
contributed to the development of this framework and new
applications of Dixmier and more general singular traces \cite
{AGPS}, \cite{AGPS1}, \cite {AGPS2}, \cite{AzSu}, \cite{BeF}, \cite{Co5}
\cite{CP1}, \cite{CP2}, \cite{CPRS1}, \cite{CPRS2}, \cite{CPRS3},
\cite{CPS1}, \cite{CPS2}, \cite{DPSS}, \cite {DFWW}, \cite {GI1},
\cite {GI2}, \cite {KS}, \cite {Pr}, \cite {T1}, \cite {T2}, \cite
{Suk}.

Our exposition is organized around three issues. The first issue
dominates the early part of the article (culminating in Section 6)
and gives the characterization of measurable operators in the
sense of \cite{Co4} Section IV.2.  We approach the topic from the
study of the general theory of singular symmetric functionals
(this and its preliminaries occupy Sections 2 to 5). Our
exposition is based on the approach articulated in~\cite{DPSS,
CPS2}, which considers such traces as a special class of
continuous linear functionals on the corresponding operator
ideals.  Many features of the theory may be well understood, even
in the most trivial situation, when the von Neumann algebras in
question are commutative. In this situation, the theory of Dixmier
traces roughly corresponds to the theory of symmetric functionals
on rearrangement invariant function
spaces~\cite{DPSSS1,DPSSS2,DPSS} and allows an alternative
treatment based on the methods drawn from real analysis.

In preparation for the rest of the paper we then introduce the
notions of semifinite spectral triple, type $II$ spectral flow and
some notation for cyclic cohomology. This leads us into the second
issue, occupying Sections 10 and 12, where we describe the
expression of the Dixmier trace in terms of residues of zeta
functions. The key observation, made in Section 4, that enables us
to prove a considerable generalisation of Proposition IV.2.4 of
\cite{Co4} and also the measurability theorem of Section 6,  is
the existence of two kinds of Dixmier trace. One is associated
with the multiplicative group of the positive reals and its
invariant mean and the second, which is naturally associated with
the zeta function, arises from the invariant mean of the additive
group of the reals. The relationship between the two captures the
formula for the Dixmier trace in terms of the zeta function.

These results use the language of spectral triples and lead
naturally to the third issue, namely some applications, which
begin in Section 11. We give a proof of (an extension of) the heat
semigroup formula of \cite{Co4} (pp 563) for the Dixmier trace. We
then give a formula for the index of generalised Toeplitz
operators (Section 13), describe a special case of the Wodzicki
residue formula of Benameur and Fack \cite{BeF} for
pseudodifferential operators along the leaves of a foliation and
discuss in Section 15 a similar formula for pseudodifferential
operators with almost periodic coefficients. The formula of
Section 13 has, as a consequence, the index theorem of Lesch for
Toeplitz operators with noncommutative symbol and this is
described in Section 16. An extension of Theorem IV.2.8 of
\cite{Co4} on a residue formula for the Hochschild class of the
Chern character of Fredholm modules is given in Section 17 (see
also \cite{BeF}).

Finally, following \cite{AzSu}, we show how Lidskii's formula may
be extended to Dixmier traces in the von Neumann setting in
Section 18. In a series of corollaries we explain its relevance to
the question of measurability. In \cite{Co4}, except for the case
of pseudodifferential operators, measurability results are proved
for positive operators. The approach of \cite{AzSu} allows one to
address the problem of removing the positivity assumption.

We present below a short list of symbols and terminology used in this
paper with the indication of the place where these symbols and
notations are introduced:
\begin{itemize}
\item Symmetric and rearrangement invariant (r.i.) functionals and
rearrangement invariant spaces $E(J)$
   (Section 2, Definition~\ref{SandRIf});
\item Marcinkiewicz spaces, $M(\psi)$, $\clL^{(1,\infty)}$,
$\clL^{(p,\infty)}$
   (Section~\ref{MFandSS});
\item Generalized singular value function $\mu_{(\cdot)}(x)$
   (Section~\ref{PrelimDilInvState});
\item Semifinite von Neumann algebra $\mathcal N$, faithful normal
semifinite trace on $\mathcal N$, $\tau$. $\tau$-measurable
operators $\widetilde{\mathcal N}$, the fully symmetric operator
space associated to $({\mathcal N},\tau)$ and Banach function
space $E$ is denoted $E({\mathcal N},\tau) $ (Section 2.3); \item
Operator Marcinkiewicz spaces $M(\psi)\nt$ (Section 2.3); \item
Symmetric functionals on a fully symmetric space $E=E(0,\infty)$:
$E^*_{sym}$ (Section 3.1);  \item Sets of states $BL(\mbR)$,
$D(\mbR_+^*)$, $BL(\mbR_+)$
   (Section~\ref{PrelimDilInvState});
\item Dixmier traces $\tau_\omega$
(Definitions~\ref{dixtrd},\ref{saaodtd}),  traces $F_{\mbL},
F_\LL$ (Section 5.1); \item Connes-Dixmier traces $\tau_\omega$
with $\omega\in CD(\mbR_+^*)$ (Section~\ref{Connes-Dixmier
traces}); \item Measurable operators (Definitions~\ref{mesopd1}
   and~\ref{mesopd2});
\item Spectral triples $(\clA,\clH,\clD)$, grading operator $\Gamma$,
$\clI$-summability
   (Definition~\ref{SpectTripIsumd});
\item $QC^k$ semifinite spectral triples (Definition~\ref{qck});
\item Spectral dimension of $(\clA,\clH,\clD)$
   (Definition~\ref{dimension});
\item Spectral flow $sf(\{F_t\})$, $sf(D,uDu^*)$
   (Section~\ref{spectralflow});
\item $(b,B)$-cochain, $(b,B)$-cocycle, $(b^T, B^T)$-chain, $(b^T,
   B^T)$-cycle (Section~\ref{spectralflow});
   \item Zeta functions $\zeta(s)$, $\zeta_A(s)$ (Section 10.2);
\item Bohr compactification $\R^n_B$, $C^*$-algebra ${\mathcal AP}(\R^n)$
(Section~\ref{ShFr}).
\end{itemize}

\noindent{\bf Acknowledgement} The authors
would like to thank N.A. Azamov and A.A.
Sedaev, for comments, assistance and criticism.

\let\le\leq
\let\ge\geq
\def\x{x}
\def\y{y}
\def\z{z}


{

\setlength{\belowdisplayshortskip}{1\belowdisplayskip}
\newcommand{\no}{\hspace*{\indentation}}
\renewcommand{\arraystretch}{1.75}
\setlength{\arrayrulewidth}{0.5\arrayrulewidth}

\newcommand{\tab}[1]{\hspace*{#1\tablength}}
\newcommand{\nm}[1]{\mbox{\ensuremath{\| #1 \|}}}
\newcommand{\implies}{\ensuremath{\Rightarrow}}
\newcommand{\fa}{\ensuremath{\ \, \forall \,}}
\renewcommand{\iff}{\ensuremath{\; \Longleftrightarrow \;}}
\renewcommand{\sp}[1]{\ensuremath{\mathrm{sp}(#1)}}
\newcommand{\gap}[1]{\ensuremath{\, #1 \,}}
\newcommand{\inset}[2]{\ensuremath{\{ #1 \, | \, #2 \} }}
\newcommand{\inprod}[2]{\ensuremath{\langle #1 , #2 \rangle}}
\newcommand{\anti}[2]{\ensuremath{ \{ #1 , #2 \} }}
\newcommand{\CC}{\ensuremath{\mathbb{C}}}
\newcommand{\RR}{\ensuremath{\mathbb{R}}}
\newcommand{\TT}{\ensuremath{\mathbb{T}}}
\newcommand{\SB}{\ensuremath{\mathbb{S}}}
\newcommand{\ZZ}{\ensuremath{\mathbb{Z}}}
\newcommand{\HH}{\ensuremath{\mathbb{H}}}
\newcommand{\NN}{\ensuremath{\mathbb{N}}}
\newcommand{\iny}{\ensuremath{\infty}}
\newcommand{\bd}[1]{ \ensuremath{r_{\NN}(\ensuremath{#1})}}
\newcommand{\bp}{\ensuremath{p}}
\newcommand{\nmt}{\widetilde {\mathcal {M}}}
\newcommand{\nmm}{\mathcal {M}}
\newcommand{\emt}{E(\mathcal {M},\tau )}

\newcommand{\noline}{\vspace*{1\parskip}}
\newcommand{\preskip}{\vspace*{1\belowdisplayskip}}
\newcommand{\display}[1]{$$#1$$}
\newcommand{\smdisplay}[1]{\\[4pt] \no \centerline{#1} \\[4pt]}
\newcommand{\text}[1]{\mbox{\normalfont #1}}
\newcommand{\mod}[1]{\ensuremath{
         \text{ \hspace*{-3pt} mod \hspace*{-3pt} } #1}}

\section{Preliminaries: spaces and functionals}
\label{PrelimOnSpAndF}
\medskip \noindent Consider a Banach space $(E, \Vert \cdot \Vert_E)$ of
real valued Lebesgue measurable functions (with identification
$\lambda$ a.e.) on the interval $J=[0,\infty)$ or else on
$J=\mbN$. Let $x^*$ denote the non-increasing, right-continuous
rearrangement of $|x|$ given by \display{x^*(t) = \inf \inset{s
\geq 0}{\lambda(\{ |x|> s\}) \leq t} , \ t > 0,} where $\lambda$
denotes Lebesgue measure. Then $E$ will be called rearrangement
invariant if

\smallskip\noindent (i).\quad $E$ is an ideal lattice, that is if
$y\in E$, and $x$ is any measurable function on $J$ with $0\leq
|x|\leq |y|$, then $x\in E$ and $\Vert x\Vert_E\leq \Vert
y\Vert_E$;

\smallskip\noindent (ii).\quad if $y\in E$ and
  if $x$ is any
measurable function on $J$ with $x^*=y^*$, then $x\in E$ and
$\Vert x\Vert_E= \Vert y\Vert_E$.

\medskip \noindent In the case $J=\mbN$,
it is convenient to identify $x^*$ with the rearrangement of the
sequence $|x|=\{|x_n|\}_{n=1}^\infty$ in the descending order. For
basic properties of rearrangement invariant (=r.i.) spaces we
refer to the monographs \cite{KPS}, \cite{LT}, \cite{LT2}. We note
that for any r.i. space $E=E(J)$ the  following continuous
embeddings hold
$$
L_1\cap L_\infty(J)\subseteq E\subseteq L_1 + L_\infty(J).
$$

\medskip \noindent
The r.i. space $E$ is said to be fully symmetric Banach space if
it has the additional property that if $y\in E$ and $L_1 +
L_\infty(J)\ni x\prec\prec y$, then $x\in E$ and $\Vert
x\Vert_E\leq \Vert y\Vert_E$. Here, $x\prec\prec y$ denotes
submajorization in the sense of Hardy-Littlewood-P\'olya:
$$\int_0^tx^*(s)ds\leq  \int_0^ty^*(s)ds,\quad \fa t>0.
$$

\medskip \noindent
A classical example of non-separable fully symmetric function and
sequence spaces $E(J)$ is given by Marcinkiewicz spaces.

\subsection{Marcinkiewicz function and sequence spaces}\label{MFandSS}

\medskip \noindent
Let $\Omega$ denote the set of concave functions $\psi : [0,\iny)
\to [0,\iny)$ such that $\lim_{t \to 0^+} \psi(t) = 0$ and
$\lim_{t \to \iny} \psi(t) = \iny$. Important functions belonging
to $\Omega$ include $t$, $\log(1+t)$, $t^\alpha$ and
$(\log(1+t))^{\alpha}$ for $0 < \alpha < 1$.
Let $\psi \in \Omega$. Define the weighted
mean function
$$
a\left( x,t\right) =\frac{1}{\psi \left(
t\right) } \int_{0}^{t}x^{\ast }\left( s\right) ds\quad t>0
$$
  and denote by
$M(\psi)$ the Marcinkiewicz space of measurable functions $x$ on
$[0,\iny)$ such that
\begin{equation}\label{WeieghedMeanF}
  \nm{x}_{M(\psi)} := \sup_{t
 > 0} a\left( x,t\right) = \nm{a\left( x,\cdot\right)}_\iny
  < \iny.
\end{equation}
  The definition of the Marcinkiewicz sequence space
  $(m(\psi),\nm{x}_{m(\psi)})$ is similar,
$$
m(\psi)=\left\{x=\{x_n\}_{n=1}^\infty\ :\
\|x\|_{m(\psi)}:=\sup_{N\ge1}\frac{1}{\psi(N)}
\sum_{n=1}^Nx_n^*<\infty\right\}.
$$

\smallskip
\noindent The norm closure of $M(\psi) \cap L^1([0,\iny))$
(respectively, of $\ell_1=\ell_1(\mbN)$) in $M(\psi)$
(respectively, in $m(\psi)$) is denoted by $M_1(\psi)$
(respectively, $m_1(\psi)$). For every $\psi\in \Omega$, we have
$M_1(\psi)\neq M(\psi)$. The Banach spaces
$(M(\psi),\nm{.}_{M(\psi)})$, $(m(\psi),\nm{.}_{m(\psi)})$,
$(M_1(\psi),\nm{.}_{M(\psi)})$, $(m_1(\psi),\nm{.}_{m(\psi)})$ are
  examples of fully symmetric spaces
\cite{LT}, \cite{KPS}.

\indent Let $M_+(\psi)$ (respectively, $m_+(\psi)$) denote the set
of positive functions of $M(\psi)$ (respectively, $m(\psi)$). For
every $x\in M(\psi)$, we write $x=x_+-x_-$, where
$x_+:=x\chi_{\{t:x(t)>0\}}$ and $x_-:=x-x_+$.
The spaces
$$
\clL^{(1,\infty)}:=M(\log(1+t)) \cap L_\infty\quad
\text{and}\quad \clL^{(p,\infty)}:=M(t^{1-\frac 1p}) \cap L_\infty, \quad
1<p<\infty
$$
  play very important part in the sequel. Note that these spaces are
  still Marcinkiewicz spaces. Indeed, $\clL^{(1,\infty)}$
  (respectively, $\clL^{(p,\infty)}$, $p>1$) may be
  identified with the space $M(\psi_1)$
(respectively $M(\psi_p)$, $p>1$), where
$$
\psi_1(t)=\begin{cases}
t\cdot\log2,&0\le t\le 1\\
\log(1+t),& 1\le t <\infty
\end{cases},
$$
respectively,
$$
\psi_p(t)=\begin{cases}
t,&0\le t \le 1\\
t^{1-\frac{1}{p}},& 1\le t <\infty
\end{cases}.
$$
The (Marcinkiewicz) norm given by  formula~\eqref{WeieghedMeanF} on
the space $\clL^{(p,\infty)}$ is denoted by $\|\cdot\|_{(p,\infty)}$,
$1\le p <\infty$.
\subsection{Singular symmetric  functionals on Marcinkiewicz spaces.}
\label{SSFonMS}

\begin{defn}[cf. \cite{DPSS}, Definition 2.1]\label{SandRIf}
A positive functional $f\in M(\psi)^*$ is said to be symmetric
(respectively, r.i.) if $f(x) \leq f(y)$ for all $x,y \in
M_+(\psi)$ such that $x\prec\prec y$ (respectively, $x^*=y^*$).
Such a functional is said to be supported at infinity (or
singular) if $f(|x|) = 0$ for all $x \in M_1(\psi)$ (equivalently,
$f(x^*\chi_{[0,s]})=0$, for every $x\in M(\psi)$ and the indicator
function $\chi_{[0,s]}$ of the interval $[0,s]$ for all $s>0$).
\end{defn}
\smallskip
\indent
The following theorem completely characterizes Marcinkiewicz spaces
admitting non-trivial symmetric functionals.
\begin{thm}[{\cite{DPSS,DPSSS1,DPSSS2}}]\label{MSCwithNTF} If $\psi\in
\Omega$, then
\begin{thlist}
\item a non-zero symmetric functional on $M(\psi)$ (respectively,
   $m(\psi)$) supported at infinity exists if and only if
\begin{equation} \liminf_{t\to \infty} \frac{\psi (2t)}{\psi (t)}=1,
\label{eq1.1}
\end{equation}
\item a non-zero symmetric functional on $M(\psi)$ supported at zero
   exists if and only if
\begin{equation}
  \liminf_{t\downarrow 0} \frac{\psi (2t)}{\psi (t)}=1. \label{eq1.1s}
\end{equation}
\end{thlist}
\end{thm}
  Thus, for example,
$(\clL^{(1,\infty)})^*_{sym,\infty}\ne\{0\}$, whereas
$(\clL^{(p,\infty)})^*_{sym,\infty}=\{0\}$, for all $1<p<\infty$.
The conditions~\eqref{eq1.1} and~\eqref{eq1.1s} admit the following
geometric interpretation. Let us denote by $N(\psi)$ the norm closure
in $M(\psi) $ of the (order) ideal
$$
N(\psi) ^0:=\{f\in M(\psi) :f^*(\cdot) \leq k\psi
   ^\prime (\frac{\cdot}{k})\ {\rm for\ some } \ k\in {\mbN}\}.
$$
Clearly, $N(\psi)$ is a Banach function space (a subspace of
$M(\psi)$) and is rearrangement invariant. Assuming (for
simplicity) that $\psi$ is linear in a neighbourhood of $0$, the
space $N(\psi)$ is fully symmetric (and thus coincides with
$M(\psi)$) if and only if~\eqref{eq1.1} fails. In other words,
$M(\psi)^*_{sym,\infty}\ne\{0\}$ if and only if $N(\psi)\ne
M(\psi)$. For this and other geometric interpretations of
conditions~\eqref{eq1.1} and~\eqref{eq1.1s} we refer the reader
to~\cite{BM},\cite[II.5.7]{KPS}  and \cite{Ru}. For various
constructions of singular symmetric functionals on $M(\psi)$ (and
more generally on fully symmetric spaces and their non-commutative
counterparts) we refer to \cite{DPSS}, \cite{DPSSS1},
\cite{DPSSS2}. Constructions relevant to our main topic will be
reviewed below, in Section~\ref{concConstOfSSF}.

Our focus on symmetric functionals supported at infinity is explained by
the numerous applications of their non-commutative counterparts in
non-commutative geometry. Non-commutative
analogues of symmetric functionals supported at zero can be thought of
as ``Dixmier traces associated with von Neumann algebras of type
$II_1$'' and have not
found any applications to date.
\subsection{Symmetric operator spaces and functionals.}
Here, we extend the ideas of the previous sections to the
setting of (noncommutative) spaces of measurable operators. We
denote by $\clN $ a semifinite von Neumann algebra on the
Hilbert space $\clH $, with a fixed faithful and normal
semifinite trace $\tau $.
We shall be mainly concerned with $\tau(1)=\infty$, where $1$ is the
identity in $\clN$.
  A linear operator $x$:dom$(x)\to \clH  $, with domain dom$(x)\subseteq
\clH $, is called
  affiliated with $\clN $ if $ux=xu$ for all unitary $u$ in the
commutant $\clN '$ of $\clN $. The closed and densely defined
operator $x$, affiliated with  $\clN $, is called
$\tau$-measurable if for every $\epsilon >0$ there exists an
orthogonal projection $p\in \clN$ such the $p({\clH
})\subseteq$dom$(x)$ and $\tau(1-p)<\epsilon$. The set of all
$\tau$-measurable operators is denoted $\clNt$.

\medskip\indent We next recall the notion of generalized singular
value function~\cite{FK,F}. Given a self-adjoint operator $x$ in ${\mathcal
H}$, we denote by $e^x(\cdot )$ the spectral measure of $x$.  Now
assume that $x$ is $\tau$- measurable.  Then $e^{\vert
x\vert}(B)\in \clN $ for all Borel sets $B\subseteq{\mbR}$, and there
exists $s>0$ such that $\tau (e^{\vert
x\vert}(s,\infty))<\infty$. For $t\geq 0$, we define
$$\mu_t(x)=\inf\{s\geq 0 : \tau (e^{\vert x\vert}(s,\infty))\leq t\}.$$
The function $\mu(x):[0,\infty)\to [0,\infty ]$ is called the {\it
generalized singular value function} (or decreasing rearrangement)
of $x$; note that $\mu_t(x)<\infty$ for all $t>0$.  For the basic
properties of this singular value function we refer the reader to
\cite{FK}.

\medskip\indent If we consider $\clN  =L_\infty([0,\infty),m)$, where $m$
denotes Lebesgue measure on $[0,\infty)$, as an abelian von
Neumann algebra acting via multiplication on the Hilbert space
${\clH }=L^2([0,\infty),m)$, with the trace given by
integration with respect to $m$, it is easy to see that the set of
all $\tau$-measurable operators affiliated with $\clN $
consists of all measurable functions on $[0,\infty)$ which are
bounded except on a set of finite measure, and that the
generalized singular value function $\mu(f)$ is precisely the
decreasing rearrangement $f^*$.

\medskip\indent If $\clN =\clL(\clH)$ (respectively, $\ell_\infty(\mbN)$)
and $\tau $ is the standard trace $Tr$ (respectively, the counting
measure on $\mbN$), then it is not difficult to see that $\clNt
=\clN$. In this case, $x\in \clN$ is compact if and only if $\lim
_{t\to \infty }\mu _t (x)=0$; moreover,
$$
  \mu _n(x)=\mu _t(x), \quad t\in [n,n+1),\quad  n=0,1,2,\dots ,
$$
and the sequence $\{\mu _n(x)\}_{_{n=0}}^{\infty }$ is just the
sequence of eigenvalues of $\vert x\vert $ in non-increasing order
  and counted according to multiplicity.

\medskip \indent Given a semifinite von Neumann algebra
$(\clN,\tau)$ and a fully symmetric Banach function space
$(E,\Vert\cdot\Vert_E)$ on $([0,\infty),m)$, we define the
corresponding non-commutative space $\ent $ by setting
$$E(\clN,\tau) = \{ x\in\clNt : \mu (x)\in E\}.$$
Equipped with the norm $\Vert x\Vert _{_{E(\clN,\tau)}}:= \Vert
\mu (x)\Vert_E$, the space $(E(\clN,\tau),\Vert\cdot\Vert
_{_{E(\clN,\tau)}})$ is a Banach space and is called the
(non-commutative) fully symmetric operator space associated with
$(\clN,\tau)$ corresponding to $(E,\Vert\cdot\Vert_E)$. If
$\clN=\ell_\infty(\mbN)$, then the space $E(\clN,\tau)$ is simply
the (fully) symmetric sequence space $\ell_E$, which may be viewed
as the linear span in $E$ of the vectors $e_n=\chi_{_{[n-1,n)}}$,
$n\ge 1$ (see e.g. \cite{LT}).
In the case $\nt=(\clL(\clH),Tr)$, we denote $E\nt$ simply by
$E(\clH)$. Note, that the latter space coincides with the
(symmetrically-normed) ideal of compact operators on $\clH$
associated with (symmetric) sequence space $\ell_E$ (see
e.g~\cite{GK}).

We shall be mostly concerned with fully symmetric operator spaces
$E\nt$ and $E(\clH)$ when $E=M(\psi)$, in particular, when
$E=\clL^{(p,\infty)}$, $1\le p<\infty$. We refer to the spaces
$M(\psi)\nt$ as to operator Marcinkiewicz spaces. Sometimes, for
brevity, we shall omit the symbols $\nt$ and $\clH$ from the
notations and this should not cause a confusion.

Further references to the theory of fully symmetric operator spaces can
be found in \cite{DPSS, Sch90DAN, CS2}.

\begin{defn}
  A linear functional $\phi \in
E(\clN,\tau )^*$ is called symmetric (respectively,
r.i.) if $\phi $ is positive, (that is, $\phi
(x)\geq 0$ whenever $0\leq x\in E(\clN,\tau )$) and $\phi
(x)\leq \phi (y)$ whenever $\mu (x)\prec\prec \mu (y)$ (respectively,
$\phi (x)= \phi (y)$ whenever $x,y\ge 0$ and $\mu (x)= \mu (y)$).
\end{defn}

\medskip \indent For a given $x\in\clNt$, the set
$\Omega(x)=\{y\in\clNt\ :\ y\prec\prec x\}$ is called the orbit of
the operator $x$. If
$x\in L_1(\clN, \tau)+\clN$, then the set $\Omega(x)$ is
conveniently described in terms of absolute contractions.
Denote by $\Sigma$ the set of all linear operators $T: L_1\nt+\clN\to
L_1\nt +\clN$ such that $T(a)\in L_1\nt$ (respectively, $\clN$) if $a\in
L_1\nt$ (respectively, $\clN$) and such that $\|T\|_{L_1\nt\to
L_1\nt}\le1$, $\|T\|_{\clN\to\clN}\le 1$. It follows from~\cite{DDP2}
that $y\prec\prec x$, $x\in L_1\nt+\clN$, $y\in \clNt$ if and only
if there exists $T\in \Sigma$ such that $T(x)=y$. Thus,
$$
\Omega (x)=\{Tx\ :\ T\in \Sigma \}.
$$
If $E\nt$ is a fully symmetric operator space, we have $\Omega
(x)\subseteq E\nt$ for every $x\in E\nt$ and therefore a bounded
positive linear functional $\varphi$ on $E\nt$ is symmetric if and only
if $\varphi (|Tx|)\le \varphi (x)$ for every $T\in \Sigma $ and $0\le
x\in E\nt$.

\medskip \indent Now we assume that $\alpha :\clN\rightarrow
\clN$ is a $\ast $%
-automorphism which is in addition trace preserving, that is, $\tau
\left( \alpha \left( a\right) \right) =\tau \left( a\right) $ for all
$0\leq a\in \clN$. It is easy to see that such an automorphism
extends (uniquely) to a $\ast $-automorphism $
\tilde{\alpha}:{\clNt}\rightarrow\clNt$, which is rearrangement preserving,
that is, $\mu \left(
   \tilde{\alpha}\left( x\right) \right) =\mu \left( x\right) $ for all
$x\in \widetilde{\clN}$. Thus, we can view rearrangement
invariant functionals on $E(\clN,\tau )$ as positive
functionals which are invariant with respect to the action of the
group of all trace preserving $\ast$-automorphisms of $\clN$, which is a
subgroup of $\Sigma $. If $\clN
=\clL(\clH)$, then every rearrangement invariant
functional on $E(\clH)$ is simply a trace (i.e. unitary invariant positive
functional on $\clL(\clH))$, which extends to a
continuous linear functional on $E(\clH)$. Clearly, every
symmetric functional is rearrangement invariant. However,
{\em a priori},
it is not clear (see~\cite{T1}) whether there exists (for example on the ideal
$\clL^{(1,\infty)}(\clH)$) a rearrangement invariant singular functional,
which is not necessarily symmetric, or
whether there exists a trace on
$\clL^{(1,\infty)}(\clH)$, which does not coincide with a Dixmier trace (see
Section~\ref{concConstOfSSF} below). Very recently the first example
of such a trace has appeared in~\cite{KS}.  In fact, if
\begin{equation}\label{eq4nn}
\lim_{t\to\infty}\frac{\psi(2t)}{\psi(t)}=1,
\end{equation}
  then there exists a
non-zero trace (or r.i. functional) on every operator
Marcinkiewicz space $M(\psi)\nt$ which vanishes on $N(\psi)\nt$.
In particular, if $\psi(t)=\log(1+t)$, then such a functional
vanishes on every operator $0\le x\in\clL^{(1,\infty)}\nt$ with
$\mu_t(x)=\frac{1}{t}$ (or, if $\clN=B(\clH)$ then on every
compact operator $0\le x\in \clL^{(1,\infty)}(\clH)$, such that
$\mu_n(x)=\frac{1}{n}$, $n\ge1$). It is not clear yet whether r.i.
functionals which are not symmetric exist on every $M(\psi)\nt$
with $\psi$ satisfying condition~\eqref{eq1.1}.
}

\section{General facts about symmetric functionals.}\label{GenFactsOnSF}
Let
$E=E(0,\infty)$ be a fully symmetric space. By $E_+$ we denote the
set of all  nonnegative functions from $E$ and by $E^*_{sym}$ the
set of all symmetric functionals of $E$.

A positive linear functional $\phi$ on $E$ is called \newn{normal} (or
\newn{order continuous}),
     if from $f_n\downarrow 0$ it follows that $\phi(f_n)\downarrow 0.$
\begin{prop} \cite{DPSS}
    If a functional $\phi\in E^*_{sym}$ is order continuous, then $E\subset
    L^1[0,\infty),$ and
$\phi$ is proportional to the integral against the
    Lebesgue measure.
\end{prop}
We define the dilation operator $D_s$ as in \cite{KPS} by $D_s
f(t) = f(t/s)$. Note that $D_s$ is a bounded operator on $E$ and
$\norm{D_s}_{E\to E} \leq \max\set{1,s}$, moreover $(D_s f)^* =
D_s f^*$ for any function $f\in E.$ The following result is
established in~\cite[Proposition~2.3]{DPSS} under the assumption
that $\phi$ is symmetric, however the proof holds also for r.i.
functionals.
\begin{prop}\label{gfasfp1}
If $0\le\phi\in E^*$ is rearrangement invariant, then
$\phi(D_sf)=s\phi(f)$ for all $f\in E$ and $s>0$.
\end{prop}

A positive  $\phi\in E^*$ is said to be \newn{singular}, if from
$0\leq \phi' \leq \phi$, $\phi'\in E^*,$ $\phi'$ is
order-continuous, it follows that $\phi' = 0.$

\begin{prop} \cite{DPSS} (i)
Every symmetric functional on $E$ can be uniquely decomposed into the sum of
a normal functional and a singular symmetric functional.
Moreover, the  normal functional is zero unless  $E\subseteq L^1[0,\infty).$\\
(ii) Any singular symmetric functional can be uniquely decomposed into the
sum of singular symmetric functionals, supported at zero and at infinity.\\
(iii) The set of symmetric functionals forms a lattice.
\end{prop}
The following result shows that every symmetric functional on $E$ admits
a ``natural extension'' up to a symmetric functional on $\ent$ for every
semifinite von Neumann algebra $\nt$. By $E\nt_{+}$ we denote the set
of all positive operators from $E\nt$.
{
\def\p{\,\prec\kern-1.2ex \prec\,}
\def\ch{\raise 0.5ex \hbox{$\chi$}}
\def\m{\ ${\cal M}$}
\def\n{\ ${\cal N}$}
\def\nm{{\cal M}}
\def\mt{\ $\widetilde {\cal M}$}
\def\nmt{\widetilde {\cal M}}
\def\nn{{\cal N}}
\def\nnt{\widetilde {\cal N}}
\def\r{{\mbR}^+}
\def\h{H({\cal M})}
\def\g{G({\cal M})}
\def\up{M^{(p)}(E)}
\def\do{M_{(p)}(E)}
\def\doq{M_{(q)}(E)}
\def\em{E({\cal M},\tau )}
\def\mpc{$(L^p(\nm),{\cal M})$}
\def\mpcn{(L^p(\nm),{\cal M})}
\def\mp+{$L^p(\nm)+{\cal M}$}
\def\mp+n{L^p(\nm)+{\cal M}}
\def \cf {\rlap {$ \prec $}{$ \relbar$}}
\def \cfn {\rlap {\prec }{\relbar }}
\def \ycfx {$y {\rlap {\prec }{\relbar}} x$}
\def \ycfpx {$y {\rlap {\prec }{\relbar}} ^px$}

\begin{thm}[\cite{DPSS}]\label{gfasft7} Let $\phi_0\in E_{sym}^* $.
   If $\phi (x):=\phi _0(\mu (x))$, for all $x\in E\nt_{+}$, then
   $\phi $ extends to a symmetric functional $0\leq \phi \in
   E(\clN,\tau )^*$.
\end{thm}
\begin{proof} It clearly suffices to show that $\phi $ is additive on
   $E\nt_+$.  Let $x,y\in E\nt_+$.  Since $\mu (x+y)\p \mu
   (x)+\mu (y)$ (\cite[Theorem 4.4]{FK}), and since $\phi _0$ is
   symmetric, it follows that $$\phi (x+y)=\phi _0(\mu (x+y))\leq \phi
   _0(\mu (x)+\mu (y))=\phi (x)+\phi (y).$$ To prove the converse
   inequality, we use the easily verified fact (see
   e.g.~\cite[Proposition 1.10]{GI1}) that
\begin{equation}\label{eqdpss4.1} \int _0^t \mu_s(x)ds +\int _0^t\mu
_s(y)ds \leq \int _0^{2t}\mu
_s(x+y)ds,\quad \forall t>0.
\end{equation}
Observing that~\eqref{eqdpss4.1} is equivalent to the submajorization
\begin{equation}\label{eqdpss4.2}
\mu (x)+\mu (y)\p 2D_{\frac{1}{2}}\mu (x+y),
\end{equation}
it follows from Proposition~\ref{gfasfp1} that
\begin{equation*}
\begin{split}
\phi (x)+\phi (y)&=\phi _0(\mu (x)+\mu (y))\\
&\leq 2\phi _0(D_{\frac{1}{2}}\mu (x+y))=\phi _0(\mu (x+y))=\phi
(x+y).\\
\end{split}
\end{equation*}
Thus $\phi $ is additive on $E\nt_+$ and this suffices to
complete
the proof of the Theorem.
\end{proof}

\begin{thm}\label{gfasft8} Let $(\clN,\tau )$ be a
semifinite
von Neumann algebra without minimal projections, and let $E$ be a fully
symmetric Banach function space on $[0,\infty )$. If $0 \leq \phi \in
E(\clN,\tau )^*$ is a symmetric functional,
then there exists a symmetric
functional $0\leq \phi_0\in E^*$ such that $\phi (x)=\phi _0(\mu (x))$
for all
$0\leq x\in E(\clN,\tau )$.
\end{thm}
\begin{proof}
It  is sufficient to show that there exists a symmetric functional
$0\leq \phi
_0\in E[0,\tau (1))^*$ satisfying $\phi (x)=\phi _0(\mu (x))$ for all
$0\leq x \in
E(\clN,\tau )$.  Let $\clM$ be the commutative von Neumann
algebra
$L_\infty [0,\tau (1))$, with trace given by integration.  The algebra
$\clMt$ may be identified with the space of all
measurable functions
on $[0,\tau (1))$ which are bounded except on a set of finite measure.
Since
$\clM$ does not contain any minimal projections, there exists a
positive
rearrangement-preserving algebra homomorphism $J:\clMt\to
\clNt$ (\cite[Lemma~4.1]{CS2}, \cite[Theorem~3.5]{DDP1}).  Let
$0\leq \phi\in E\nt^*$ be
symmetric. For $f\in E[0,\tau (1))$, define $\phi _0(f):=\phi
(Jf)$.  It is clear that $0\leq \phi_0 \in E[0,\tau (1))^*$ is
symmetric.
Moreover,
if $0\leq x\in E\nt$, then $\mu (J(\mu(x)))=\mu (x)$ and
hence
$\phi (x)=\phi (J(\mu (x)))=\phi _0(\mu (x))$.
\end{proof}

It is even easier to see that the correspondence between the sets
$E^*_{sym}$ given in Theorems~\ref{gfasft7} and~\ref{gfasft8} also
exists for the set of symmetric functionals $(\ell_E)^*_{sym}$ and
$(E(\clH))^*_{sym}$. Furthermore, as the following result shows there exists
a simple connection also between the sets $(\ell_E)^*_{sym}$ and $E^*_{sym}$.

\begin{thm}[\cite{DPSS}] Let $E$ be a fully symmetric Banach function space on
   $[0,\infty )$ and let $E(\clH)$ be the corresponding ideal of compact
   operators on infinite-dimensional Hilbert space $\clH$.  If $0\leq \phi
\in E\nt^*$ is a symmetric
   functional, then there exists $\phi _0\in E^*_{sym}$ such that $\phi
   (x)=\phi _0(\mu (x))$ for all $x\in E\nt_+$.
\end{thm}

}

Now, let us consider the following question naturally arising from
Theorem~\ref{gfasft7} and suggested in~\cite{T1}. Suppose that $E$
(respectively, $\ell_E$) is a
fully symmetric function (respectively, sequence) space and $\varphi_0$
is a positive r.i. functional on $E$ (respectively,
$\ell_E$). Is the functional $$\varphi (x):=\varphi _0(\mu(x))$$ additive
on $E\nt$ (respectively, $E(\clH)$)? Very recently this question has
been answered in the affirmative in~\cite{KS}, provided that $\clN$ is a
factor. The proof given in~\cite{KS} is based on the deep results
from~\cite{DFWW} (see also~\cite{K1,DK2}) concerning the structure of
the commutator space~ $[I,J]$ spanned by all commutators $[T,S]$, $T\in
I$, $S\in J$ (here, $I$ and $J$ are ideals of compact operators from
$\clL(\clH)$).
Note that every r.i. functional on the ideal $I$  vanishes on
$[I,\clL(\clH)]$. Many important results from~\cite{DFWW,K1,DK2}
admit an extension to the case of general semifinite von Neumann
algebras and their ideals~\cite{DK,Fa04FA}
\section{Preliminaries on dilation and translation invariant states.}
\label{PrelimDilInvState}
A construction of Dixmier traces $\tau_\omega$ depends crucially
on the choice of the ``invariant mean'' $\omega$. In this section
we recall and review the most important classes of such means. We
denote by $\ell_\infty=\ell_\infty(\mbN)$ the Banach space of all
bounded sequences of complex numbers. By a state on a unital
$C^*$-algebra we mean a positive linear functional with value 1 on
the unit of the algebra. We recall that a positive linear
functional ${{\LL}}\in \ell_{\infty }^{\ast }$ is called a Banach
limit if ${\LL}$ is translation invariant and ${{\LL}}\left(
\mathbf{1}\right) =1$ (here, $\mathbf{1}=(1,1,1\dots)$). A Banach
limit ${\LL}$ satisfies in particular ${{\LL}}\left( \xi \right)
=0$ for all $\xi \in c_{0}$ ($=$ all sequences from $\ell_\infty$
converging to zero). We denote the collection of all Banach limits
on $\ell_{\infty }$ by $BL(\mbN) $. Note that $\left\| \LL\right\|
=1$ for all $ \LL\in BL(\mbN) $.

\bigskip
We recall that sequence $\xi =\left\{ \xi _{n}\right\}
_{n=1}^{\infty}\in \ell_{\infty }$ is said to be {\it almost
convergent}  to $\alpha\in \mbR$, denoted $\displaystyle
F\text{-}\lim_{n\rightarrow \infty }\xi _{n}=\alpha $ if and only
if ${{\LL}}\left( \xi \right) =\alpha$ for all ${{\LL}}\in BL(\mbN)$.
The notion of an almost convergent sequence is due to
G.G. Lorentz~\cite{Lor},
  who showed that the sequence $\left\{ \xi _{n}\right\}
_{n=1}^{\infty}$ is almost convergent to $\alpha$ if and only if
the equality $\displaystyle
\lim_{p\rightarrow \infty }\frac{\xi _{n}+\xi _{n+1}+\cdots +\xi _{n+p-1}}{p}%
=\alpha $ holds uniformly for $n=1,2,\dots\ $. We denote by $ac$
(respectively, $ac_0$) the set of all almost convergent
(respectively, all almost convergent to $0$) sequences from
$\ell_\infty$. Clearly, $ac$ and $ac_0$ are closed subspaces in
$\ell_\infty$.
We define the shift operator $\Shift\colon \ell_\infty\to \ell_\infty,$
the  Ces\`{a}ro operator $\Cesaro:\ell_\infty\to \ell_\infty$
and dilation operators $\Dilation_n:\ell_\infty \to \ell_\infty$ for $n\in
\mathbb{N}$ by formulas
\begin{gather*} 
   \Shift\Seque\x =(\x_2,\x_3,\x_4,\dots). \\
   \Cesaro\Seque\x = (x_1,\frac {x_1+x_2}2,\frac {x_1+x_2+x_3}3,\dots), \\
   \Dilation_n\Seque\x=(\underbrace{\x_1,\ldots,\x_1}_n,\underbrace{\x_2,\ldots,\x_2}_n,\ldots),
\end{gather*}
for all $\x=\Seque \x \in \ell_\infty.$

Each of the above operators is positive and leaves invariant the unit
element $\mathbf{1}$ of the algebra $\ell_\infty$ and consequently is
bounded with the norm equal to $1$.  Moreover, $\sequen \Dilation n$
is an abelian semigroup. The main tool in our construction of various
classes of invariant means is the well known Markov-Kakutani fixed point
theorem.

\begin{thm}[Markov-Kakutani] Let $F$ be a locally convex Hausdorff space
and let $K$ be a non-empty convex compact subset of $F$. Let $\clT$ be an
abelian semigroup of linear continuous operators on $F$ such that
$S(K) \subseteq K$ for all $S\in \clT.$ Then there exists $\x\in K$ such
that $Sx=\x$ for all $S\in\clT.$
\label{MarkovKakutani}
\end{thm}

It is easy to see that the set of all fixed points from Theorem
\ref{MarkovKakutani} forms a convex compact subset of $K.$

We shall be applying the Markov-Kakutani theorem in the setting when
$F=(\ell_\infty)^*$, $(L_\infty(\mbR))^*$, $(L_\infty(\mbR_+^*))^*$ equipped
with the weak $*$-topology. For simplicity of exposition we present
the proofs only for the first case.
\begin{lemma}[\cite{DPSSS2,CPS2}] \label{DilCesShiftLemma} The following is
true. \\
   $(i)$  $\Dilation_n\Shift=\Shift^n\Dilation_n \quad  \forall \ n\geq 1$; \\
   $(ii)$  $\Cesaro\Shift\x-\Shift\Cesaro\x\in c_{0}  \quad  \forall \ \x
   \in \ell_\infty$; \\
   $(iii)$  $\Cesaro\Dilation_n\x-\Dilation_n\Cesaro \x\in c_{0} \quad
   \forall \ \x\in \ell_\infty.$
\end{lemma}
\begin{proof}
The proof of (i) is straightforward. For the proof of (ii) note that for
all $\x\in \ell_\infty$
$$
   |(\Cesaro\Shift\x)_k -( \Shift\Cesaro\x)_k| =
\left|\frac{1}{k+1}\frac{\x_2 +\cdots +\x_{k+1}}{k}-\frac{1}{k+1}\x_1\right|
    \leq \frac{2}{k+1}\|\x\|_\infty,
$$
which shows that $\Cesaro\Shift\x-\Shift\Cesaro\x\in c_0.$

We now indicate the proof of (iii). Let $n\geq 1$ and $\x\in \ell_\infty.$
For $1\leq k\in \mathbb{N}$ there exist $l\geq 1$ and $1\leq r\leq n$ such
that $k=(l-1)n+r$. Hence
$$
   (\Cesaro\Dilation_n\x)_k = \frac 1k \sums{i=1}^{k}(\Dilation_n\x)_i
   = \frac{n}{k}\sums{j=1}^{l-1}\x_j + \frac{r}{k}\x_l
$$
and
$$
   (\Dilation_n\Cesaro\x)_k = (\Cesaro x)_l =
   \frac{1}{l}\sums{j=1}^l\x_j.
$$
Noting that, $nl-k=n-r\leq n$ and $rl-k\leq nl-k \leq n,$ it follows that
$$
    \left|(\Cesaro\Dilation_n\x)_k - (\Dilation_n\Cesaro\x)_k \right|
    =\left|\frac{nl-k}{kl}\sums{j=1}^{l-1}\x_j +\frac{rl-k}{kl} \x_l\right|
   \leq \frac{n}{k}\left( \frac{1}{l}\sums{j=1}^{l-1}|\x_j|\right) + \frac
nk |\x_l| \leq  \frac{2n}{k}\|\x\|_\infty.
$$
This shows that $\Cesaro\Dilation_n\x-\Dilation_n\Cesaro\x \in
\ell^{(1,\infty)}\subseteq c_0.$
\end{proof}

\begin{thm} 
\label{InvState}
There exists a state $\tilde\omega$ on $\ell_\infty$ such that for all
$n\geq 1$
$$
   \tilde\omega\circ\Shift = \tilde\omega\circ\Cesaro =
\tilde\omega\circ\Dilation_n =\tilde\omega.
$$
\end{thm}
\begin{proof}
Let $K=\left\{ 0\leq \varphi \in (\ell_\infty)^* \colon \varphi
(\mathbf{1})=1, \ \Shift^*\varphi =\varphi \right\}.$ Since $K$
contains ordinary Banach limits it is not empty. It is clear that
$K$ is convex and $*$-weakly compact. We claim that
$\Dilation_n^*(K) \subseteq K$. Indeed, by Lemma
\ref{DilCesShiftLemma}(i) above we know that
$\Shift^*\Dilation_n^*=\Dilation_n^*(\Shift^*)^n$, hence for
$\varphi \in K$
$$
   \Shift^*( \Dilation_n^*\varphi) = \Dilation_n^*(\Shift^*)^n\varphi
   =\Dilation_n^*\varphi,
   $$
   which implies that $\Dilation_n^*\varphi \in K$. Therefore we may
   apply Theorem~\ref{MarkovKakutani} to the set
   $K$ and the abelian semigroup $\seque {\Dilation^*}.$ Consequently
   the set
$$
   K_1=\left\{ 0\leq \varphi \in (\ell_\infty)^\ast:
\varphi(\mathbf{1})=1,\ \Shift^*\varphi=\varphi,\
   \Dilation_n^*\varphi =\varphi, \ n\geq 1\right\}
$$
is non-empty and again it is clear that $K_1$ is convex and $*$-weakly compact.

We now show that $\Cesaro^*(K_1) \subseteq K_1$. To this end,
first observe that $\varphi(\z) =0$ for all $\z\in c_{0}$ and all
$\varphi \in K_1$ (as $\Shift^*\varphi =\varphi$). Given $ \varphi
\in K_1$ it follows from Lemma \ref{DilCesShiftLemma}(iii) that
$$
   (\Dilation_n^*\Cesaro^*\varphi)(\x)
   -(\Cesaro^*\Dilation_n^*\varphi)(\x) = \varphi
(\Cesaro\Dilation_nx-\Dilation_n\Cesaro x) =0
$$
for all $\x\in \ell_\infty$ and so
$$
\Dilation_n^*( \Cesaro^*\varphi ) =\Cesaro^*( \Dilation_n^*\varphi )
=\Cesaro^*\varphi.
$$
Similarly it follows from Lemma \ref{DilCesShiftLemma}(ii) that
$\Shift^*( \Cesaro^*\varphi
) =\Cesaro^*\varphi $ for all $\varphi \in K_1$. Consequently, $%
\Cesaro^*( K_1) \subseteq K_1$. Applying
Theorem~\ref{MarkovKakutani} to the set $K_1$ and the semigroup
$\left\{(\Cesaro^*)^n\right\}_{n=0}^\infty$, we may conclude that
there exists $\tilde\omega\in K_1$ such that
$\Cesaro^*(\tilde\omega) = \tilde\omega,$ by which the proof is
complete.
\end{proof}

We define the isomorphism $L: L_\infty(\mbR)\to L_\infty(\mbR^*_+)$
by $L(f) = f\circ \log$. Firstly, we define the Cesaro means (transforms)
on $L_\infty(\mbR)$ and $L_\infty(\mbR^*_+)$, respectively by:
$$H(f)(u)=\frac{1}{u}\int_0^u f(v) dv\quad \text{for}\quad f\in
L_\infty(\mbR),\ u\in \mbR$$
and,
$$M(g)(t)=\frac{1}{\log t}\int_1^t g(s) \frac{ds}{s}\quad \text{for}\quad g\in
L_\infty(\mbR^*_+),\ t>0.$$

A brief calculation yields for $g\in  L_\infty(\mbR^*_+)$,
$$LHL^{-1}(g)(r) = \frac{1}{\log r}\int_0^{\log r} g(e^u)du
= \frac{1}{\log r}\int_1^{r} g(v)\frac{dv}{v}= M(g)(r),$$
i.e  $L$ intertwines the two means.

We shall now consider analogues of the operators $\Shift, \Dilation_n$
and $\Cesaro$ acting on $L_\infty(\mbR)$ and $L_\infty(\mbR_+^*)$.

\begin{defn}\label{familiesOfSelf-mapsDef}
Let $T_b$ denote translation by $b\in \mbR$, $D_a$ denote
dilation by $\frac{1}{a}\in \mbR^*_+$ and let $P^a$ denote exponentiation
by $a\in \mbR^*_+$. That is,
   \begin{eqnarray*}
     T_b(f)(x) &=& f(x+b)\quad \text{for}\quad f\in L_\infty(\mbR),\\
     D_a(f)(x) &=& f\left(\frac{x}{a}\right)\quad \text{for}\quad f\in
L_\infty(\mbR),\\
     P^a(f)(x) &=& f(x^a)\quad \text{for}\quad f\in L_\infty(\mbR^*_+).
   \end{eqnarray*}
\end{defn}

Some of the basic relations between these $L_\infty$ spaces and their
self-maps are provided for easy access by the following proposition,
whose proof is similar to Lemma \ref{DilCesShiftLemma}.

\begin{prop}[\cite{CPS2}] \label{DualityProp}
$L_\infty(\mbR)$ together with the self-maps, $D_a$, $T_b$, and $H$
($a>0,b\in\mbR$) is related to $L_\infty(\mbR^*_+)$ together with the
self-maps, $P^a$, $D_a$, and $M$ ($a>0$) via the isomorphism
$$L: L_\infty(\mbR)\to L_\infty(\mbR^*_+)$$ and the following identities:\\
(1) $LD_{\frac{1}{a}}L^{-1} = P^a$ for $a>0$,\\
(2) $LT_bL^{-1} = D_{(\exp(b))^{-1}}$ for $b\in\mbR$,
\\
(3) $LHL^{-1} = M$,\\
(4) $D_aH = HD_a$ and $P^aM = MP^a$ for $a>0,$\\
(5) $\liminfty{t} (HT_b-T_bH)f(t)=0$ for $f\in L_\infty(\mbR)$ and
     $b\in\mbR,$\\
(6) $\liminfty{t} (MD_a-D_aM)f(t)=0$ for $f\in L_\infty(\mbR^*_+)$ and
     $a>0.$
\end{prop}

\begin{prop}[\cite{CPS2}]\label{cps2prop}
If a continuous functional $\tilde\omega$ on $L_\infty(\mbR)$ is invariant
under the Cesaro operator $\Cesaro,$
the shift operator $\Shift_a$ or the dilation operator $\Dilation_a$ then
$\tilde\omega\circ L^{-1}$
is a continuous functional on  $L_\infty(\mbR^*_+)$ invariant under
$\CesaroRPlus,$ the dilation operator $\Dilation_a$ or
$P^a$ respectively. Conversely, composition with $L$ converts an
$\CesaroRPlus,$ $\Dilation_a$ or $P^a$ invariant
continuous functional on $L_\infty(\mbR^*_+)$
into an $\Cesaro,$ $\Shift_a$ or $\Dilation_a$ invariant continuous
functional on  $L_\infty(\mbR)$.
\end{prop}

We denote by $C_0(\mbR)$ (respectively, $C_0(\mbR^*_+)$) the continuous
functions
on $\mbR$ (respectively, $\mbR^*_+$) vanishing at infinity.

The proof of the following theorem is similar to that of Theorem
\ref{InvState}.
\begin{thm}[\cite{CPS2}] \label{DualityThm}
There exists a state $\tilde\omega$ on $L_\infty(\mbR)$
satisfying the following conditions:\\
(1) $\tilde\omega(C_0 (\mbR)) \equiv 0$.\\
(2) If $f$ is real-valued in $L_\infty(\mbR)$
then
$$ess \liminf\limits_{t\to\infty} f(t) \leq\tilde\omega(f)\leq ess
\limsup\limits_{t\to\infty}
f(t).$$
(3) If the essential support of $f$ is compact then $\tilde\omega(f)=0.$\\
(4) For all $a>0$ and $c\in\mbR$
$$
    \tilde\omega = \tilde\omega\circ T_c = \tilde\omega\circ D_a =
\tilde\omega\circ H.
$$
\end{thm}
Combining Theorem~\ref{DualityThm} and Proposition~\ref{cps2prop}, we
obtain
\begin{cor} \label{DualityCor}
There exists a state $\omega$ on $L_\infty(\mbR^*_+)$
satisfying the following conditions:\\
(1) $\omega(C_0(\mbR^*_+)) \equiv 0$.\\
(2) If $f$ is real-valued in $L_\infty(\mbR^*_+)$
then
$$ess\liminf\limits_{t\to\infty} f(t) \leq\omega(f)\leq ess
\limsup\limits_{t\to\infty} f(t).$$
(3) If the essential support of $f$ is compact then $\omega(f)=0.$\\
(4) For all $a,c>0$
$$
    \omega = \omega\circ D_c = \omega\circ P^a = \omega\circ M.
$$
\end{cor}


The results given in Theorem~\ref{DualityThm} and
Corollary~\ref{DualityCor} allow one to exercise an alternative
approach to the theory of Dixmier~\cite{Dix66CR} and Connes-Dixmier
traces~\cite{Co4}.  Whereas Dixmier's original approach is based on
the use of dilation invariant functionals, we replace the latter with
Banach limits ($=$ translation invariant functionals)
and make use of the  well-developed theory of
almost convergent sequences.

We introduce the following notation
\begin{align*}
BL(\mbR) =\{&\text{the set of all states $\tilde \omega$ on
   $L_\infty(\mbR)$ satisfying conditions (1)--(3) of
   Theorem~\ref{DualityThm}}\\
&\text{such that $\tilde \omega \circ T_c= \tilde
   \omega $ for every $c\in\mbR$}\},\\
D(\mbR_+^*) =\{&\text{the set of all states $\omega$ on
   $L_\infty(\mbR_+^*)$ satisfying conditions (1)--(3) of
   Corollary~\ref{DualityCor}}\\
&\text{ such that $\omega \circ D_c=
   \omega $ for every $c\in\mbR_+^*$}\},\\
BL(\mbR_+) =\{&\text{the set of all states $\tilde \omega$ on
   $L_\infty(\mbR_+)$ satisfying conditions (1)--(3) of
   Corollary~\ref{DualityCor}}\\
&\text{such that $\tilde \omega \circ T_c= \tilde
   \omega $ for every $c\in\mbR_+$}\}.\\
\end{align*}

The following simple remark plays an important role in the sequel.
If $\omega\in D(\mbR_+^*)$, then $\mbL:=\omega \circ L$ belongs to
$BL(\mbR)$. If $\mbL\in BL(\mbR_+)$, then $ \tilde \omega:= \mbL
\circ L^{-1}$ belongs to $D(\mbR_+^*)$. Finally, note that the
isomorphism $L: L_\infty(\mbR)\to L_\infty(\mbR^*_+)$ sends the
space $C_b(\mbR)$ of all bounded continuous functions on $\mbR$
onto the space $C_b(\mbR_+)$ of all bounded continuous functions
on $\mbR_+$. Thus, one can reformulate all the results from
Propositions \ref{DualityProp}, \ref{cps2prop}, Theorem
\ref{DualityThm} and Corollary \ref{DualityCor} for the
spaces of continuous bounded functions on $\mbR$ and $\mbR_+$.

\section{Concrete constructions of singular symmetric
functionals.}\label{concConstOfSSF}
\subsection{Dixmier traces}
If $\omega $ is a state on $\ell_\infty$ (respectively, on
$L_\infty(\mbR)$, $L_\infty(\mbR_+^*)$), then we shall frequently denote
its value on the element $\{x_i\}_{i=1}^\infty$ (respectively, $f\in
L_\infty(\mbR)$, $L_\infty(\mbR_+^*)$) by $\omega\mbox{-}\lim_{i\to\infty} x_i$
(respectively, $\omega\mbox{-}\lim_{t\to\infty}f(t)$).
Recall that Theorems~\ref{InvState}, \ref{DualityThm} and
Corollary~\ref{DualityCor} guarantee the existence of translation and/or
dilation invariant states on $\ell_\infty$, $L_\infty(\mbR)$,
$L_\infty(\mbR_+^*)$.
For simplicity, we explain the construction of Dixmier traces for the
ideal of compact operators $(\clL^{(1,\infty)}(\clH),
\|\cdot\|_{(1,\infty)})$ defined in Section~\ref{MFandSS}.
\begin{defn}\label{dixtrd}
Let $\omega$ be a $D_2$-invariant state on $\ell_\infty$. \index{Dixmier trace}
\newn{Dixmier trace} of $T \in \clL_+^{(1,\infty)}(\clH)$ is a number
$$
\tau_\omega(T) := \limo{N}\ILogN \sums{n=1}^N \mu_n(T).$$
\end{defn}
\begin{rems*}
We have deliberately chosen $\omega$ to satisfy only the
   dilation invariance assumption in the proof
   below, even though Dixmier originally imposed on $\omega$ the
   assumption of dilation and translation invariance.
We shall discuss
   differences below.
\end{rems*}
\begin{prop}\label{TrisLinear} $\tau_\omega(S+T) = \tau_\omega(S) +
\tau_\omega(T)  \quad \forall  S, \ T
\in \DixIdealPos{1}$.
\end{prop}
\begin{proof}
Set, for brevity
$$
\sigma_N (X):=\sum_{n=1}^N \mu_n(X),\quad X\in \clL^{(1,\infty)}(\clH)
$$
and note that
$$
\sigma_N (X)=\sup\{Tr(XP)\ :\ P= P(\clH)\ \text{is a projection
and}\ \dim P(\clH)=N\}.
$$
(see e.g.~\cite{GK, FK}). For a given $\epsilon >0$, let projections
$P_1$ and $P_2$ satisfy the conditions $\dim P_1(\clH)=\dim P_2(\clH)=N$
and $Tr(SP_1)> \sigma_N (S)- \epsilon $, $Tr(TP_2) > \sigma_N (T)-
\epsilon $. Setting $P:=P_1\vee P_2$, we have
$$
Tr((S+T)P)= Tr(SP)+Tr(TP)\ge
Tr(SP_1)+Tr(TP_2)
 > \sigma_N (S)+\sigma_N (T)-2 \epsilon .
$$
Since $\dim P(\clH)\le 2N$ and $ \epsilon $ is an arbitrary positive
number, we have
$$
\sigma_{2N} (S+T)\ge \sigma_N (S) + \sigma_N (T),\quad N\ge1.
$$
Setting, for brevity,
  $$\alpha_N=\ILogN \sigma_N(S), \quad \beta_N=\ILogN \sigma_N(T), \quad
\gamma_N=\ILogN \sigma_N(S+T).$$
  we restate the above inequality as
  $$
\frac{ \log (2N+1) }{ \log (N+1) } \gamma_{2N} \ge \alpha _N + \beta _N,\quad
N\ge1.
$$
Assume, for a moment, that we know
\begin{equation}\label{eq1new}
\omega\text{-} \lim_{N\to\infty} \gamma_{2N} = \omega \text{-}
\lim_{N\to\infty}
\gamma _N.
\end{equation}
Noting that $\{ \gamma_{2N} \}_{N\ge1}\in \ell_\infty$ and so
$\frac{\log (1+2N) }{\log (1+N) }\gamma_{2N} -\gamma_{2N} \to0$, we infer from
the above inequality that
$$
\tau_\omega (S+T)= \omega \text{-} \lim_{N\to\infty} \gamma _ N
\ge \omega\mbox{-}\lim_{N\to\infty} \alpha _N + \omega\mbox{-}
\lim_{N\to\infty}
\beta _N
=\tau_\omega (S)+ \tau_\omega (T).
$$
Since the converse inequality $\tau_\omega (S+T)\le \tau_\omega (S)+
\tau_\omega (T)$ follows immediately from the well-known inequality
$\sigma_N (S+T)\le \sigma_N (S)+ \sigma_N (T)$ (see~\cite{GK, FK}), the
proof is completed. It remains to explain equality~\eqref{eq1new}.

Note that $D_2$-invariance of $\omega $ immediately implies that
$$ \omega (\{\gamma_{2N} \}_{N=1}^\infty)
=\omega (D_2\{ \gamma_{2N} \}_{N=1}^\infty)
=\omega (\{\gamma_2,\gamma_2,\gamma_4,\gamma_4,\gamma_6,\gamma_6,\ldots\})
$$
and therefore, in order to prove~\eqref{eq1new} it is sufficient to
verify that
$$
D_2\{\gamma_{2N}\}_{N=1}^\infty-\{\gamma_N\}_{N=1}^\infty=
\{\gamma_2,\gamma_2,\gamma_4,\gamma_4,\gamma_6,\gamma_6,\ldots\}
-\{\gamma_1,\gamma_2,\gamma_3,\gamma_4,\gamma_5,\ldots\}\in c_0,
$$
or equivalently, that $\gamma_{2N} - \gamma_{2N-1} \to0$, as
$N\to\infty$. For $N\ge2$, we have
$$
\gamma_{2N} -\gamma_{2N-1} =\left( \frac{1}{\log 2N}-\frac{1}{\log
(2N-1)}\right)\sigma_{2N-1} (T+S)+\frac{1}{\log 2N} \cdot\mu_{2N}(T+S).
$$
It is obvious that the second summand above tends to $0$ as
$N\to\infty$. Noting that the condition $T+S\in \clL^{(1,\infty)}(\clH)$
guarantees $\sigma_{2N-1} (T+S)=O(\log (2N -1))$ and that
$\frac{1}{\log 2N}-\frac{1}{\log (2N-1)}=o \left( \frac{1}{\log
(2N-1)}\right)$, we see that the first summand also tends to $0$ as
$N\to\infty$.
\end{proof}
\begin{rems}\label{IsometricEmbeddingRem}
   Consider an isometric embedding $i: \ell_\infty\to L_\infty[0,\infty)$
   given by $\{x_j\}_{j=1}^\infty\stackrel{i}{\mapsto} \sum_{j=1}^\infty x_j
   \chi_{[j-1,j)}$, where $\chi_{[j-1,j)}$ is the characteristic
   function of the interval $[j-1,j)$. Observe that if
   $i(\{\gamma_N\}_{N\ge1})=f$, then $i(\{\gamma_{2N}\}_{N\ge1})=f\circ
   D_{\frac{1}{2}}$. Therefore, the proof of~\eqref{eq1new} above would
   become trivial, if $\omega$ were a $D_{\frac{1}{2}}$-invariant state
   on $L_\infty[0,\infty)$.
\end{rems}
\begin{defn}\label{saaodtd}
   Dixmier trace of a self-adjoint operator $T\in\clL^{(1,\infty)}(\clH)$
   is $\tau_\omega(T) :=
   \tau_\omega(T_+) - \tau_\omega(T_-)$ and Dixmier trace of an arbitrary
operator
   $T\in\clL^{(1,\infty)}(\clH)$ is
   $\tau_\omega(T) := \tau_\omega(Re(T)) + i\tau_\omega(Im(T)).$
\end{defn}
\begin{prop}\label{DTisSST}
   Dixmier trace $\tau_\omega$ is a (singular) symmetric functional on
   $\clL^{(1,\infty)}(\clH)$, such that $\tau_\omega(ST)=\tau_\omega(TS)$ for
   every $T\in\clL^{(1,\infty)}(\clH)$, $S\in B(\clH)$.
\end{prop}
\begin{proof} The inequality $|\tau_\omega(x)|\le\|x\|_{(1,\infty)}$,
     $x\in\clL^{(1,\infty)}(\clH)$, and further that $\tau_\omega$ is a
     symmetric functional follows immediately from
     Definition~\ref{dixtrd} and
     Corollary~\ref{DualityCor}. Since every operator $S\in\clL(\clH)$ is a
     linear combination of four unitary operators \cite[VI.6]{RS1}, it
     is sufficient to prove the equality $\tau_\omega(UT) =
\tau_\omega(TU)$ for a
     unitary $U.$ Since every operator from $\clL(\clH)$ is a linear
     combination of positive operators, it is sufficient to prove the
     last equality for positive $T$'s.  In this case, the latter
     equality follows immediately from the fact
     $\mu_n(UTU^*)=\mu_n(UT)=\mu_n(TU)=\mu_n(T), \ \forall n\ge 1$.
\end{proof}

Definitions~\ref{dixtrd} and~\ref{saaodtd} extend to Marcinkiewicz
spaces  $\clL^{(1,\infty)}\nt$ and further to Marcinkiewicz spaces
$M(\psi)\nt$, where $\psi\in\Omega$ satisfies
condition~\eqref{eq4nn}. More precisely, fix an arbitrary state
$\omega$ on $L_\infty(\mbR_+^*)$ satisfying conditions (1)--(3) of
Corollary~\ref{DualityCor} which is $D_{\frac{1}{2}}$-invariant.
Setting
\begin{equation}
\label{Troeq} \tau_\omega(x):=\omega\text{-}
\lim_{t\to\infty}a(x,t),\quad 0\le x\in M(\psi)\nt
\end{equation}
and repeating a slightly modified argument (see the details
in~\cite[p.~51]{DPSS}) from the proof of
Propositions~\ref{TrisLinear} and \ref{DTisSST}, we obtain an
additive homogeneous functional on $M(\psi)\nt_+$, which extends
to a symmetric functional on $M(\psi)\nt$ by linearity.
In the sequel, we refer to any functional $\tau_\omega $ defined
in~\eqref{Troeq}, where $\omega \in D(\mbR_+^*)$ as a Dixmier
trace.

Finally, we note that the duality between the dilation invariant
functionals on $L_\infty(\mbR_+^*)$ and translation invariant
functionals on $L_\infty(\mbR)$ allows an alternative definition of
Dixmier traces.
For simplicity, we consider this definition only for the space
$\clL^{(1,\infty)}\nt$, where $\nt$ is an arbitrary semifinite von
Neumann algebra.

Let $\mbL$ (respectively, $\LL$) belong to $BL(\mbR_+)$
(respectively, $BL(\mbN)$). We set
\begin{equation*}
\begin{split}
F_{\mbL}(T)&:=\mbL\text{-}\lim_{t\to\infty}
\frac{1}{\log(1+e^t)}\int_0^{e^t}\mu_s(T)ds,\\
F_\LL(T)&:=\LL\text{-}\lim_{N\to\infty} \frac{1}{\log(1+e^N)}
\sum_{n=1} ^{[e^N]}\mu_n(T),
\end{split}\quad T\in \clL^{(1,\infty)}\nt
\end{equation*}
where $[e^N]$ is the integral part of $e^N$.

\begin{thm}\label{thn5.5}[{\cite{CPS2,DPSSS1,DPSSS2,LSS}}] For every semifinite
von Neumann algebra $\nt$ and arbitrary  states $\mbL\in
BL(\mbR_+)$ and $\LL\in BL(\mbN)$, the functionals $F_\mbL$ and
$F_\LL$ are symmetric functionals on $\clL^{(1,\infty)}\nt$.
\end{thm}
The following result shows that the class of Dixmier traces on
$\clL^{(1,\infty)}\nt$ coincides with the sets of functionals
$\{F_\LL\ :\ \LL\in BL(\mbN)\}$ and $\{F_\mbL\ :\ \mbL\in
BL(\mbR_+)\}$. The proof of the first equality below follows from
the remarks at the end of the preceding section.

\begin{thm}[{\cite[Theorems~2.3,6.2]{LSS}}]\label{th5.6new}
For every semifinite von Neumann algebra $\nt$, we have
$$
\{\tau_\omega\ | \ \omega\in D(\mbR_+^*)\} =\{F_\mbL\ |\ \mbL\in
BL(\mbR_+)\}=\{F_\LL\ |\ \LL\in BL(\mbN)\}.
$$
\end{thm}
The detailed study of the class of concave functions
$\psi\in\Omega$ for which analogues of Theorems~\ref{thn5.5}
and~\ref{th5.6new} hold for similarly defined classes of symmetric
functionals on Marcinkiewicz spaces $M(\psi)\nt$ is contained
in~\cite{DPSSS1,DPSSS2,LSS}.

\subsection{Connes-Dixmier traces} \label{Connes-Dixmier traces}
We have shown in the preceding subsection that with every
state~$\omega\in D(\mbR_+^*)$ (respectively,
  $\mbL \in BL(\mbR_+),BL(\mbR)$, $\LL \in  BL(\mbN)$) there exists an
associated Dixmier trace~$\tau_\omega$ (respectively, a symmetric
functional~$F_\mbL$, $F_\LL$). It is possible to isolate various
subsets in the sets of states~$D(\mbR_+^*)$, $BL(\mbR_+)$,
$BL(\mbR)$, $ BL(\mbN)$ and relate with them corresponding subsets
of traces. For example, let us consider the sets of all
$H$-invariant (respectively, $M$-invariant) states
on~$L_\infty(\mbR)$ (or~$\ell_\infty$) (respectively,
$L_\infty(\mbR_+^*)$).  It is easy to see that
\begin{equation}
   \label{FirstEmbed}
   \{ \omega:\ \
   \omega\ \ \text{is a $H$-invariant state on~$L_\infty(\mbR)$
     (resp.~$\ell_\infty$)}\} \subsetneqq \{\omega:\ \ \omega \in BL(\mbR)\
   \ (\text{resp.~$ BL(\mbN)$})\}
\end{equation}
and
\begin{equation}
   \label{SecondEmbed}
   \{\omega:\ \ \text{$\omega$ is an $M$-invariant state
     on~$L_\infty(\mbR_+^*)$}\} \subsetneqq \{\omega:\ \ \omega \in
     D(\mbR_+^*)\}.
\end{equation}
Indeed, suppose that~$0\le \omega\in \ell_\infty^*$ is such
that~$\omega({\bf 1}) = 1$ and~$\omega(Hx) = \omega(x)$, for every~$x \in
\ell_\infty$.  To prove~$\omega(Tx) = \omega(x)$, $x\in\ell_\infty$,
it is sufficient to show that~$\omega(HTx) = \omega(Hx)$,
$x\in\ell_\infty$.  However, a straightforward calculation yields $$
(HTx)_N - (Hx)_N = \frac{x_2+\ldots+x_{N+1}} N - \frac{x_1+ \ldots +
   x_{N}} N = \frac {x_{N+1} - x_1} N \rightarrow 0,\ \ \text{as}\ \ N
\rightarrow \infty. $$  Similarly, it can be shown that for
every~$x\in L_\infty(\mbR)$ and~$b \in \mbR$ we have $$ \lim_{t
   \rightarrow \infty} (HT_b x)(t) - (Hx)(t) = 0 $$  which
establishes~\eqref{FirstEmbed} and the inclusion~\eqref{SecondEmbed}
follows from~\eqref{FirstEmbed} via Corollary~\ref{DualityCor}.  Alain
Connes in~\cite{Co4} suggested to work with the set of states
on~$L_\infty(\mbR_+^*)$, which is  larger then the set on the
left hand side of~\eqref{SecondEmbed}.  Namely, let us consider the
following class of states on~$L_\infty(\mbR_+^*)$
$$ CD(\mbR_+^*) :=
\{\tilde \omega = \gamma
\circ M\ :\ \text{$\gamma$ is an arbitrary singular state
   on~$C_b[0, \infty)$}\}. $$  It is still easy to verify
that~$CD(\mbR_+^*) \subsetneqq D(\mbR_+^*)$, and then infer the
proper inclusion $CBL(\mbR_+) \subsetneqq BL(\mbR_+)$, where $$
CBL(\mbR_+) := \{\mbL = \gamma \circ H\ : \ \text{$\gamma$ is an
arbitrary singular
   state on~$C_b(\mbR_+)$}\},
$$  from
Proposition~\ref{DualityProp}.  We refer to the subclass of Dixmier traces
$$
\{\tau_\omega\ : \omega \in CD(\mbR_+^*)\}
$$
  as the class of
Connes-Dixmier traces.  The following theorem shows that an analogue
of Theorem~\ref{th5.6new}  also holds for
the class of Connes-Dixmier traces.
\begin{thm}
   [{\cite[Theorems~5.6, 6.2]{LSS}}]
   \label{MustBeReLabeled}
   For every semi-finite von Neumann algebra~$\nt$, we have $$ \{
   \tau_\omega:\ \ \omega \in CD(\mbR_+^*)\} = \{F_\mbL:\ \ \mbL \in
   CBL(\mbR_+)\}. $$
\end{thm}

We complete this subsection with the remark that there is another
natural subclass of Dixmier traces which is associated with the
subset of states on $L_\infty(\mbR_+^*)$ appearing in
Corollary~\ref{DualityCor}
$$
\{\omega \in L_\infty(\mbR_+^*)^*:\ \ \text{$\omega$ is an
   $M$-invariant and $P^a$-invariant state on~$L_\infty(\mbR_+^*)$, $a
   >0$}\},
$$
or equivalently, with the set
$$
\{\mbL\in L_\infty(\mbR)^*\ :\ \text{$\mbL$ is an $H$-invariant and
   $D_a$-invariant state on $L_\infty(\mbR)$, $a>0$}\}.
$$
The class of Dixmier traces
associated with the latter set is further referred to as the class of
maximally invariant Dixmier traces. Clearly, the latter class is
contained in the class of Connes-Dixmier traces.  Maximally invariant
Dixmier traces are used in Sections~\ref{DTandRofZF} and~\ref{HSFsec}.

\subsection{Rearrangement invariant functionals and singular traces.}
\label{sec:RIfuncsAndSingTraces}

The class of Dixmier traces (and its subclasses) studied in the
preceding subsections is a special subclass of the general class of
(singular) symmetric functionals on the (fully symmetric) operator
space~$\clL^{(1, \infty)} \nt$ (or, on operator ideal~$\clL^{(1,
   \infty)} (\clH) $).  The latter class is, in its turn, a subclass of
the class of r.i.\ functionals on~$\clL^{(1, \infty)} \nt$.  In
the case of the operator ideal~$\clL^{(1, \infty)} (\clH)$, the
latter class may be viewed as the set of all (singular) traces
on~$\clL(\clH)$, which take finite values on the elements
from~$\clL^{(1, \infty)}(\clH)$. A literature devoted to general
singular traces on~$\clL(\clH)$ is tremendous. We limit our list
of papers from this area to the following articles \cite{GI1, GI2,
AGPS, AGPS1, AGPS2, CGS, Va}. In many cases, the theory of
singular traces runs in parallel with the theory of symmetric
functionals. For instance, results of~\cite{AGPS,Va}
(respectively, \cite{GI1}) concerning  necessary and sufficient
conditions on a positive compact operator $T$ (respectively,
positive $\tau$-measurable operators $T\in L_1\nt +\clN$) for the
existence of a singular trace which takes a finite non-zero value
on $T$ are in fact very close to the result of
Theorem~\ref{MSCwithNTF}. It should be noted that the
results~\cite{GI1} are stated for general $\tau$-measurable
operators and not just for the operators from $L_1\nt+\clN$,
however not all the results there treating this general case are
supplied  with a reliable proof
(e.g.~\cite[Proposition~3.3]{GI1}).

Finally, we point out at the important connections between the theory
of general singular traces and the study of the structure of
commutator spaces for operators acting on a Hilbert space (see the
comments made at the end Section~\ref{GenFactsOnSF}). The latter study is
related
to cyclic homology and algebraic $K$-theory of operator ideals and is
beyond the scope of the present survey. For details, we
recommend~\cite{DFWW, K1, DK, DK2, Fa04FA}.
\section{Class of measurable elements.}
In this section, we briefly review the notion of measurable operators
introduced by A.Connes~\cite{Co4}.

\begin{defn}\label{mesopd1} $T\in\DixIdeal{1}$ is called
\newn{(Dixmier)-measurable}
if $\tau_\omega(T)$ does not depend on the choice
of $\omega\in D(\mbR_+^*)$.
\end{defn}
\begin{defn}[\cite{Co4}]\label{mesopd2}
   $T\in\DixIdeal{1}$ is called \newn{(Connes-Dixmier)-measurable} if
   $\tau_\omega(T)$ does not depend on the choice of
   $\omega\in CD(\mbR_+^*)$.
\end{defn}
\begin{rems}\label{rem6.3}
\begin{thlist}
\item It is obvious that the sets of measurable operators defined
   above are linear spaces, which are, in fact, closed subspaces of
   $\clL^{(1,\infty)}(\clH)$. However, these subspaces are not order
   ideals, in other words, the fact that a self-adjoint operator $A$ is
   measurable does not necessarily imply that $A_+$ and $A_-$ are
   measurable operators.  Example. Take a positive non-measurable
   diagonal operator $A = \diag\{a_1,a_2,a_3,\ldots\}$ from
   $\DixIdeal{1}.$ Define a diagonal operator $B$ by
   $B=\diag\{a_1,-a_1,a_2,-a_2,\ldots\}$. Evidently, $B$ is measurable,
   moreover $\tau_\omega(B)=0$ for all $\omega.$ However, the positive
   and negative parts of $B$ are not measurable.
\item Definitions of Dixmier and Connes-Dixmier measurable
operators naturally extend~\cite{LSS} to Marcinkiewicz spaces
$\clL^{(1,\infty)}\nt$, where $\nt$ is an arbitrary semifinite von
Neumann algebra, and further, to operator Marcinkiewicz spaces
$M(\psi)\nt$, for all $\psi$ satisfying condition~\eqref{eq1.1}.
\end{thlist}
\end{rems}

It is obvious, that every Dixmier-measurable operator is also
Connes-Dixmier measurable. Our objective in the present section is to describe
the classes of positive Connes-Dixmier measurable operators and positive
Dixmier-measurable operators and to show that these two classes actually
coincide. To this end, we will need two auxiliary results.
\begin{thm}\cite[section 6.8]{Hardy} \label{HardyThm}
   Let $b(t)$ be a positive piecewise differentiable function such that
   $tb'(t) > -H$ for some $H>0$ and all $t>C,$ where $C$ is a constant.
   Then
   $$
      \liminfty{t} \frac 1t \int_0^t b(s)\,ds = A
     \quad\text{for some}\quad A\geq 0 \quad\text{if and only if}\quad
      \liminfty{t} b(t) = A.
   $$
\end{thm}
Imitating the Lorentz definition of almost convergent sequences, a
positive function $f\in C_b[0,\infty)$ is said to be \newn{almost
   convergent} if all  states from $BL(\mbR_+)$ take the same value
on this function.
\begin{thm}\label{LorentzThm}{\cite[Theorem~3.3]{LSS}}
   If a function $f\in C_b[0,\infty)$ is almost convergent to a number $A$
then, the following limit
   $$
      \lim_{t\to\infty} \frac 1t \int_0^{t} f(s)\,ds
   $$
exists and is equal to $A$.
\end{thm}
\begin{proof}
Suppose the result is false. Then there exists a constant $c\ne A$ such
that $(Hf)(t_n)\to c$ for some sequence $t_n\uparrow \infty$. Take the
unit ball $B$ of $C_b([0,\infty))^*$ and consider the sequence of
functionals $\sigma_{t_n}(x)=x(t_n)$, $n\ge1$ from $B$. Since $B$ is
weak$^*$-compact, this sequence has a limit point $V\in B$. It is easy
to see that $V\ge0$, $V(1)=1$, $V(p)=\lim_{n\to\infty} p(t_n)=0$ for
every $p\in C_0[0,\infty)$ and also that $V(H(f))=\lim_{n\to\infty}
H(f)(t_n)=c$. Define the functional $L$ on $L_\infty(\mbR_+)$ by setting
$L(x):=V(H(x))$. It is easy to verify
that $L$ is a state from $BL(\mbR_+)$  and
that $V(f)=c\ne A$. Thus, the supposition that the result does not hold
is false.
\end{proof}

The following theorem is the main result of this section.
\begin{thm} \cite{LSS}\label{th6.5}
   A positive operator $T$ from $\clL^{(1,\infty)}\nt$ is
Dixmier-measurable if and only if the limit
   $$
      \liminfty{t} \ILog \DixInt{T}
   $$
   exists.
\end{thm}
\begin{proof}
   The ``if'' part of the assertion is trivial.  Now, fix an operator
   $T\in\clL^{(1,\infty)}_+\nt$ such that for
   $g(t):=\frac{1}{\log(1+t)}\int_0^t\mu_s(T)ds$, we have
   $\tau_\omega(T)=\omega\text{-}\lim_{t\to\infty} g(t)=A\ge0$ for every
   $\omega\in CD(\mbR_+^*)$. It follows from the remarks made in
   Section~\ref{Connes-Dixmier traces} that for all
  $\mbL\in BL(\mbR_+)$, we have
   $Tr_\mbL(T):=\mbL\text{-}\lim_{\lambda\to\infty} g(e^\lambda)=A$,
   and therefore, by Theorem~\ref{LorentzThm}, we obtain
   $$
     \liminfty{u} \frac 1u \int_0^u\left(\frac 1{\log(1+e^\lambda)}
     \int_0^{e^\lambda} \mu_s(T)\,ds\right)\,d\lambda=A.
   $$
   Setting,
   $$
      b(\lambda) := \frac 1{\log(1+e^\lambda)} \int_0^{e^\lambda}
      \mu_s(T)\,ds,\quad\lambda>0
   $$
   we have
   \begin{multline*}
  \lambda b'(\lambda)
     \geq \lambda\frac d{d\lambda}\left(\frac 1{\log(1+e^\lambda)}\right)
     \int_0^{e^\lambda} \mu_s(T)\,ds
     =  - \frac {\lambda e^\lambda}{(1+e^\lambda)\log^2(1+e^\lambda)}
     \int_0^{e^\lambda} \mu_s(T)\,ds \\
     \geq - \frac {\lambda}{\log(1+e^\lambda)} \cdot \frac
     {1}{\log(1+e^\lambda)}\int_0^{e^\lambda} \mu_s(T)\,ds
     \geq - \DixNorm{T}.
   \end{multline*}
   Applying Theorem~\ref{HardyThm}, we now infer that
   $\lim_{\lambda\to\infty} b(\lambda)=A$, and therefore
   $\lim_{t\to\infty} g(t)=A$.
\end{proof}
We shall now show that a similar argument as in the proof above yields a
stronger result.
\begin{thm}[\cite{LSS}]\label{thTisCDmeasurable}  A positive operator $T$
from $\DixIdealN{1}$
is Connes-Dixmier-measurable if and only if the limit
  $$
      \liminfty{t} \ILog \DixInt{T}
   $$
   exists.
\end{thm}
\begin{proof}
We need only to show the ``only if'' part. We shall use the notations
$g(\cdot)$ and $b(\cdot)$ introduced in the proof of Theorem~\ref{th6.5}.
Suppose that $T\in \clL^{(1,\infty)}_+\nt$ satisfies the equality
$A=\gamma\circ M(g)=\gamma\circ LHL^{-1}(g)$ for every state $\gamma \in
C_b^*[0,\infty)$ vanishing on $C_0(0,\infty)$.

Note that if $\gamma$ is dilation invariant state, then
$\gamma\circ L$ is a translation invariant state. This remark (and
the fact that $L: L_\infty(\mbR)\to L_\infty(\mbR_+^*)$ is an
isomorphism) show that $HL^{-1}(g)=Hb$ is almost convergent.
Applying Theorem~\ref{LorentzThm} to the function $Hb$, we see
that the limit $\lim_{t\to\infty} HH(b)(t)$ exists. Assume, for a
moment, that we have already verified the assumption of
Theorem~\ref{HardyThm} for the function $Hb$. Then, we infer from
that theorem that the limit $\lim_{t\to\infty} H(b)(t)$ also
exists, and repeating the application of the same theorem (as in
the proof of Theorem~\ref{th6.5}), we conclude that there exists
the limit $\lim_{t\to\infty} b(t)$, and hence the limit
$\lim_{\lambda\to\infty} g(\lambda)$.

It remains to verify the assumption of Theorem~\ref{HardyThm} for
$Hb$. We have for all $\lambda\ge1$
\begin{align*}
\lambda(Hb)'(\lambda)&=\frac{\lambda b(\lambda)-\int_0^\lambda
   b(s)ds}{\lambda}
\ge-\|b\|_\infty\ge -\|T\|_{(1,\infty)}.
\end{align*}
\end{proof}
\begin{cor} The set of all positive Dixmier measurable operators and the
set of all positive Connes-Dixmier measurable operators coincide.
\end{cor}
The following questions are open:
\begin{thlist}
\item Do the spaces of Dixmier measurable and Connes-Dixmier
measurable operators coincide? \item What is the description of
the set of all (positive) operators, which are measurable with
respect to the set of all maximally invariant Dixmier traces?
\item What is the description of the set of all (positive)
operators measurable with respect to the set of all symmetric
functionals?
\end{thlist}
One might be tempted to pose question (iii) in greater generality
and ask whether the concept of operators measurable with respect
to the set of all r.i. functionals can be sensibly formulated.
There is an example due to \cite{KS} (see the end of Section
\ref{PrelimOnSpAndF} above) which answers this in the negative.

For results extending Theorems~\ref{th6.5}
and~\ref{thTisCDmeasurable} to Marcinkiewicz spaces $M(\psi)\nt$,
with $\psi\in\Omega$ satisfying condition~\eqref{eq1.1s}, we refer
the reader to~\cite{LSS}.

\section{Norming properties of Dixmier and Connes-Dixmier functionals}

A reader may have an impression that Dixmier and Connes-Dixmier traces
form a very ``thin'' subset of the unit sphere of the dual space. Such an
impression is wrong as established by Theorems~\ref{NormingProperty}
and~\ref{NormingPropertyForConnesDixmierTraces} below.  We shall need
the following theorem of Sucheston.
\begin{thm}[\cite{Su67AMM}] \label{SuchestonThm} For $x\in \ell_\infty$
   $$
     \sup_{\LL\in BL(\mbN)}  \LL(x) = \liminfty{n}\left(\sups{m} \frac
       1n \sums{j=1}^n x_{m+j}\right).
   $$
\end{thm}
The following proposition easily follows from its ``commutative''
counterpart, which, in its turn, can be obtained from~\cite{KPS}.
\begin{prop}\label{distanceTh}
   Let $T\in\DixIdealN{1}.$ The distance from $T$ to the subspace
   $\clL_0^{(1,\infty)}(\clN,\tau)$ in the norm $\DixNorm{\cdot}$ is
   equal to
   $$
   \rho(T):=\limsup\limits_{t\to\infty} \ILog\DixInt{T}.
   $$
\end{prop}
Recall, that in the special case
$\clL^{(1,\infty)}\nt=\clL^{(1,\infty)}(\clH)$, the space
$\clL_0^{1,\infty}(\clH)$ is the closed linear span in
$\clL^{(1,\infty)}(\clH)$ of the set of all finite-dimensional
operators.
\begin{thm} \cite{LSS} \label{NormingProperty}
   Let $T\in\DixIdealN{1}.$ The distance from $T$ to the subspace
   $\clL_0^{(1,\infty)}(\clN,\tau)$ in the norm $\DixNorm{\cdot}$ is
   equal to $\sup\{\tau_\omega(|T|)\ :\ \omega\in D(\mbR_+^*)\}$.
\end{thm}
\begin{proof}
   It is sufficient to consider the case $T\ge0$.  We note first that
   $$
   \sup\{\tau_\omega(T)\ :\ \omega\in D(\mbR_+^*)\}=\sup\{Tr_\mbL(T)\
   :\ \mbL\in BL(\mbR_+)\}=\sup\{Tr_\LL(T)\ :\ \LL\in BL(\mbN)\}.
   $$
   It is clear that
   $$
   q(T):=\sup_{\LL\in BL(\mbN)}\Tr_\LL(T) \ \leq \ \limsup\limits_{t\to\infty}
   \ILog\DixInt{T}\ (=\rho(T)).
   $$

   We have to prove the reverse inequality $q(T)\geq \rho(T).$ By
   Sucheston's theorem \ref{SuchestonThm}, it is enough to prove that \
   $\forall\,\epsilon>0 \ \exists\, N\in\mbN \ \forall\, n\geq N \
   \exists\, m\in \mbN$ \ such that
   $$
   \frac 1n \sums{j=1}^n \frac 1{m+j} \ints0^{e^{m+j}} \mu_s(T)\,ds
   \geq \rho(T)-\epsilon.
   $$
   For this purpose, it is enough to put $N=1$ and to take $m$ such
   that $\frac m{m+n}>\frac {\rho(T)-\epsilon}{\rho(T)-\epsilon/2}$ and
\begin{equation}\label{500}
   \frac 1{\log(1+e^m)} \int_0^{e^m} \mu_s(T)\,ds >
   \rho(T)-\epsilon/2.
\end{equation}
   Then
   \begin{multline*}
     \frac 1n \sums{j=1}^n \frac 1{m+j} \ints0^{e^{m+j}} \mu_s(T)\,ds
     \geq \frac 1n \sums{j=1}^n \frac 1{m+n} \ints0^{e^m} \mu_s(T)\,ds
     = \frac m{m+n}\cdot\frac 1m \ints0^{e^m} \mu_s(T)\,ds >
     \rho(T)-\epsilon.
   \end{multline*}
To verify that selection of $m$ satisfying~\eqref{500} is feasible,
first we locate a sequence $1\le t_1< t_2\ldots\uparrow\infty$, such
that
\begin{equation}\label{501}
  \lim_{k\to\infty}\frac{1}{\log(1+t_k)}\int_0^{t_k}\mu_s(T)ds>\rho(T)-\epsilon/4,\quad k\ge1,
\end{equation}
(this may be done due to Proposition~\ref{distanceTh}). For every $k$,
we define $m_k\in\mbN$, so that  $e^{m_k-1}\le t_k\le e^{m_k}$. Then
\begin{align*}
\frac{1}{\log(1+t_k)}\int_0^{t_k}\mu_s(T)ds&
\le\frac{1}{\log(1+e^{m_k-1})}\int_0^{e^{m_k}}\mu_s(T)ds\\
&=\frac{\log(1+e^{m_k})}{\log(1+e^{m_k-1})}\cdot\frac{1}{\log(1+e^{m_k})}\int_0^{e^{m_k}}\mu_s(T)ds.
\end{align*}
Since $\frac{\log(1+e^{m_k})}{\log(1+e^{m_k-1})}\to 1$, we see
that~\eqref{500} follows from~\eqref{501}.
\end{proof}

A natural question is whether the norming property remains true for
the class of Connes-Dixmier traces, as answered below.
\begin{thm} \cite{LSS} \label{NormingPropertyForConnesDixmierTraces}
Let $T\in\DixIdealN{1}.$ The distance from $T$ to the subspace
$\clL_0^{(1,\infty)}(\clN,\tau)$ in the norm $\DixNorm{\cdot}$ is
equivalent to
$\sups{\omega}\tau_\omega(T),$ where the supremum is taken over all
singular states $\omega=\gamma\circ \CesaroRPlus,$ where $\gamma$
is a singular state on $C_b[0,\infty).$
\end{thm}
As in the preceding section, the results given in
Theorems~\ref{NormingProperty}
and~\ref{NormingPropertyForConnesDixmierTraces} admit an extension
to Marcinkiewicz spaces $M(\psi)\nt$ with $\psi\in\Omega$
satisfying condition~\eqref{eq1.1s}.

We finish this section with the comment that it is not clear yet,
whether the difference in the results of Theorems~\ref{NormingProperty}
and~\ref{NormingPropertyForConnesDixmierTraces}
signify that the set of all Dixmier traces is different from the set
of all Connes-Dixmier traces.
{

\newcommand{\npartial}{\slash\!\!\!\partial}
\newcommand\Trig{\operatorname{Trig}}
\newcommand\maF{\clF }
\newcommand\CAP{{\mathcal AP}}
\newcommand\ep{\epsilon}
\newcommand\cM{\clM }
\newcommand\cnn{\clN }
\newcommand{\IS}{\mathbb{S}}
\newcommand\cS{\clS }
\newcommand\cA{\clA }
\newcommand{\Heis}{\operatorname{Heis}}
\newcommand{\Solv}{\operatorname{Solv}}
\newcommand{\Spin}{\operatorname{Spin}}
\newcommand{\SO}{\operatorname{SO}}
\newcommand{\ind}{\operatorname{ind}}
\newcommand{\Index}{\operatorname{index}}
\newcommand{\ch}{\operatorname{ch}}
\newcommand{\rank}{\operatorname{rank}}

\newenvironment{alphlist}{
\renewcommand{\theenumi}{{\rm\alph{enumi}}}
\renewcommand{\labelenumi}{\theenumi)}
\renewcommand{\labelenumii}{\theenumii)}
\begin{enumerate}}{\end{enumerate}}
\newcommand\End{\operatorname{End}}
\newcommand{\abs}[1]{\lvert#1\rvert}
  \newcommand{\A}{\clA }
         \newcommand{\D}{\clD }
\newcommand{\HH}{\clH }
\newcommand\wM{\widetilde{M}}
         \newcommand{\B}{\clB }
         \newcommand{\K}{\clK }
\newcommand{\oo}{\clO }\newcommand\maH{\clH }
         \newcommand{\coker}{{\mbox coker}}
         \newcommand{\p}{\partial}
         \newcommand{\dd}{|\D|}
         \newcommand{\n}{\parallel}
\newcommand{\bma}{\left(\begin{array}{cc}}
\newcommand{\ema}{\end{array}\right)}
\newcommand{\bca}{\left(\begin{array}{c}}
\newcommand{\eca}{\end{array}\right)}
\newcommand{\sr}{\stackrel}
\newcommand{\da}{\downarrow}
\newcommand{\tD}{\tilde{\D}}
         \newcommand{\R}{\mathbb R}
         \newcommand{\C}{\mathbb C}
         \newcommand{\h}{\mathbb H}
\newcommand{\Z}{\mathbf Z}
\newcommand{\N}{\mathbf N}
\newcommand{\tto}{\longrightarrow}
\newcommand{\ben}{\begin{displaymath}}
         \newcommand{\een}{\end{displaymath}}
\newcommand{\be}{\begin{equation}}
\newcommand{\ee}{\end{equation}}

         \newcommand{\bean}{\begin{eqnarray*}}
         \newcommand{\eean}{\end{eqnarray*}}
\newcommand{\nno}{\nonumber\\}
\newcommand{\bea}{\begin{eqnarray}}
         \newcommand{\eea}{\end{eqnarray}}
\newcommand{\supp}[1]{\operatorname{#1}}
\newcommand{\norm}[1]{\parallel\, #1\, \parallel}
\newcommand{\ip}[2]{\langle #1,#2\rangle}
\setlength{\parskip}{.3cm}
\newcommand{\nc}{\newcommand}
\nc{\nt}{\newtheorem}
\nc{\gf}[2]{\genfrac{}{}{0pt}{}{#1}{#2}}
\nc{\mb}[1]{{\mbox{$ #1 $}}}
\nc{\real}{{\mathbb R}}
\nc{\comp}{{\mathbb C}}
\nc{\ints}{{\mathbb Z}}
\nc{\Ltoo}{\mb{L^2({\mathbf H})}}
\nc{\rtoo}{\mb{{\mathbf R}^2}}
\nc{\slr}{{\mathbf {SL}}(2,\real)}
\nc{\slz}{{\mathbf {SL}}(2,\ints)}
\nc{\su}{{\mathbf {SU}}(1,1)}
\nc{\so}{{\mathbf {SO}}}
\nc{\hyp}{{\mathbb H}}
\nc{\disc}{{\mathbf D}}
\nc{\torus}{{\mathbb T}}
\newcommand{\tk}{\widetilde{K}}
\newcommand{\boe}{{\bf e}}\newcommand{\bt}{{\bf t}}
\newcommand{\vth}{\vartheta}
\newcommand{\CGh}{\widetilde{\CG}}
\newcommand{\db}{\overline{\partial}}
\newcommand{\tE}{\widetilde{E}}
\newcommand{\tr}{\mbox{tr}}
\newcommand{\ta}{\widetilde{\alpha}}
\newcommand{\tb}{\widetilde{\beta}}
\newcommand{\txi}{\widetilde{\xi}}
\newcommand{\hV}{\hat{V}}
\newcommand{\la}{\langle}
\newcommand{\ra}{\rangle}
\newcommand{\lp}{L_p(\clN ,\tau )}
\newcommand{\vlp}{\Vert _{_{L_p(\clN ,\tau )}}}
\newtheorem*{rem*}{Remark}
\section{Fredholm modules and spectral triples}

\subsection{Notation and definitions}
Let
  $\clN $ be a semifinite von Neumann algebra on a separable Hilbert
space $\mathcal H$
and let $\lp$ be a non-commutative $L_p$-space associated with
$(\clN ,\tau )$, where $\tau$ is a faithful, normal semifinite trace on
$\clN $.
Let $\clA $ be a unital Banach
$*$-algebra which is represented in $\clN $
via a continuous $*$-homomorphism $\pi$ which,
without loss of generality, we may assume to be faithful.
Where no confusion arises we suppress $\pi$ in the notation.
The fundamental objects of our analysis are explained in the following
definition.

Let $\clK _{\mathcal N }$ be the
$\tau$-compact operators in $\clN $
(that is the norm closed ideal generated by the projections
$E\in\mathcal N$ with $\tau(E)<\infty$).

\begin{defn}\label{SpectTripIsumd}
\begin{deflist}
\item A semifinite odd spectral triple $(\A,\HH,\D)$ is given by a Hilbert
   space $\HH$, a $*$-algebra $\A\subset \clN$ where $\clN$ is a
   semifinite von Neumann algebra acting on $\HH$, and a densely
   defined unbounded self-adjoint operator $\D$ affiliated to $\clN$
   such that
\begin{deflist}
\item $[\D,a]$ is densely defined and extends to a bounded operator
   for all $a\in\A$
\item $(\lambda-\D)^{-1}\in\K_\clN$ for all $\lambda\not\in{\R}$
\end{deflist}
\item We say that $(\A,\HH,\D)$ is even if in addition there is a
   ${\Z}_2$-grading such that $\A$ is even and $\D$ is odd. That is an
   operator $\Gamma$ such that $\Gamma=\Gamma^*$, $\Gamma^2=1$, $\Gamma
   a=a\Gamma$ for all $a\in\A$ and $\D\Gamma+\Gamma\D=0$.
\item If $\mathcal I$ is a symmetrically normed ideal in ${\mathcal
     K}_\clN $ then we say that the spectral triple $(\clA ,\clN ,\clD
   )$ is $\mathcal I$-summable if $(1+\clD ^2)^{-1/2}\in \mathcal I$.
\end{deflist}
\end{defn}
\begin{rems}
In \cite{BeF} the terminology for the concept we have just introduced is
`von Neumann spectral triple'.
The two most important special cases of this
definition are when
\begin{alphlist}
\item $\mathcal I$  is the ideal $L^p$ in which case we say $({\mathcal
A},\clN ,\clD )$  is $p$-summable and
\item the case
where the Dixmier trace was first evident namely $\mathcal I=
\clL^{p,\infty}$. In this case, we say that $({\mathcal A},\clN ,\clD )$ is 
$(p, \infty)$-summable.
\end{alphlist}
There is a third case which is not relevant for the discussion here
and that is the notion of theta summability \cite{Co2,Co3,Co4}.
\end{rems}
\noindent{\bf Note}.
In this paper, for simplicity of exposition,
we will deal only with unital algebras
$\clA \subset \clN$ where the identity of $\mathcal A$ is that of $\clN$.
We have adopted a notational convention correlated to the context in which
we are working.
A calligraphic $\D$ will always denote
an unbounded self-adjoint operator forming part of a
semifinite spectral triple $(\A,\HH,\D)$. A roman $D$ will
denote a self-adjoint operator on a Hilbert space, usually with some side
conditions.

In the original paper on noncommutative geometry by Alan Connes
\cite{Co1} the notion of spectral triple was introduced as an
`unbounded Fredholm module' (see also \cite{Co4}). The study of
semifinite spectral triples was initiated in \cite{CP1} in the
context of spectral flow \cite{APS1,APS3,Ph}
 and carried further in \cite{BeF}
  with a focus on applications to foliations.
In \cite{Co1} the notion of Fredholm module was introduced
based on the idea of a Kasparov module. This may be
generalised to the semifinite case as follows.

\begin{defn} (\cite{CP1}, \cite{Suk}) A
bounded $p$-summable pre-Breuer-Fredholm module for $\clA $, is a
pair $(\clN ,F_0)$ where $F_0$ is a bounded self-adjoint operator
in $\clN $ satisfying:

(1) $|{\bf 1}-F_0^2|^{1/2}$
  belongs to  $\lp$; and
\hfill\break
\noindent (2)
$\clA _p:=\{a\in \clA \ |\ [F_0,a]\in \lp\}$) is a dense
$*$-subalgebra of $\clA $.

When $F_0^2={\bf 1}$ we drop the prefix \lq  pre-'.
\end{defn}

\noindent In the special case when $\clN  =\clL (\clH )$ and $\tau
$ is the
  standard trace Tr, we omit \lq\lq Breuer" from the
  definition.
In this case, the non-commutative $L_p$-space
coincides with the Schatten-von Neumann ideal $ \clC _p$ of compact
operators and the definition originates from Definition 3 from
\cite{Co4} p.290.
\subsection{Bounded versus unbounded}
The relationship between the bounded and unbounded pictures
is worth some comment. There are two approaches to this question,
that in \cite{SWW} (and explained in more detail in \cite{GVF})
  and a more general approach via
perturbation theory in \cite{CPS1}. In the case of even spectral
triples this matter was settled by Connes (see \cite{Co1} $I$.6).
Let $D$ and $D_0$ be unbounded self adjoint operators on $\clH$
differing by a bounded operator in $\mathcal N$. We study directly
the map $\phi$ defined by $ \phi(D)= D(1+D^2)^{-1/2} $
and the difference $ \phi(D)-\phi (D_0). $ In \cite{CPS1} a number
of results were established, the
  one most relevant to this survey being:
\begin{thm}
With the assumptions above on $D$ and $D_0$ and with $D-D_0\in \mathcal N$
and $1<p<\infty$
we have
$$
\Vert \phi(D)-\phi(D_0)\vlp
\leq {Z}_p\max \{ \Vert D-D_0\Vert ^{1/2}, \ \Vert D-D_0\Vert \}
\cdot \Vert (1+D^2)^{-1/2}\vlp.
$$
for some positive constant ${Z}_p$ which depends
on $p$ only.
\end{thm}

 From this theorem (cf \cite{SWW}) one obtains the
\begin{cor} If
$1<p<\infty$ and  $(\clA ,\clN  , \D_0)$ is an odd semifinite
$p$-summable spectral triple
for the Banach $*$-algebra
  then $(\clN  , {\rm sign}(\D_0))$ is an odd  bounded
$p$-summable Breuer-Fredholm module for  $\clA $.
\end{cor}

\subsection{More on Semifinite Spectral Triples}
\begin{defn}\label{qck} A semifinite spectral triple $(\A,\HH,\D)$ is $QC^k$
for $k\geq 1$
($Q$ for quantum) if for all $a\in\A$
the operators $a$ and $[\D,a]$ are in the domain of $\delta^k$, where
$\delta(T)=[\dd,T]$ is the partial derivation on $\clN$ defined by $\dd$.
\end{defn}

\noindent{\bf Notes}.
(i) The notation is meant to be analogous to the classical case, but
we introduce
the $Q$ so that there is no confusion between quantum differentiability of
$a\in\A$ and classical differentiability of functions.\\
(ii) By a partial derivation we mean that $\delta$ is defined on
some subalgebra of $\clN$ which need not be (weakly) dense in
$\clN$. More precisely, $dom\,\delta=
\{T\in\clN:\delta(T)\mbox{ is bounded}\}$. \\

\noindent{\bf Observation} If $T\in\clN $, one can show that
$[\dd,T]$ is bounded if and only if $[(1+\D^2)^{1/2},T]$ is
bounded, by using the functional calculus to show that
$\dd-(1+\D^2)^{1/2}$ extends to a bounded operator in $\clN$. In
fact, writing $\dd_1=(1+\D^2)^{1/2}$ and $\delta_1(T)=[\dd_1,T]$
we have $dom\,\delta^n=dom\,\delta_1^n$ for all $n.$

\begin{proof} Let $f(\D)=(1+\D^2)^{1/2}-\dd$, so, as noted above, $f(\D)$
extends to a bounded operator in $\clN$. Since
$$ \delta_1(T)-\delta(T)=[f(\D),T]$$ is always bounded, $dom\,\delta=
dom\,\delta_1$. Now $\delta\delta_1=\delta_1\delta$, so
\begin{eqnarray*}
\delta^2_1(T)-\delta^2(T)&=&\delta_1(\delta_1(T))-\delta_1(\delta(T))+
\delta_1(\delta(T))-\delta(\delta(T))\\
&=&[f(\D),\delta_1(T)]+[f(\D),\delta(T)].\end{eqnarray*} Both
terms on the right hand side are bounded, so $dom\,\delta^2=
dom\,\delta_1^2$. The proof proceeds by induction.
\end{proof}

Thus the condition defining $QC^k$ can be replaced by \ben
a,[\D,a]\in\bigcap_{k\geq n\geq 0}dom\,\delta_1^n\ \ \ \forall
a\in\A.\een This is important as we {\em do not assume at any
point that $\dd$ is invertible}.

If $(\A,\HH,\D)$ is a $QC^k$ spectral triple, we may endow the algebra
$\A$ with the topology determined by the seminorms
\ben a\tto \n \delta^n(a)\n+\n\delta^n([\D,a])\n,\ \ \ n=0,1,2,...,k \een
By \cite[Lemma 16]{R}
we may, without loss of generality, suppose that $\A$ is
complete in the resulting topology by completing if necessary. This completion
is stable under the holomorphic functional calculus, so we
have a sensible spectral theory and
$K_*(\A)\cong K_*(\bar\A)$ via inclusion, where $\bar\A$
is the $C^*$-completion of $\A$.
  A $QC^\infty$ spectral triple is one that is  $QC^k$ for all $k=1,2,\ldots$.

\noindent{\bf Observation}
If $T\in\clN$ and $[\D,T]$ is bounded, then $[\D,T]\in\clN$.

\begin{proof} Observe that $\D$ is affiliated with $\clN$, and so commutes with
all projections in the commutant of $\clN$, and the commutant of $\clN$
preserves
the domain of $\D$. Thus if $[\D,T]$ is bounded, it too commutes with all
projections in the commutant of $\clN$, and these projections preserve the
domain of $\D$, and so $[\D,T]\in\clN$.
\end{proof}
Similar comments apply to $[\dd,T]$, $[(1+\D^2)^{1/2},T]$ and
combinations such as $[\D^2,T](1+\D^2)^{-1/2}$.
We will often simply write $\clL^1$ for the trace  ideal in order
to simplify the notation, and denote the norm on $\clL^1$ by
$\n\cdot\n_1$. Note that in the case where $\clN \neq\clL (\clH
)$, $\clL^1$ need not be
   complete in this norm but it is complete in the norm $||.||_1 +
||.||_\infty$.
  (where $||.||_\infty$ is the uniform norm).
\subsection{Summability and Dimension}

Finite summability conditions on a spectral triple
give a half-plane where the function
\be z\mapsto \tau((1+\D^2)^{-z})\label{zeta}\ee
  is well-defined and holomorphic.

\begin{defn}\label{dimension}
If $(\A,\HH,\D)$ is a $QC^\infty$ spectral triple,
we call
\ben p=\inf\{a\in{\R}:\tau((1+\D^2)^{-a/2})<\infty\}\een
the {\bf spectral dimension} of $(\A,\HH,\D)$.
\end{defn}
\section{Spectral Flow}\label{spectralflow}

One of the main motivations for extending
the study of spectral triples to the semifinite von Neumann setting
is the study of type $II$ spectral flow (this concept is due to Phillips
\cite{Ph,Ph1}). Let
$\pi:\clN \to \clN /{\mathcal K_{\clN}}$
be the canonical mapping. A Breuer-Fredholm operator is one that
maps to an invertible operator under $\pi$.
A full discussion of Breuer-Fredholm theory in a semifinite
von Neumann algebra is contained in \cite{CPRS2} extending
the discussion of the
Appendix to \cite{PR} and \cite{B1,B2}.
  As usual $D$ is an unbounded densely defined self-adjoint
Breuer-Fredholm operator  on $\HH$ (meaning $D(1+D^2)^{-1/2}$ is
bounded and Breuer-Fredholm in $\mathcal N$)
with $(1+D^2)^{-1/2}\in \clK _\mathcal N$.
For a unitary $u\in \mathcal N$ such that $[D,u]$ is a bounded operator,
  the path
$$D_t^u:=(1-t)\,D+tuDu^*$$
of unbounded self-adjoint
Breuer-Fredholm operators is continuous in the sense that
$$F_t^u:=D_t^u\left(1+(D_t^u)^2\right)^{-\frac{1}{2}}$$
is a norm continuous path of self-adjoint Breuer-Fredholm operators in
$\mathcal N$ \cite{CP1}.
Recall that the Breuer-Fredholm index of
a Breuer-Fredholm operator $F$ is defined by
$$ind(F)=\tau(Q_{ker F})-\tau(Q_{coker F})$$
where $Q_{ker F}$ and $Q_{coker F}$ are the projections onto
the kernel and cokernel of $F$.

\begin{defn}
If $\{F_t\}$ is a continuous path of self-adjoint Breuer-Fredholm
operators in $\mathcal N$, then the definition of the {\bf
spectral flow} of the path, $sf(\{F_t\})$ is based on the
following sequence of observations in \cite{P1} and
\cite{Ph1} (see also~\cite{BCPRSW,BLP}):

\noindent 1. While the function $t\mapsto sign(F_t)$ is typically
discontinuous, as is the
projection-valued mapping $t\mapsto P_t=\frac{1}{2}(sign(F_t)+1)$,
for $F_t=2P_t-1$ with $P_t$ the non-negative
spectral projection, $t\mapsto \pi(P_t)$ is continuous.\\
2. If $P$ and $Q$ are projections in $\mathcal N$
then $PQ:Q\HH \to P\HH$ is a Breuer-Fredholm operator if and only if
$||\pi(P)-\pi(Q)||<1$ in which case
$ind(PQ)\in {\R}$ is well-defined. \\
3. If we partition the parameter interval of $\{F_t\}$ so
that the $\pi(P_t)$ do not vary much in norm on each subinterval
of the partition then
$$sf(\{F_t\}):=\sum_{i=1}^n ind(P_{t_{i-1}}P_{t_i})$$
is a well-defined and (path-) homotopy-invariant number which
agrees with the usual notion of spectral flow in the type $I_\infty$
case.\\
4. For $D$ and $u$ as above, we define the {\it spectral flow} of
the path $D_t^u:=(1-t)\,D+tuDu^*$ to be the spectral flow of the
path $F_t$ where
$F_t=D_t^u\left(1+(D_t^u)^2\right)^{-\frac{1}{2}}$. We denote this
by \ben sf(D,uDu^*)= sf(\{F_t\}),\een and observe that this is an
integer in the $\clN =\clL (\clH )$ case and a real number in the
general semifinite case.
\end{defn}

Special cases of spectral flow in a semifinite von Neumann algebra were
discussed in
\cite{M,P1,P2,Ph,Ph1}.

Let $P$ denote the projection onto the nonnegative spectral
subspace of $D$. The spectral flow along $\{D_t^u\}$ is equal to
$sf(\{F_t\})$ and by \cite{CP1} this is the Breuer-Fredholm index
of $PuPu^*$. (Note that $sign F^u_1=2uPu^*-1$ and that for this
special path we have $P-uPu^*$ is compact so $PuPu^*$ is certainly
Breuer-Fredholm from $uPu^*\HH$ to $P\HH$.) Now, \cite[Appendix
B]{PR}, we have $ind (PuPu^*)= ind (PuP)$.

The operator $PuP$ is known as a generalised Toeplitz operator.
   Formulae for its index in terms of the Dixmier trace are discussed
in Sections 13 and 16. These Dixmier trace formulae follow from
the analytic formulae for spectral flow discovered in
\cite{G},\cite{CP1},\cite{CP2} which we will now explain.

\subsection{Spectral Flow Formulae}

We now introduce the spectral flow formula of Carey and Phillips,
\cite{CP1,CP2}.
This formula starts with a semifinite spectral triple $(\A,\HH,\D)$
and computes the spectral flow from $\D$ to $u\D u^*$,
where $u\in\A$ is unitary with $[\D,u]$ bounded, in the case
where $(\A,\HH,\D)$ is of dimension $p\geq 1$.
Thus for any $n>p$ we have by the extension of Theorem 9.3 of \cite{CP2}
to the case of general semifinite von Neumann algebras (see~\cite{CPRS3}):
\be sf(\D,u\D
u^*)=\frac{1}{C_{n/2}}\int_0^1\tau(u[\D,u^*](1+(\D+tu[\D,u^*])^2)^{-n/2})dt,
\label{basicformula}\ee
with $C_{n/2}=\int_{-\infty}^\infty(1+x^2)^{-n/2}dx$.
This real number $sf(\D,u\D u^*)$
recovers the pairing of the $K$-homology class
$[\D]$ of $\mathcal A$ with the $K_1(\clA )$ class $[u]$
(see below). There is a geometric way to view this formula
due originally to Getzler \cite{G}. It is shown in
\cite {CP2} that
the functional $X\mapsto \tau(X(1+(\D+X)^2)^{-n/2})$ on $\clN _{sa}$
determines an exact one-form
on an affine space $\D+\clN_{sa}$. Thus (\ref{basicformula})
represents the integral of this one-form along the path
$\{\D_t=(1-t)\D+ tu\D u^*\}$
provided one appreciates that $\dot\D_t=u[\D,u^*]$ is a tangent vector
to this path.
Moreover this formula is
scale invariant. By this we mean that if we replace $\D$ by $\epsilon\D$, for
$\epsilon>0$, in the right hand side of (\ref{basicformula}), then the left
hand
side is unchanged as is evident from the definition
of spectral flow. This is because spectral flow only involves the {\it
phase} of $\mathcal D$
which is the same as the phase of $\epsilon\mathcal D$.

\subsection{Relation to Cyclic Cohomology}
To place some of the discussion in its correct context and to
prepare for later applications we need
to discuss some aspects of cyclic
cohomology.  We will use the normalised $(b,B)$-bicomplex
(see \cite{Co4,Lo}).

We introduce the following
linear spaces.
Let $C_m=\A\otimes \bar\A^{\otimes m}$
where $\bar\A$ is the quotient $\A/\C I$ with $I$ being the
identity
element of $\A$ and (assuming with no loss of generality that $\A$ is complete
in the $\delta$-topology) we employ the projective tensor product.
Let $C^m=Hom(C_m,\C)$ be the linear space of continuous multilinear functionals
on $C_m$.
We may define the $(b,B)$ bicomplex using these spaces
(as opposed to $C_m=\A^{\otimes m+1}$ et cetera) and the resulting cohomology
will be the same. This follows because the bicomplex defined using
$\A\otimes \bar\A^{\otimes m}$ is quasi-isomorphic to that defined
using  $\A\otimes \A^{\otimes m}$.

A normalised
$\mathbf{(b,B)}${\bf-cochain}, $\phi$ is a finite collection of continuous
multilinear
  functionals on $\A$,
$$\phi=\{\phi_m\}_{m=1,2,...,M}\mbox{ with }\phi_m\in C^m.$$
It is a (normalised) $\mathbf{(b,B)}${\bf-cocycle} if, for all
$m$, $b\phi_m+B\phi_{m+2}=0$ where $b: C^m\to C^{m+1}$, $B:C^m\to
C^{m-1}$ are the coboundary operators given by

\begin{multline}\label{900}(B\phi_m)(a_0,a_1,\ldots,a_{m-1})
=\sum_{j=0}^{m-1} (-1)^{(m-1)j}\phi_m(1,a_j,a_{j+1},\ldots,a_{m-1},a_0,
\ldots,a_{j-1})\\
(b\phi_{m-2})(a_0,a_1,\ldots,a_{m-1})=\\
\sum_{j=0}^{m-2}(-1)^j\phi_{m-2}(a_0,a_1,\ldots,a_ja_{j+1},\ldots,a_{m-1})
+(-1)^{m-1}\phi_{m-2}(a_{m-1}a_0,a_1,\ldots,a_{m-2})\end{multline}
  We write $(b+B)\phi=0$ for brevity.
Thought of as functionals on $\A^{\otimes m+1}$
a normalised cocycle will satisfy
$\phi(a_0,a_1,\ldots,a_n)=0$ whenever any $a_j=1$ for $j\geq 1$.
An {\bf odd} ({\bf even}) cochain has $\{\phi_m\}=0$ for $m$ even (odd).

  Similarly, a $\mathbf{(b^T,B^T)}${\bf-chain}, $c$ is a (possibly infinite)
  collection
$c=\{c_m\}_{m=1,2,...}$ with
$c_m\in C_m$.
The $(b,B)$-chain $\{c_m\}$
is a $\mathbf{(b^T,B^T)}${\bf-cycle} if $b^Tc_{m+2}+B^Tc_m=0$ for all $m$. More
briefly, we write $(b^T+B^T)c=0$.
Here $b^T,B^T$
are the boundary operators of cyclic homology,
and are the transpose of the coboundary operators
$b,B$ in the following sense.

The pairing between a $(b,B)$-cochain $\phi=\{\phi_m\}^M_{m=1}$ and a
$(b^T,B^T)$-chain
$c=\{c_m\}$
is given by
\ben \langle \phi,c\rangle= \sum_{m=1}^M\phi_m(c_m).\een
This pairing satisfies
\ben \langle (b+B)\phi,c\rangle=\langle\phi,(b^T+B^T)c\rangle.\een

We have the relations
$$b^2=B^2=0=bB+Bb=(b+B)^2$$
so that we may define the
cyclic cohomology of $\A$ as the cohomology
of the total $(b,B)$-complex.

In this survey the main application we will make of
the Dixmier trace to cyclic cohomology is the discussion in Section 17 of
the formula of A. Connes for the Hochschild class of the Chern character.
Here we use the unrenormalised complex
$\tilde C^m=\A^{\otimes m+1}$. The definition of the
Hochschild coboundary $b$ on $\tilde C^m$ involves
the same formula~\eqref{900}.
The Hochschild cohomology, denoted
$HH^*(\A,\A^*)$, is then the cohomology of the complex
$(\tilde C^*(\A),b)$.

The Hochshild boundary on the complex $\tilde C_m=\A^{\otimes m+1}$
is the operator $b^T$.
We say that $c$ is a Hochschild cycle if $b^Tc=0$.
When the Hochschild homology
is well-defined we denote it by $HH_*(\A)$.

One can interpret spectral flow (in the type $I$ case) as the
pairing between an odd $K$-theory class represented by a unitary
$u$, and an odd $K$-homology class represented by $(\A,\HH,\D)$,
\cite[Chapter III,IV]{Co4}. This point of view also makes sense in
the general semifinite setting, though one must suitably interpret
$K$-homology, \cite{CPRS1,CP2}. A central feature of \cite{Co4} is
the translation of the $K$-theory pairing to cyclic theory in
order to obtain index theorems. One associates to a suitable
representative of a $K$-theory class, respectively a $K$-homology
class, a class in periodic cyclic homology, respectively a class
in periodic cyclic cohomology, called a Chern character in both
cases. (We will not digress here to discuss the periodic theory
and the periodicity operator referring instead to \cite{Co4} and
\cite{GVF}.) The principal result is then \be sf(\D,u\D
u^*)=\langle [u],[(\A,\HH,\D)]\rangle=-\frac{1}{\sqrt{2\pi i}}
\langle [Ch_*(u)],[Ch^*(\A,\HH,\D)]\rangle,\label{indpair}\ee
where $[u]\in K_1(\A)$ is a $K$-theory class with representative
$u$ and $[(\A,\HH,\D)]$ is the $K$-homology class of the spectral
triple $(\A,\HH,\D)$. On the right hand side, $Ch_*(u)$ is the
Chern character of $u$. We recall that the Chern character of a
unitary $u$ is the following (infinite) collection of odd chains
$Ch_{2j+1}(u)$ satisfying $bCh_{2j+3}(u)+BCh_{2j+1}(u)=0$, \ben
Ch_{2j+1}(u)=(-1)^jj!u^*\otimes u\otimes u^*\otimes\cdots\otimes
u\ \ \ (2j+2\ \ \mbox{entries}) .\een We have used the notation
$[Ch_*(u)]$ for the periodic cyclic homology class. Similarly
$[Ch^*(\A,\HH,\D)]$ is the periodic cyclic cohomology class of the
Chern character of $(\A,\HH,\D)$.

\section{The Dixmier trace and residues of the zeta function}\label{DTandRofZF}

Many applications of the Dixmier trace rely on being able to
calculate it by taking a residue of an associated zeta function.
The main result in this direction in the type $I$ case is
Proposition IV.2.4 of \cite{Co4}. Recent advances \cite{CPS2} have
extended this to the general semifinite von Neumann setting.

\subsection{Preliminaries}


First it is useful to have an estimate on the singular values
of the operators in $\clL ^{(1,\infty)}$.

\begin{lemma}\label{2.1}
For $T\in\clL ^{(1,\infty)}$ positive
there is a constant $K>0$ such that for each
$p\geq 1$,
$$\int_0^t\mu_s(T)^pds \leq K^p\int_0^t\frac{1}{(s+1)^p}ds.$$
\end{lemma}

\begin{proof}  By \cite[Lemma 2.5 (iv)]{FK},  for all $0\leq T\in \clN $ and all continuous increasing
functions $f$ on $[0,\infty)$ with $f(0)\ge 0$, we have
$\mu_s(f(T))=f(\mu_s(T))$ for all $s>0$. Combining this fact with
well-known result of Hardy-Littlewood-P\'olya (see e.g. \cite{F},
Lemma 4.1), we see that $T_1\prec\prec T_2$, $0\leq T_1,T_2\in
\clN $ implies $T_1^p\prec\prec T_2^p$ for all $p\in (1,\infty)$.
Now, by definition of $\clL ^{(1,\infty)}$ the singular values of
$T$ satisfy $\int_0^t\mu_s(T)ds=O(\log t)$ so that for some $K>0$,
$$\int_0^t\mu_s(T) ds \leq K \int_0^t\frac{1}{(s+1)}ds,\quad \forall t>0.$$
In other words $\mu_s(T)\prec\prec K/(1+s)$ and the assertion of lemma
follows immediately.
\end{proof}
In the next theorem we use dilation invariant states from
Theorem~\ref{DualityThm}.
\begin{thm}\label{2.2}(weak$^*$-Karamata theorem)
Let $\tilde\omega\in L_\infty(\mbR)^*$
be a dilation invariant state and
let $\beta$ be a real valued, increasing,
right continuous function on $\mbR_+$
which is zero at zero
and such that the integral
$h(r)=\int_0^\infty e^{-\frac{t}{r}}d\beta(t)$
converges for all $r>0$ and
$C=\tilde\omega\mbox{-}\lim_{r\to\infty} \frac{1}{r}h(r)$ exists.
Then
$$\tilde\omega\mbox{-}\lim_{r\to\infty} \frac{1}{r}h(r)
=\tilde\omega\mbox{-}\lim_{t\to\infty}\frac{\beta(t)}{t}.$$
\end{thm}

\begin{rems}
The classical Karamata theorem states, in the notation of the
theorem, that if the ordinary limit
$\lim_{r\to\infty} \frac{1}{r}h(r)=C$ exists
then $C=\lim_{t\to\infty}\frac{\beta(t)}{t}$.
The proof of this classical result is obtained by replacing,
in the proof of Theorem \ref{2.2},
$\tilde\omega\mbox{-}\lim$ throughout by the ordinary limit.
\end{rems}
\begin{proof}
Let
$$g(x)=\left\{\begin{array}{ll} x^{-1} & \mbox{for } e^{-1}\leq x\leq 1\\
                                 0       & \mbox{for } 0\leq x< e^{-1}
\end{array}\right.$$
so that $g$ is right continuous at $e^{-1}$.
Then for $r>0$, $t\to e^{-t/r}g(e^{-t/r})$ is left continuous
at $t=r$. Thus the Riemann-Stieltjes integral
$\int_0^\infty e^{-t/r}g(e^{-t/r})d\beta(t)$
exists for each $r>0$.
We claim that for any polynomial $p$
$$\tilde\omega\mbox{-}\lim_{r\to\infty}\frac{1}{r}\int_0^\infty
e^{-t/r}p(e^{-t/r})d\beta(t) =C\int_0^\infty
e^{-t}p(e^{-t})dt.$$ To see this first compute
for $p(x)=x^n$,
$$\frac{1}{r}\int_0^\infty e^{-t/r}e^{-nt/r}d\beta(t)
=\frac{1}{r}\int_0^\infty e^{-(n+1)t/r}d\beta(t).$$
Therefore
$$\frac{1}{n+1}\tilde\omega\mbox{-}\lim_{r\to\infty}\frac{1}{r/(n+1)}
\int_0^\infty e^{-(n+1)t/r}d\beta(t) =\frac{C}{n+1}$$
by dilation invariance of $\tilde\omega$.
Thus $$\tilde\omega\mbox{-}\lim_{r\to\infty}
\frac{1}{r}\int_0^\infty e^{-t/r}e^{-nt/r}d\beta(t)
=C\int_0^\infty
e^{-t}(e^{-t})^ndt.$$ Since $\tilde\omega$ is linear
the claim follows for all $p$.

Choose sequences of polynomials
$\{p_n\},\ \{P_n\}$ such that for
all $x\in [0,1]$
$$-1\leq p_n(x)\leq g(x) \leq P_n(x)\leq 3$$
and such that $p_n$ and $P_n$ converge a.e. to $g(x)$.
Then since $\tilde\omega$ is positive it preserves order:
$$C\int_0^\infty
e^{-t}p_n(e^{-t})dt=\tilde\omega\mbox{-}\lim_{r\to\infty}\frac{1}{r}\int_0^\infty
e^{-t/r}p_n(e^{-t/r})d\beta(t)\leq
\tilde\omega\mbox{-}\lim_{r\to\infty}\frac{1}{r}\int_0^\infty
e^{-t/r}g(e^{-t/r})d\beta(t)$$
$$\leq\ldots\leq C\int_0^\infty
e^{-t}P_n(e^{-t})dt.$$
By the Lebesgue Dominated Convergence Theorem both
$\int_0^\infty
e^{-t}p_n(e^{-t})dt$ and $\int_0^\infty
e^{-t}P_n(e^{-t})dt$ converge to $\int_0^\infty
e^{-t}g(e^{-t})dt$ as $n\to\infty$.
But a direct calculation yields $\int_0^\infty
e^{-t}g(e^{-t})dt=1$ and
$$\int_0^\infty e^{-t/r}g(e^{-t/r})d\beta(t)=\beta(r).$$
Hence
$$C=
\tilde\omega\mbox{-}\lim_{r\to\infty}\frac{1}{r}\int_0^\infty
e^{-t/r}g(e^{-t/r})d\beta(t)
=\tilde\omega\mbox{-}\lim_{r\to\infty}\frac{\beta(r)}{r}.$$
\end{proof}

\noindent Recall that for any $\tau$-measurable operator $T$, the distribution
function of $T$ is defined by
$$
\lambda _t(T):=\tau (\chi_{(t,\infty)}(|T|)),\quad t>0,
$$
where $\chi_{(t,\infty)}(|T|)$ is the spectral projection of $|T|$
corresponding to the interval $(t,\infty)$ (see [FK]). By Proposition 2.2
of [FK],
$$\mu_s(T)=\inf\{t\ge 0\ :\ \lambda_t(T)\leq s\}
$$
we infer that for any  $\tau$-measurable operator $T$, the distribution
function $\lambda _{(\cdot)}(T)$ coincides with the (classical) distribution
function of $\mu_{(\cdot)}(T)$. From this formula and the fact that $\lambda$
is right-continuous, we can easily see that for $t>0$, $s>0$
$$s\geq\lambda_t\Longleftrightarrow \mu_s\leq t.$$
Or equivalently,
$$s < \lambda_t \Longleftrightarrow \mu_s > t.$$
Using Remark 3.3 of [FK] this implies that:
\begin{equation} \int_0^{\lambda _t} \mu_s(T)ds=\int_{[0,\lambda _t)}
\mu_s(T)ds=
\tau (|T|\chi_{(t,\infty)}(|T|)),\quad t>0.\label{equation*}\end{equation}

\begin{lemma}
For $T\in\clL ^{(1,\infty)}$ and
$C>\Vert T\Vert _{(1,\infty)}$ we have eventually
$$
\lambda _{\frac{1}{t}}(T)\leq Ct\log t.
$$
\end{lemma}

\begin{proof} Suppose not and there exists $t_n\uparrow \infty$ such that
$\lambda _{\frac{1}{t_n}}(T)> Ct_n\log t_n$ and so for $s\leq Ct_n\log t_n$
we have $\mu_s(T)\geq\mu_{Ct_n\log t_n}(T) > \frac{1}{t_n}$. Then for
sufficiently
large $n$
$$
\int_0^{ Ct_n\log t_n}\mu_s(T)ds> \frac{1}{t_n}\cdot  Ct_n\log t_n=
C\log t_n.
$$
Choose $\delta>0$ with $C-\delta>\Vert T\Vert_{(1,\infty)}$.
  Then for sufficiently
large $n$
$$
C\log t_n=(C-\delta)\log t_n +\delta\log t_n>
\Vert T\Vert _{(1,\infty)}\log( Ct_n)
+\Vert T\Vert _{(1,\infty)}\log(\log (t_n+1))
$$
$$
=\Vert T\Vert _{{(1,\infty)}}\log( Ct_n\log (t_n+1)).$$
This is a contradiction with the inequality
$\int_0^{t}\mu_s(T)ds\leq \Vert T\Vert_{{(1,\infty)}}\log (t+1)$,
which holds for any $t>0$ due to the
definition of the norm in $\clL ^{(1,\infty)}$.
\end{proof}

\noindent An assertion somewhat similar to Proposition \ref{2.4}
  below was formulated
in \cite{Pr} and supplied with an incorrect proof. We use a different approach.
\bigskip

\begin{prop}\label{2.4}
   For $T\in\clL ^{(1,\infty)}$ positive let $\omega$ be a state on
   $L_\infty(\mbR^*_+)$ satisfying all the conditions of
   Corollary~\ref{DualityCor}. For every $C>0$
$$
\tau_\omega(T)=
\omega\mbox{-}\lim_{t\to\infty}\frac{1}{\log(1+t)}\int_0^t \mu_s(T) ds
=\omega\mbox{-}\lim_{t\to\infty}\frac{1}{\log(1+t)}
\tau(T\chi_{(\frac{1}{t},\infty)}(T))
$$
$$
=\omega\mbox{-}\lim_{t\to\infty}\frac{1}{\log(1+t)}\int_0^{Ct\log t}\mu_s(T) ds
$$
and if one of the $\omega\mbox{-}$limits is a true limit then so are the
others.

\end{prop}

\begin{proof} We first note that
$$
\int_0^t\mu_s(T)ds\leq \int_0^{\lambda _{\frac{1}{t}}(T)} \mu_s(T)ds+1,
\quad t>0.
$$
Indeed, the inequality above holds trivially if
$t\leq\lambda _{\frac{1}{t}}(T)$. If $t> \lambda _{\frac{1}{t}}(T)$, then
$$
\int_0^t\mu_s(T)ds=\int_0^{\lambda _{\frac{1}{t}}(T)} \mu_s(T)ds+\int_
{\lambda _{\frac{1}{t}}(T)}^t \mu_s(T)ds.
$$
Now $s> \lambda _{\frac{1}{t}}(T)$ implies that $\mu_s(T)\leq \frac{1}{t}$
so we have
$$
\int_0^t\mu_s(T)ds\leq \int_0^{\lambda _{\frac{1}{t}}(T)}
\mu_s(T)ds+\frac{1}{t}(t-\lambda _{\frac{1}{t}}(T))\leq \int_0^{\lambda
_{\frac{1}{t}}(T)} \mu_s(T)ds+1.
$$
Using this observation and lemma  above we see that for
$C>\Vert T\Vert _{{(1,\infty)}}$ and any fixed
$\alpha >1$  eventually
$$
\int_0^t\mu_s(T)ds \leq \int_0^{\lambda _{\frac{1}{t}}(T)} \mu_s(T)ds+1\leq
\int_0^{Ct\log t} \mu_s(T) ds +1 \leq\int_0^{t^\alpha}
\mu_s(T) ds+1
$$
and so eventually
$$
\frac{1}{\log(1+t)}\int_0^t\mu_s(T)ds
\leq\frac{1}{\log(1+t)}( \int_0^{\lambda _{\frac{1}{t}}(T)} \mu_s(T)ds+1)
\leq\frac{1}{\log(1+t)}(\int_0^{Ct\log t} \mu_s(T) ds +1)
$$
$$
\leq\frac{\log(1+t^\alpha)}{\log(1+t)\log(1+t^\alpha)}
(\int_0^{t^\alpha} \mu_s(T) ds+1).
$$
Taking the $\omega$-limit we get
$$
\tau_\omega(T) \leq\omega\mbox{-}\lim_{t\to\infty}\frac{1}{\log(1+t)}
\int_0^{\lambda _{\frac{1}{t}}(T)} \mu_s(T)ds
\leq\omega\mbox{-}\lim_{t\to\infty}\frac{1}{\log(1+t)}\int_0^{Ct\log t}
\mu_s(T) ds
$$
$$
\leq\omega\mbox{-}\lim_{t\to\infty}\frac{\alpha}{\log(1+t^\alpha)}
\int_0^{t^\alpha} \mu_s(T) ds=\alpha\tau_\omega(T)
$$
where the last line uses Corollary~\ref{DualityCor}~(4).
Since this holds for all $\alpha>1$ and using~\ref{equation*} we get the
conclusion
for $\omega$-limits and $C>\Vert T\Vert _{{(1,\infty)}}$.
The assertion for an arbitrary $0<C\leq \Vert T\Vert _{{(1,\infty)}}$
follows immediately by noting that for $C'> \Vert T\Vert _{{(1,\infty)}}$
one has eventually
$$
\int_0^t\mu_s(T)ds\leq\int_0^{Ct\log t} \mu_s(T) ds\leq\int_0^{C't\log t}
\mu_s(T) ds.
$$

To see the last assertion of the Proposition suppose that
$\lim_{t\to\infty}\frac{1}{\log(1+t)}\int_0^t\mu_s(T)ds=A$
then by the above argument we get
$$
A\leq\liminf_{t\to\infty}\frac{1}{\log(1+t)}
\tau(T\chi_{(\frac{1}{t},\infty)}(T))
\leq\limsup_{t\to\infty}\frac{1}{\log(1+t)}
\tau(T\chi_{(\frac{1}{t},\infty)}(T))
\leq\alpha A
$$
for all $\alpha>1$ and hence
$\lim_{t\to\infty}\frac{1}{\log(1+t)}
\tau(T\chi_{(\frac{1}{t},\infty)}(T))=A$.
On the other hand if the limit
$\lim_{t\to\infty}\frac{1}{\log(1+t)}\tau(T\chi_{(\frac{1}{t},\infty)}(T))$
exists and equals $B$ say then
$$
\limsup_{t\to\infty}\frac{1}{\log(1+t)}\int_0^t \mu_s(T) ds\leq B\leq
\alpha\liminf_{t\to\infty}\frac{1}{\log(1+t)}\int_0^t \mu_s(T) ds
$$
for all $\alpha>1$ and so
$$
\lim_{t\to\infty}\frac{1}{\log(1+t)}\int_0^t \mu_s(T) ds=B
$$
as well. The remaining claims follow similarly.
\end{proof}

\subsection{ The zeta function and the Dixmier trace}

This subsection is motivated by Proposition IV.2.4 of
\cite{Co4}. We will describe several generalisations of this result to the
von Neumann setting. Our approach will be analogous
though somewhat different to
that in \cite{Co4}.

The zeta function of a  positive $T\in \clL ^{(1,\infty)}$
is given by
$\zeta(s)=\tau(T^s)$
while for $A\in\mathcal N$ we set
$\zeta_A(s)=\tau(AT^s).$
We are interested in the asymptotic behaviour of $\zeta(s)$
and $\zeta_A(s)$ as $s\to 1$.

Now it is elementary to see that the discussion of singular
traces is relevant because
by Lemma \ref{2.1} we have for some $K>0$ and all $s>1$
  $$\tau(T^s) =\int_0^\infty \mu_r(T^s) dr=\int_0^\infty \mu_r(T)^s dr
$$
$$\leq \int_0^\infty \frac{K^s}{(1+r)^s} dr = \frac{K^s}{s-1}.$$
{}From this it follows that $\{(s-1)\tau(T^s)\vert \ s>1\}$
is bounded. Now for $A$ bounded
$|(s-1)\tau(AT^s)|\leq ||A||(s-1)\tau(T^s)$
so that $(s-1)\tau(AT^s)$ is also bounded
and hence for any
$\tilde\omega \in L_\infty(\mbR)^*$ satisfying conditions (1), (2) and (3)
of Theorem \ref{DualityThm}
\begin{equation}\tilde\omega\mbox{-}\lim_{r\to
\infty}\frac{1}{r}\tau(AT^{1+\frac{1}{r}})
\label{equation3.1}\end{equation}
  exists.

Here $r\to \frac{1}{r}\tau(AT^{1+\frac{1}{r}})$
is defined as a function on all of $\mbR$ by extending it to be identically
zero
for $r<1$. One might like to think of (\ref{equation3.1})
as  $\tilde\omega\mbox{-}\lim_{s\to 1}(s-1)\tau(AT^s)$ but this of course does
not (strictly speaking) make sense whereas if $\lim_{s\to
   1}(s-1)\tau(AT^s)$
exists then it is $\lim_{r\to\infty}\frac{1}{r}\tau(AT^{1+\frac{1}{r}})$.

In the following theorem we will take $T\in\clL ^{(1,\infty)}$ positive,
$||T||\leq 1$ with spectral resolution $T=\int \lambda dE(\lambda)$. We would
like to integrate with respect to $d\tau(E(\lambda))$; unfortunately, these
scalars $\tau(E(\lambda))$ are, in general,
  all infinite. To remedy this situation, we instead
must integrate with respect to the increasing (negative) real-valued function
$N_T(\lambda)=\tau(E(\lambda)-1)$ for $\lambda >0$. Away from $0$, the
increments $\tau(\triangle E(\lambda))$ and $\triangle N_T(\lambda)$ are, of
course, identical.

\begin{thm}\label{3.1}
   For $T\in\clL ^{(1,\infty)}$ positive, $||T||\leq 1$ and
   $\tilde\omega\in L_\infty(\mbR)^*$ satisfying all the conditions of
   Theorem~\ref{DualityThm}, let $\tilde\omega=\omega\circ L$ where $L$
   is given in section~\ref{PrelimDilInvState} (prior to
   Definition~\ref{familiesOfSelf-mapsDef}), then we have:
$$\tau_\omega(T)=\tilde\omega\mbox{-}\lim\frac{1}{r}\tau(T^{1+\frac{1}{r}}).$$
If $\lim_{r\to\infty}\frac{1}{r}\tau(T^{1+\frac{1}{r}})$ exists then
$$\tau_\omega(T)=\lim_{r\to\infty}\frac{1}{r}\tau(T^{1+\frac{1}{r}})$$
for an arbitrary dilation invariant functional
$\omega\in L_\infty(\mbR^*_+)^*$.
\end{thm}

\begin{proof}
By (\ref{equation3.1}) we can apply the weak$^*$-Karamata theorem
to $\frac{1}{r}\tau(T^{1+\frac{1}{r}})$. First write
$\tau(T^{1+\frac{1}{r}})=\int_{0^+}^1 \lambda^{1+\frac{1}{r}}
dN_T(\lambda)$. Thus setting $\lambda=e^{-u}$
$$\tau(T^{1+\frac{1}{r}})=\int_0^\infty e^{-\frac{u}{r}}d\beta(u)$$
where $\beta(u)=\int_u^0 e^{-v}dN_T(e^{-v})=-\int_0^u e^{-v}dN_T(e^{-v})$.
Since the change of variable $\lambda=e^{-u}$ is strictly decreasing, $\beta$
is, in fact, nonnegative and increasing.
By the weak$^*$-Karamata theorem
  applied to $\tilde\omega\in L_\infty(\mbR)^*$
$$\tilde\omega\mbox{-}\lim_{r\to \infty} \frac{1}{r}\tau(T^{1+\frac{1}{r}})
=\tilde\omega\mbox{-}\lim_{u\to\infty}
\frac{\beta(u)}{u}.$$

Next with the substitution $\rho=e^{-v}$ we get:
\begin{equation}\tilde\omega\mbox{-}\lim_{u\to\infty}\frac{\beta(u)}{u}=
\tilde\omega\mbox{-}\lim_{u\to\infty}
\frac{1}{u}\int^1_{e^{-u}}
\rho dN_T(\rho).  \label{equation3.2}\end{equation}

Set $f(u)=\frac{\beta(u)}{u}$.
We want to make the change of variable $u=\log t$
or in other words to consider $f\circ \log = Lf$.
We use the discussion in section~\ref{PrelimDilInvState}
which tells us that if we start with an $M$ invariant
functional $\omega\in L_\infty(\mbR_+^*)^*$
then the functional $\tilde\omega=\omega\circ L$
is $H$ invariant as required by the
theorem. Then we have
$$\tilde\omega\mbox{-}\lim_{r\to \infty} \frac{1}{r}\tau(T^{1+\frac{1}{r}})
=\tilde\omega\mbox{-}\lim_{u\to\infty}\frac{\beta(u)}{u}=\tilde\omega\mbox{-}\lim_{u\to\infty}
f(u)=\omega\mbox{-}\lim_{t\to\infty}Lf(t)=\omega\mbox{-}\lim_{t\to\infty}
\frac{1}{\log t}\int_{1/t}^1\lambda dN_T(\lambda).$$

Now, by Proposition \ref{2.4}
$$\omega\mbox{-}\lim_{t\to\infty}
\frac{1}{\log t}\int_{1/t}^1\lambda dN_T(\lambda)
=\omega\mbox{-}\lim_{t\to\infty} \frac{1}{\log
t}\tau(\chi_{(\frac{1}{t},1]}(T)T)
=\tau_\omega(T).$$
This completes the proof of the first part of the theorem.

The proof of the second part is similar. Using the classical
Karamata theorem (see the remark following the statement of
Theorem \ref{2.2}) we obtain the following analogue
of~(\ref{equation3.2}):
$$\lim_{r\to\infty}\frac{1}{r}\tau(T^{1+r})=\lim\frac{\beta(u)}{u}=
\lim_{u\to\infty}\frac{1}{u}\int_{e^{-u}}^1\rho dN_T(\rho).$$
Making the substitution $u=\log t$ on the right hand side we have
$$\lim_{u\to\infty}\frac{1}{u}\int_{e^{-u}}^1\rho dN_T(\rho)=
\lim_{t\to\infty}\frac{1}{\log t}\int_{\frac{1}{t}}^1\lambda dN_T(\lambda)
=\tau_\omega(T)$$ where in the last equality we need only
dilation invariance of the state $\omega\in L_\infty(\mbR^*_+)^*$
and not the full list of conditions of Corollary~\ref{DualityCor}
\end{proof}

The map on positive $T\in \clL ^{(1,\infty)}$ to $\mbR$
given by $T\to \tau_\omega(T)$ can be extended
by linearity to a $\mbC$ valued functional on
all of $\clL ^{(1,\infty)}$.
Then the functional
\begin{equation}
A\mapsto\tau_\omega(AT) \label{equation**}
\end{equation}
  for $A\in\mathcal N$ and fixed
$T\in\clL ^{(1,\infty)}$ is well defined. We intend to study the
properties of (\ref{equation**}). Part of the interest in this
functional stems from the following result (see \cite{CGS}) as
well as the use of the Dixmier trace in noncommutative geometry
\cite{Co4}.

\begin{lemma}\label{3.2} Let $T\in\clL ^{(1,\infty)}$, then\\
(i) For $A\in\mathcal N$  we have
$$\tau_\omega(AT)=\tau_\omega(TA).$$
(ii) Assume that $D_0$ is an unbounded self adjoint
operator affiliated with $\mathcal N$ such that\\
$T=(1+D_0^2)^{-1/2}\in\clL ^{(1,\infty)}$.
If $[A_j,|D_0|]$ is a bounded operator for $A_j\in\clN ,\ j=1,2$ then
$$\tau_\omega(A_1A_2T)=\tau_\omega(A_2A_1T).$$
\end{lemma}

\begin{proof}(i) This is proposition A.2 of \cite{CM}.
The proof is elementary, first show that
$\tau_\omega(UTU^*)=\tau_\omega(T)$
then use linearity to extend to arbitrary
$T\in\clL ^{(1,\infty)}$.
Replace $T$ by $TU$ then use linearity again.\\
(ii) We remark that  $[A_j,|D_0|]$ defining a bounded operator
means that the $A_j$ leave $dom(|D_0|)=dom(D_0)$ invariant and
that $[A_j,|D_0|]$ is bounded on this domain (see
\cite[Proposition 3.2.55]{BR} and its proof for equivalent but
seemingly weaker conditions). As $|D_0|-(1+D_0^2)^{1/2}$ is
bounded, $[A_j, (1+D_0^2)^{1/2}]$ defines a bounded operator
whenever $[A_j,|D_0|]$ does. As $T^{-1} = (1+D_0^2)^{1/2}$ and
$T:\mathcal H \to dom(T^{-1})$, we see that the formal
calculation:
$$[A_j,T]=A_jT-TA_j=T(T^{-1}A_j-A_jT^{-1})T=T[T^{-1},A_j]T$$
makes sense as an everywhere-defined operator on $\mathcal H$. That is,
$$[A_j,T]=T[(1+D_0^2)^{1/2},A_j]T\in(\clL ^{(1,\infty)})^2\subseteq
\clL ^1.$$ Then
we have, using part (i),
$$\tau_\omega(A_1A_2T)=\tau_\omega(A_2A_1T)-\tau_\omega([A_1,T]A_2).$$
So then
$$\tau_\omega(A_1A_2T)
=\tau_\omega(A_2A_1T)-\tau_\omega(T[(1+D_0^2)^{1/2},A_1]TA_2).$$
Since the operator in the last term is trace class we are done.
\end{proof}

We will consider spectral triples which involves choosing a
subalgebra of $\mathcal N$ on which~(\ref{equation**}) will define
a trace.
\begin{thm}\label{3.8}
Let $A\in\mathcal N$, $T\geq 0$, $T\in\clL ^{(1,\infty)}.$\\
(i) If $\lim_{s\to 1^+}(s-1)\tau(AT^s)$
exists then it is equal to $\tau_\omega(AT)$ where we choose $\omega$ as in
the proof of Theorem \ref{3.1}.\\
(ii) More generally, if we choose functionals $\omega$ and $\tilde\omega$ as in
the proof of Theorem \ref{3.1} then
$$\tilde\omega\mbox{-}\lim_{r\to\infty}
\frac{1}{r}\tau(AT^{1+\frac{1}{r}})=\tau_\omega(AT).$$
\end{thm}

\begin{proof} For part (i) we first assume that $A$ is self adjoint.
Write $A=a^+-a^-$ where $a^\pm$ are positive.
Choose $\tilde\omega$ as in the proof of
Theorem \ref{3.1}, then
\begin{eqnarray*}
\lim_{s\to 1^+}(s-1)\tau(AT^s)&=&
\tilde\omega\mbox{-}\lim_{r\to\infty}\frac{1}{r}
\tau(AT^{1+\frac{1}{r}})\\
&=&\tilde\omega\mbox{-}\lim_{r\to\infty}\frac{1}{r}\tau(a^+T^{1+\frac{1}{r}})
-\tilde\omega\mbox{-}\lim_{r\to\infty}\frac{1}{r}\tau(a^-T^{1+\frac{1}{r}})\\
&=& \tau_\omega(a^+T)-\tau_\omega(a^-T)\\
&=&\tau_\omega(AT).
\end{eqnarray*}
Here the third equality uses a technical result (Proposition 3.6)
from \cite{CPS2} and then Theorem \ref{3.1}.
The reduction from the general case to the self-adjoint case now follows
in a similar way.

For part (ii), we assume that $A$ is positive.
 From Lemma \ref{3.2}(i) and
Theorem \ref{3.1},
and Proposition 3.6 of \cite{CPS2} we have
\begin{eqnarray*}
\tau_\omega(AT)&=&\tau_\omega(A^{1/2}TA^{1/2})=\tilde\omega\mbox{-}
\lim_{r\to \infty}\frac{1}{r}\tau((A^{1/2}TA^{1/2})^{1+\frac{1}{r}})\\
&=&\tilde\omega\mbox{-}\lim_{r\to \infty}\frac{1}{r}\tau(AT^{1+\frac{1}{r}}).
\end{eqnarray*}
For general $A$ we reduce to the case $A$ positive as in the proof of part (i).
\end{proof}

\section{The heat semigroup formula}\label{HSFsec}

This Section is inspired by \cite{Co4} and again is taken from
\cite{CPS2}.
We retain the assumption  $T\geq 0$ throughout
and define $e^{-T^{-2}}$ as the operator that is zero on $\ker T$
and on $\ker T^\perp$ is defined by the functional
calculus. We remark that if $T\geq 0$, $T\in\clL ^{(p,\infty)}$
for some $p\geq 1$ then  $e^{-tT^{-2}}$ is trace class for all $t>0$.

Our aim in this Section is to prove the following

\begin{thm}\label{4.1}
If $A\in\mathcal N$, $T\geq 0$, $T\in\clL ^{(1,\infty)}$
then,
$$\omega\mbox{-}\lim_{\lambda\to\infty}\lambda^{-1}\tau(Ae^{-\lambda^{-2}T^{-2}})=
\Gamma(3/2)\tau_\omega(AT)$$
for $\omega\in L_\infty(\mbR^*_+)^*$ satisfying the conditions of
Corollary
\ref{DualityCor}.
\end{thm}

Let $\zeta_A(p+\frac{1}{r})=\tau(AT^{p+\frac{1}{r}})$.
Notice that $\frac{1}{2}\Gamma(\frac{p}{2})
\tilde\omega\mbox{-}\lim_{r\to\infty} \frac{1}{r}\zeta_A(p+\frac{1}{r})$
always exists. Hence we can reduce the hard part of the
proof of Theorem \ref{4.1} to the following preliminary result.

\begin{prop}\label{4.2}
If $A\in\mathcal N,A\geq 0$, $T\geq 0$, $T\in\clL ^{(p,\infty)}$, $1\leq p<\infty$
then, choosing $\omega$ and $\tilde\omega$
as in the proof of Theorem~\ref{3.1}, we have
$$\omega\mbox{-}\lim_{\lambda\to\infty} \frac{1}{\lambda}\tau(Ae^{-T^{-2}\lambda^{-2/p}})
=\frac{1}{2}\Gamma(\frac{p}{2})
\tilde\omega\mbox{-}\lim_{r\to\infty} \frac{1}{r}\zeta_A(p+\frac{1}{r}).$$
\end{prop}

\begin{proof}
We have, using the Laplace transform,
$$T^s = \frac{1}{\Gamma(s/2)}\int_0^\infty t^{s/2 -1}e^{-tT^{-2}}dt.$$
Then
$$\zeta_A(s)=\tau(AT^s)=
\frac{1}{\Gamma(s/2)}\int_0^\infty t^{s/2 -1}\tau(Ae^{-tT^{-2}})dt.$$
Make the change of variable $t=1/\lambda^{2/p}$ so that the preceding
formula becomes
$$\frac{p}{2}\Gamma(s/2)\zeta_A(s)=
\int_0^\infty \lambda^{-\frac{s}{p}-1}
\tau(Ae^{-\lambda^{-2/p}T^{-2}})d\lambda.$$
We split this integral into two parts,
$\int_0^1$ and $\int_1^\infty$
and call the first integral $R(r)$ where $s=p+\frac{1}{r}$.
Then
$$R(r)=\int_0^1 \lambda^{-\frac{1}{pr}-2}
\tau(Ae^{-\lambda^{-2/p}T^{-2}})d\lambda=\int_1^\infty
t^{\frac{p}{2}+\frac{1}{2r} -1}\tau(Ae^{-tT^{-2}})dt.$$
The integrand decays exponentially in $t$ as $t\to\infty$
because
$T^{-2}\geq {\Vert T^2\Vert}^{-1}\bf 1$
so that
$$\tau(Ae^{-tT^{-2}})\leq \tau(Ae^{-T^{-2}}e^{-\frac{t-1}{\Vert T^{2}\Vert}}).
$$
Then we can conclude that $R(r)$ is bounded independently of $r$ and so
$\lim_{r\to\infty}\frac{1}{r}R(r)=0$.
For the other integral
the change of variable
$\lambda = e^\mu$ gives
$$\int_1^\infty \lambda^{-\frac{1}{pr}-2}
\tau(Ae^{-\lambda^{-2/p}T^{-2}})d\lambda
=\int_0^\infty e^{-\frac{\mu}{pr}} d\beta(\mu)$$
where
$\beta(\mu)=\int_0^\mu e^{-v}\tau(Ae^{-e^{-\frac{2}{p}v}T^{-2}})dv$.
Hence we can now write
$$\frac{p}{2}\Gamma((p+\frac{1}{r})/2)\zeta_A(p+\frac{1}{r})
=\int_0^\infty e^{-\frac{\mu}{pr}} d\beta(\mu)+R(r).$$
Now consider
$$\frac{p}{2}\tilde\omega\mbox{-}\lim_{r\to\infty}\frac{1}{r}
\Gamma(\frac{p}{2}+\frac{p}{2r})\zeta_A(p+\frac{1}{r})=
\frac{p}{2}
\Gamma(p/2)\tilde\omega\mbox{-}\lim_{r\to\infty}\frac{1}{r}\zeta_A(p+\frac{1}{r}).$$
Then
$$\frac{p}{2}
\Gamma(p/2)\tilde\omega\mbox{-}\lim_{r\to\infty}\frac{1}{r}\zeta_A(p+\frac{1}{r})
=p\tilde\omega\mbox{-}\lim_{r\to\infty}\frac{1}{pr}
\int_0^\infty e^{-\mu/pr}d\beta(\mu)$$
(remembering that the term $\frac{1}{r}R(r)$ has limit
zero as $r\to \infty$).
By dilation invariance and
Theorem~\ref{2.2} we then have
\begin{equation}\frac{p}{2}
\Gamma(p/2)\tilde\omega\mbox{-}\lim_{r\to\infty}\frac{1}{r}
\zeta_A(p+\frac{1}{r})=p\tilde\omega\mbox{-}\lim_{\mu\to\infty}
  \frac{\beta(\mu)}{\mu}. \label{equation4.1}\end{equation}
Making the change of variable $\lambda = e^v$
in the expression for $\beta(\mu)$
we get
$$\frac{\beta(\mu)}{\mu}=\frac{1}{\mu}\int_1^{e{^\mu}}
\lambda^{-2}\tau(Ae^{-T^{-2}\lambda^{-2/p}})d\lambda$$
Make the substitution $\mu =\log t$ so the RHS becomes
$$\frac{1}{\log
t}\int_1^t\lambda^{-2}\tau(Ae^{-T^{-2}\lambda^{-2/p}})d\lambda=g_1(t)$$
This is the Cesaro mean of
$$g_2(\lambda)=\frac{1}{\lambda}\tau(Ae^{-T^{-2}\lambda^{-2/p}}).$$
So as we chose $\omega\in L_\infty(\mbR^*_+)^*$ satisfying
Corollary~\ref{DualityCor}, we have $\omega(g_1)=\omega(g_2)$.
Recalling that we choose $\tilde\omega$ to be related to $\omega$
as in Theorem \ref{3.1} and so using (\ref{equation4.1}) we obtain
$$\omega\mbox{-}\lim_{\lambda\to\infty} \frac{1}{\lambda}\tau(Ae^{-T^{-2}\lambda^{-2/p}})
=\frac{1}{2}\Gamma(\frac{p}{2})
\tilde\omega\mbox{-}\lim_{r\to\infty} \frac{1}{r}\zeta_A(p+\frac{1}{r}).$$
\end{proof}

To prove the theorem consider first the
case where $A$ is bounded, $A\geq 0$ and use
the Proposition \ref{4.2} and Theorem \ref{3.8} to assert that
$$\Gamma(3/2)\tau_\omega(AT)=
\Gamma(3/2)\tilde\omega\mbox{-}\lim_{r\to\infty}\frac{1}{r}\tau(AT^{1+\frac{1}{r}})
=\omega\mbox{-}\lim_{\lambda\to\infty}\lambda^{-1}
\tau(Ae^{-\lambda^{-2}T^{-2}}).
$$
Then for self adjoint $A$
write $A=a^+-a^-$ where $a^\pm$ are positive so that
$$\Gamma(3/2)\tau_\omega(AT)=
\Gamma(3/2)(\tau_\omega(a^+T)-\tau_\omega(a^-T))$$
$$=
\omega\mbox{-}\lim_{\lambda\to\infty}\lambda^{-1}
\tau(a^+e^{-\lambda^{-2}T^{-2}})-
\omega\mbox{-}\lim_{\lambda\to\infty}\lambda^{-1}\tau(a^-e^{-\lambda^{-2}T^{-2}})
$$
$$=
\omega\mbox{-}\lim_{\lambda\to\infty}\lambda^{-1}\tau(Ae^{-\lambda^{-2}T^{-2}}).$$
We can extend to general bounded $A$ by a similar argument.

\section{The case of $p>1$}

The results of the previous sections have analogues for $p>1$.
Some arguments are simpler
due to the fact that the singular values satisfy
for $T\in\clL ^{(p,\infty)}$ for $p>1$, $T\geq 0$
the asymptotic estimate $\mu_s(T)=O(\frac{1}{s^{1/p}})$.
Moreover $\tau(T^{p+\frac{1}{r}}) =\int_0^1 \lambda^{p+1/r} dN_T(\lambda)$
where $N_T(\lambda)=\tau(E(\lambda)-1)$ for $\lambda>0$ where
$T=\int \lambda dE(\lambda)$ is the spectral resolution for $T$.

We now establish some  $\clL ^{(p,\infty)}$
versions of our previous results.

\begin{lemma} \label{5.1} For $T\in\clL ^{(p,\infty)}$ and
$\omega$ and $\tilde\omega$ as in the proof of theorem \ref{3.1}
we have
$$p\tau_\omega(T^p)=\tilde\omega\mbox{-}\lim_{r\to\infty}
\frac{1}{r} \tau(T^{p+\frac{1}{r}}).$$
\end{lemma}

\begin{proof}
Set $\lambda=e^{-u/p}$ so that
$$\frac{1}{r} \tau(T^{p+\frac{1}{r}}) =p\frac{1}{pr}\int_0^\infty
e^{-u/rp}d\beta(u)$$
where $\beta(u)=\int_0^ue^{-v}dN_T(e^{-v/p})$.
So using dilation invariance:
$$\tilde\omega\mbox{-}\lim_{r\to\infty}\frac{1}{r} \tau(T^{p+\frac{1}{r}})
=p\tilde\omega\mbox{-}\lim_{r\to\infty}
\frac{1}{pr}\int_0^1
e^{-u/pr}d\beta(u)=p\tilde\omega\mbox{-}\lim_{u\to\infty}
\frac{\beta(u)}{u}$$
by the weak$^*$-Karamata Theorem~\ref{2.2}.
Reasoning as in the proof of Theorem \ref{3.1}  and substituting
$\lambda=e^{-v/p}$ and
$u=\log t$ we have
$$\tilde\omega\mbox{-}\lim_{u\to\infty}
\frac{\beta(u)}{u}=\omega\mbox{-}\lim_{t\to\infty}\frac{1}{\log{t}}
\int_{t^{-1/p}}^1 \lambda^p dN_T(\lambda)$$
$$=\omega\mbox{-}\lim_{t\to\infty}\frac{1}{\log{t}}\tau(\chi_{(\frac{1}{t},\infty)}
(T^p)T^p)=\tau_\omega(T^p).$$
\end{proof}

\begin{cor}
Let $T\geq 0$, $T\in\clL ^{(p,\infty)}$ then
$$\omega\mbox{-}\lim_{\lambda\to\infty} \frac{1}{\lambda}\tau(e^{-T^{-2}\lambda^{-2/p}})
=\Gamma(1+\frac{p}{2})
\tilde\omega\mbox{-}\lim_{r\to\infty} \frac{1}{r}\tau(T^{p+\frac{1}{r}}).$$
\end{cor}

\begin{proof}
Combine Proposition \ref{4.2} and Lemma \ref{5.1}.
\end{proof}

The $\clL ^{(p,\infty)}$ version of
Theorem \ref{3.8} and the following result of Connes'
are proved by following the same methods as for $p=1$.
One needs a number of straightforward extensions
of various technical results as described
in Section 5 of \cite{CPS2}.

\begin{thm}\label{5.3}
If $A$ is bounded, $T\geq 0$, $T\in\clL ^{(p,\infty)}$
for $p\geq 1$
$$\omega\mbox{-}\lim_{\lambda\to\infty}\lambda^{-1}\tau(Ae^{-\lambda^{-2/p}T^{-2}})=
\Gamma(1+p/2)\tau_\omega(AT^p).$$
\end{thm}

Finally it is now straightforward to obtain the following result.

\begin{thm}\label{5.6}
If $A$ is bounded, $T\geq 0$, $T\in\clL ^{(p,\infty)}$
and
$$\lim_{s\to p^+}(s-p)\tau(AT^s)$$
exists then it is equal to $p\tau_\omega(AT^p)$.
\end{thm}

\section{Generalised Toeplitz operators and their index}

The index theory of generalised Toeplitz operators has a long
history which may be traced in part from \cite{CDSS} and \cite{Le91JOT}.
We will explain how the results established in earlier sections may
be used to contribute to this theory.

As before $P$ denotes the projection onto the nonnegative spectral
subspace of an unbounded self adjoint operator $\D_0$ (with
bounded inverse) and we are interested in the Breuer-Fredholm
index of the operator $PuP$ acting on $P\mathcal H$. When $\clN
=\clL (\clH )$ and $\D_0$ is part of an $\clL ^{(1,\infty)}$
summable spectral triple  the most general such theorem in this
type $I$ situation is that due to \cite{CM} who show show that
$$ind(PuP)=\frac{1}{2}\tau_\omega(u[\D_0,u^*]|\D_0|^{-1}).$$ We
now show that this formula holds when $\mathcal N$ is a general
semifinite von Neumann algebra.

\begin{rems*} Of course \cite{CM} consider a much more general situation which
would apply for example in the case where $\D_0$ is part of an
$\clL ^{(p,\infty)}$ summable spectral triple. The authors obtain
a very general formula for the index of $PuP$ in terms of sums of
residues of zeta functions constructed from $\D_0$ under certain
assumptions on the analytic properties of these zeta functions.
Some, but not all, of the zeta functions in their formula are
Dixmier traces. We will not attempt to describe this result here
as there is an excellent overview in the article \cite{Hi}. The
theorem in \cite{CM} is proved for the type $I$ case and its
generalisation to the semifinite von Neumann algebra setting is
achieved in \cite{CPRS1, CPRS2}.
\end{rems*}

We need a preliminary result, \cite{CPS2} Lemma 6.1.

\begin{lemma}\label{6.1}
Let $D_0$ be an unbounded self-adjoint operator affiliated with $\mathcal N$
so that $(1+D_0^2)^{-1/2}$ is in
$\cl^{(1,\infty)}$. Let $A_t$ and $B$ be in $\mathcal N$ for $t\in[0,1]$
with $A_t$ self-adjoint and $t\mapsto A_t$ continuous. Let $D_t=D_0+A_t$ and
let $p$ be a real number with $1<p<4/3$.
Then, the quantity $$\tau\left (B(1+D_0^2)^{-p/2} - B(1+D_t^2)^{-p/2}\right )$$
is uniformly bounded independent of $t\in [0,1]$ and $p\in (1,4/3)$.
\end{lemma}

\begin{thm}\label{GTOind}
Let $(\clN ,\D_0)$ be a
$\clL ^{(1,\infty)}$-summable Breuer-Fredholm module for the unital
Banach $*$-algebra, $\mathcal A$, and let $u\in \mathcal A$
  be a unitary such that $[\D_0,u]$
is bounded. Let $P$ be the projection on the non-negative spectral subspace of
$\D_0$. Then with $\omega$ chosen as in Corollary \ref{DualityCor},
\begin{eqnarray*}
ind(PuP)=sf(\D_0,u\D_0u^*) &=& \lim_{p\to
1^+}\frac{1}{2}(p-1)\tau(u[\D_0,u^*](1+\D_0^2)^{-p/2})\\
& = & \frac{1}{2}\tau_{\omega}(u[\D_0,u^*](1+\D_0^2)^{-1/2})\\
& = & \frac{1}{2}\tau_{\omega}(u[\D_0,u^*]|\D_0|^{-1})
\end{eqnarray*}
where the last equality only holds if $\D_0$ has a bounded inverse.
\end{thm}

\begin{rems}
(1) The equality
\begin{equation}
ind(PuP)=\frac{1}{2}\tau_{\omega}(u[\D_0,u^*]|\D_0|^{-1})
\label{equation6.1}\end{equation} proved above should be compared
with Theorem IV.2.8 of \cite{Co4}. In the case where $\clN =\clL
(\clH )$ the RHS of (\ref{equation6.1}) is a Hochschild
$1-$cocycle on $\mathcal A$ which is known to equal the Chern
character of the $\clL ^{(1,\infty)}$-summable Fredholm module
$(\clA ,\D_0,\mathcal H)$.

(2) Since any $1$-summable module is clearly a $\clL ^{(1,\infty)}$
-summable module, the theorem implies that any unbounded $1$-summable module
must have a trivial pairing with $K_1(\mathcal A)$ and is therefore
uninteresting from the homological point of view.
\end{rems}
\begin{proof}
Let $\D_t^u=\D_0+tu[\D_0,u^*]$, $t\in[0,1]$ then by equation~(2) of
Section~2.1 (of~\cite{CP2}) we have for each $p>1$ that
$$
ind(PuP)=\frac{1}{C_{p/2}}\int_0^1\tau\left(u[\D_0,u^*]
(1+(\D_t^u)^2)^{-p/2}\right)dt.
$$
Now, by Lemma \ref{6.1}, we have that
$$\left|\tau\left(u[\D_0,u^*]
\left [(1+(\D_t^u)^2)^{-p/2}-(1+\D_0^2)^{-p/2}\right]\right)\right|$$
is uniformly bounded independent of $t$ and $p$ for $1<p<4/3$. Since,
$\tilde C_{p/2} \to\infty$ as $p\to 1^+$, we see that:
\begin{eqnarray*}
&&\left|ind(PuP) - \frac{1}{C_{p/2}}\tau\left(u[\D_0,u^*](1+\D_0^2)^{-p/2}
\right)\right|\\
&=&\left|\frac{1}{C_{p/2}}\int_0^1\tau\left(u[\D_0,u^*]
(1+(\D_t^u)^2)^{-p/2}\right)dt-
\frac{1}{C_{p/2}}\int_0^1\tau\left(u[\D_0,u^*](1+\D_0^2)^{-p/2}\right)
dt\right|\\
&\leq& \frac{1}{C_{p/2}}\int_0^1\left|\tau\left(u[\D_0,u^*]
\left [(1+(\D_t^u)^2)^{-p/2}-(1+\D_0^2)^{-p/2}\right]\right)\right|dt\\
&\leq& Constant/C_{p/2} \to 0.
\end{eqnarray*}
Now, it is elementary that as $p\to 1^+$
$$\frac{2}{p-1}=\int_{|x|\geq 1}\left(\frac{1}{|x|}\right)^p dx \sim
\int_{-\infty}^{\infty}\left(\frac{1}{\sqrt{1+x^2}}\right)^p dx =
\tilde C_{p/2}.$$
This ends the proof of the first equality
while the second equality follows from Theorem \ref{3.8}(i).

The third equality follows from
$$\left(\sqrt{1+\D_0^2}\right)^{-1} -|\D_0|^{-1} =
\left(\sqrt{1+\D_0^2}\right)^{-1} |\D_0|^{-1}\left(\sqrt{1+\D_0^2} +
|\D_0|\right)^{-1}.$$
\end{proof}

\section{Non-smooth foliations and pseudo-differential operators}

In this Section we give a further application of the earlier results
following  \cite{BeF} and \cite{Pr}.
The main aim of Prinzis' thesis \cite{Pr} was to establish a
Wodzicki residue formula for the Dixmier trace
of certain pseudo-differential operators associated
to actions of $\mbR^n$ on a compact space $X$.
This was greatly generalised in \cite{BeF},
however to cover that paper in detail would require
an extensive discussion of foliations, thus
we restrict ourselves here to a synopsis of the simple case in \cite{Pr}.

The set-up is the group-measure space construction of Murray-von Neumann.
Thus $X$ is a compact space equipped with a probability measure $\nu$
and a continuous free minimal ergodic action $\alpha$ of $\mbR^n$ on
$X$ leaving $\nu$ invariant. We write the action
as $x\to t.x$ for $x\in X$ and $t\in \mbR^n$. Then the crossed product
$L_\infty(X,\nu)\times_\alpha\mbR^n$ is a type $II$ factor
contained in the bounded operators on $L^2(\mbR^n, L^2(X,\nu))$. We describe
the construction.
For a function $f\in L^1(\mbR, L_\infty(X,\nu))
\subset L_\infty(X,\nu)\times_\alpha\mbR^n$
the action of $f$ on a vector $\xi$
in $L^2(\mbR^n, L^2(X,\nu))$ is defined
by twisted left convolution as follows:
$$(\tilde{\pi}(f)\xi)(s)
= \int_{\mbR^n}\alpha_s^{-1}(f(t))\xi(s-t)dt.$$
Here $f(t)$ is a function on $X$ acting as a multiplication operator on
$ L^2(X,\nu)$.
The twisted convolution algebra
$$L^1(\mbR^n, L_\infty(X,\nu)) \cap L^2(\mbR^n, L^2(X,\nu))$$ is
a dense subspace of $L^2(\mbR^n, L^2(X,\nu))$
and there is a canonical faithful, normal, semifinite trace, $Tr$,
on the von Neumann algebra that it generates. This von Neumann
algebra is $$\mathcal N=
(\tilde{\pi}(L_\infty(X,\nu)\times_{\alpha}\mbR^n))^{\prime\prime}.$$
For functions $f,g: \mbR^n\to L_\infty(X)$ which are in
$L^2(\mbR^n, L^2(X,\nu))$ and whose twisted left convolutions $\tilde{\pi}(f),
\tilde{\pi}(g)$
define bounded operators on $L^2(\mbR^n, L^2(X,\nu))$, this trace is given by:
$$Tr(\tilde{\pi}(f)^*\tilde{\pi}(g))=
\int_{\mbR^n} \int_Xf(t,x)g(t,x)^*d\nu(x)dt$$
where we think of $f,g$ as functions on $\mbR^n\times X$.

Identify
$L^2(\mbR^n)$ with $L^2(\mbR^n)\otimes 1\subset L^2(\mbR, L^2(X,\nu))$
then any scalar-valued
function $f$ on $\mbR^n$ which is the Fourier transform $f=\widehat{g}$ of a
bounded
$L^2$ function, $g$ will satisfy $f\in L^2(\mbR, L^2(X,\nu))$
and $\tilde{\pi}(f)$ will be a bounded operator.

Pseudo-differential operators are defined in terms of their
symbols. A smooth symbol of order $m$ is a function
$a:X\times\mbR^n\to\mbC$ such that for each $x\in X$
$a_x$, defined by $a_x(t, \xi) = a(t.x, \xi)$,
satisfies

\noindent(1) $\sup\{|\partial_\xi^\beta\partial_t^\gamma
a_x(t,\xi)(1+|\xi|)^{-m+|\beta|}
\ |\ (t,\xi)\in \mbR^n\times\mbR^n, \beta,\gamma\in \mbN^n, |\beta|+|\gamma|
\leq M \}<\infty$
for all $M\in\mbN$\\
(2) $\xi\to a_x(0,\xi)$ is a smooth function on $\mbR^n$ into the
space  $\clC ^\infty(X)$, the set of continuous functions $f$ on
$X$ such that $t \to(x\to f(t.x))$ is smooth on $\mbR^n$.\\

Each symbol $a$ defines a
pseudo-differential operator $Op(a)$ on $C(X)\otimes C_c^\infty(\mbR^n)$
(where $C_c^\infty(\mbR^n)$ is the space of smooth functions of compact support)
by
$$Op(a)f(x,t)=\frac{1}{(2\pi)^n}\int_{\mbR^n}e^{it\xi}
a(t.x,\xi)\hat f(x,\xi)d\xi,\ \ f\in  C(X)\otimes
C_c^\infty(\mbR^n).$$ The principal symbol of a pseudo-differential
operator $A$ on $X$ is the limit
$$\sigma_m(A)(x,\xi)=\lim_{\lambda\to\infty}\frac{a(x,\lambda\xi)}{\lambda^m},
\quad \quad (x,\xi)\in X\times\mbR^n\backslash\{0\} $$ if it
exists. We say $A$ is elliptic if its symbol $a$ is such that
$a_x$ is elliptic (that is invertible for all sufficiently large
$\xi$) for all $x\in X$. Prinzis studies invertible positive
elliptic pseudo-differential operators $A$ with a principal
symbol. Henceforth we will only consider such operators. The zeta
function of such an operator is $\zeta(z)=\tau(A^z)$ and this
exists because $A^z$ is in the trace class in $\mathcal N$
\cite{Pr} for $\Re z<-n/m$. Prinzis shows that
\begin{equation}\lim_{x\to-\frac{n}{m}^-}(x+\frac{n}{m})\zeta(x)=
-\frac{1}{(2\pi)^nm}\int_{X\times S^{n-1}}
\sigma_m(A)(x,\xi)^{-\frac{n}{m}} d\nu(x)d\xi
\label{equation7.1}\end{equation} and
that $A^{-\frac{n}{m}}\in \clL ^{(1,\infty)}(\mathcal N,\tau)$.

Now note that (\ref{equation7.1}) combined with Theorem \ref{5.6}
implies that we have the relation
$$\tau_\omega(A^{-\frac{n}{m}})=
\frac{1}{(2\pi)^nn}\int_{X\times S^{n-1}} \sigma_m(A)(x,\xi)^{-\frac{n}{m}}
d\nu(x)d\xi.$$
In other words we have a type $II$ Wodzicki residue
for evaluating the Dixmier trace of these
pseudo-differential operators.

Now, \cite{BeF} considers the case of measured foliations in the
sense of \cite{Co9}. Thus the leaves of the foliation are no
longer given by an action of $\mbR^n$ as above but are more general
submanifolds. The main result of \cite{BeF} on this topic is a
formula for the Dixmier trace of certain pseudo-differential
operators in terms of a local residue where the latter is a
generalisation of the Wodzicki residue to the foliation setting.

\section{The algebra of almost periodic functions.}

This section is adapted from \cite{BCPRSW}.

\subsection{Almost periodic pseudodifferential operators.}\label{ShFr}

We review Shubin's \cite{Shu} study of the index theory of
differential operators with almost periodic coefficients which
extends ideas of \cite{CDSS}.
  Recall that a trigonometric function is a finite linear combination of
exponential functions $e_\xi:x\mapsto e^{i<x,\xi>}$. The space $\Trig
(\R^n)$ of trigonometric functions is a $*-$subalgebra of the $C^*-$algebra
$C_b(\R^n)$ of continuous bounded functions. The uniform closure of  $\Trig
(\R^n)$ is
called the algebra of almost periodic functions and denoted $\CAP(\R^n)$.
This $C^*-$algebra is isomorphic to the algebra of
continuous functions on $\R^n_B$,
the Bohr compactification of $\R^n$.
  Addition
in $\R^n$ extends to $\R^n_B$ which is a compact abelian group containing
$\R^n$ as a dense subgroup. The
normalized Haar measure $\alpha_B$ on $\R^n_B$
is such that the family $(e_\xi)_{\xi\in \R^n}$ is orthonormal.
The measure $\alpha_B$ is given for any almost periodic function $f$ on
$\R^n$ by:
$$
\alpha_B (f) := \lim_{T\to +\infty} \frac{1}{(2T)^n} \int_{(-T,T)^n} f(x) dx.
$$
Using $\alpha_B$ one defines the Hilbert space completion $L^2(\R^n_B)$ of
$\Trig (\R^n)$. This Hilbert space  has an orthonormal basis given by
$(e_\xi)_{\xi\in \R^n}$. The Fourier transform $\maF_B:\ell^2(\R^n_d)
\longrightarrow L^2(\R^n_B)$ is given by:
$$
\maF_B (\delta_\xi) = e_\xi, \quad \text{ with }\delta_\xi (\eta) =
\delta_{\xi, \eta},
$$
where $\delta_{\xi, \eta}$ is the Kronecker symbol. We shall denote by
$\maF$ the usual Fourier transform on the abelian group $\R^n$ with its
usual Lebesgue measure.

The action of $\R^n$ on $\R^n_B$ by translations yields a topological
dynamical system whose naturally associated  von Neumann algebra is the
crossed product
$L_\infty (\R^n_B)\times \R^n$.
It is more convenient for applications to
consider the commutant of this von Neumann algebra
  denoting it by $\cnn$. It is also a  crossed product,
namely the von Neumann algebra $L_\infty (\R^n)\times \R^n_d$.
Then $\cnn$ is a type $II_\infty$ factor with a faithful normal
semifinite trace $\tau$. It can be described as the set of Borel
essentially bounded families $(A_\mu)_{\mu\in\R^n_B}$ of bounded
operators in $L^2(\R^n)$ which are $\R^n-$equivariant, i.e. such
that
$$
A_\mu = \sigma_\mu (A_0)=T_{-\mu} A_0T_{\mu}, \quad \forall \mu\in \R^n.
$$
Here and in the sequel we denote by $\sigma_\mu$ conjugation
of any operator with the translation $T_\mu$ so that
$
\sigma_\mu (B) = T_{-\mu} B T_\mu.
$
If we denote by $M_\varphi$ the operator of multiplication by a
bounded  function $\varphi$, then examples of such families are given for
any $\lambda$ by the families
$$
(\sigma_\mu(M_{e_\lambda}))_{\mu\in \R^n_B}.
$$
We choose the Fourier transform
$$\clF f(\zeta)=\frac{1}{(2\pi)^{n/2}}\int_{\R^n}e^{ix\zeta}f(x)dx.$$
Then the von Neumann algebra
$\cnn$ can be defined \cite{CDSS} as the completion
in the Hilbert space $\maH=L^2(\R^n)\otimes
L^2(\R^n_B)$ of the set of operators $\{ M_{e_\lambda} \otimes
M_{e_\lambda} , T_\lambda \otimes 1 \}$ when $\lambda$ ranges over $\R^n$.

There is a natural way to embed
the $C^*-$algebra $\CAP (\R^n)$ in $\cnn$ by setting
$$
\pi(f) := (\sigma_\mu (M_f))_{\mu\in \R^n_B}
$$
This family then
belongs to $\cnn$ and $\pi$ is clearly faithful.
Viewed as an operator on $\maH$, $\pi(f)$ is given by $\pi(f) (g) (x,\mu) =
f(x+\mu) g(x,\mu)$.
If $B=(B_\mu)_\mu$ is a positive  element of $\cnn$, then  we define the
expectation $E(B)$ as the Haar integral:
$$
E(B) := \int_{\R^n_B} B_\mu d\alpha_B(\mu).
$$
Since the family $B$ and $\alpha_B$ are translation invariant, the
operator $E(B)$ clearly commutes with the translations in
$L^2(\R^n)$ and is therefore given by a Fourier multiplier $\wM
(\varphi_B)$ with $\varphi_B$ a positive element of $L^{\infty}
(\R^n)$. Recall that the Fourier multiplier $\wM (\varphi_B)$ is
conjugation of the multiplication operator $M_\varphi$ by the
Fourier transform, i.e. $ \wM (\varphi_B) = \maF^{-1} M_\varphi
\maF.$ When the function $\varphi$ is in the Schwartz space, the
operator $\wM (\varphi_B)$ is convolution by the Schwartz function
$\frac{1}{(2\pi)^{n/2}}\maF^{-1} \varphi$.
Hence the expectation $E$ takes values in the von Neumann algebra $\wM
(L^{\infty}(\R^n))$, i.e.
$$
E: \cnn \longrightarrow \wM (L^{\infty}(\R^n)).
$$
Now, using the usual
Lebesgue integral on $\R^n$, and the normalisation
of Coburn et al\cite{CDSS}
we introduce the following definition of the trace $\tau$:
$$
\tau (B) = \int_{\R^n} \varphi_B (\zeta) d\zeta.
$$
\begin{lemma}\cite{CDSS}\
The map $\tau$ is, up to constant, the unique positive normal faithful
semifinite trace on $\cnn$.
\end{lemma}

Consider the trace on $\cnn$ evaluated on an
operator of the form $M_aK$ where $a$ is almost periodic and
$K$ is a convolution operator on
$L^2(\R^n)$ arising from multiplication by an $L^1$ function
$k$ on the Fourier transform.
We have,
$$\tau(M_aK)
=  \lim_{T\to +\infty} \frac{1}{(2T)^n} \int_{(-T, T)^n}
a(x) dx\int_{\R^n}k(\zeta)d\zeta.$$
More generally, any pseudodifferential operator $A$
on $L^2(R^n,\C^N)$
with almost periodic coefficients of nonpositive order $m$ acting on
$\C^N-$valued functions, can be viewed as a family over $\R_{B}^n$ of
pseudodifferential operators on $\R^n$. To do this first take
the symbol  $a$ of $A$, then the operator $\sigma_\mu(A)$ is the
pseudodifferential operator with almost periodic coefficients whose symbol is
$$
(x,\xi) \longmapsto a(x+\mu, \xi).
$$
When $m\leq 0$, we obtain, in this way, an element of $\cnn$. We denote by
$\Psi_{AP}^0$ the algebra of pseudodifferential operators with almost
periodic coefficients and with non positive order.
When the order $m$ of $A$ is $>0$ then the operator $A^\sharp$
given by the family $(\sigma_\mu(A))_{\mu\in \R^n_B}$ is
affiliated with $\cnn$. If the order $m$ of $A$ is $<-n$, then the
bounded operator $A^\sharp$ is trace class with respect to the
trace $\tau$ on the von Neumann algebra $\cnn\otimes
M_N(\C)$\cite[Proposition 3.3]{Shubin2} and we  have:
$$
\tau (A^\sharp) = \lim_{T\to +\infty}\frac{1}{(2T)^n}\int_{(-T,+T)^n\times
\R^n} Tr (a(x,\zeta))  dx d\zeta
$$
Indeed, the expectation $E(A^\sharp)$ is a pseudodifferential operator on
$\R^n$ with symbol denoted by $E(a)$ and is independent of the
$x-$variable, it is given by:
$$
E(a) (\zeta) = \lim_{T\to +\infty} \frac{1}{(2T)^n} \int_{(-T,+T)^n}
a(x,\zeta) dx.
$$
Hence the operator $E(A^\sharp)$ is precisely the Fourier multiplier $\wM
(E(a))$ and so:
$$
\tau (A^\sharp) = \int_{\R^n} Tr( E(a) (\zeta)) d\zeta.
$$
Let $\Psi^\infty_{AP}$ be the space of one step polyhomogeneous classical
pseudodifferential operators on $\R^n$ with almost periodic coefficients.

\begin{thm}\label{Dixmier}\
Let $A$ be a (scalar) pseudodifferential operator with almost
periodic coefficients on $\R^n$. We assume that the order $m$ of
$A$ is  $\leq -n$ and we denote by $a_{-n}$ the $-n$ homogeneous
part of the symbol $a$. Then the operator $A^\sharp$ belongs to
the Dixmier ideal $\clL^{1,\infty}(\cnn, \tau)$. Moreover, the
Dixmier trace $\tau_\omega (A^\sharp)$  of $A^\sharp$ associated
with a limiting process $\omega$ does not depend on $\omega$ and
is given by the formula:
$$
\tau_\omega (A^\sharp) = \frac{1}{n} \int_{\R^n_B \times
\IS^{n-1}}
a_{-n} (x, \zeta) d\alpha_B(x) d\zeta.
$$
\end{thm}

\begin{proof}\
We denote as usual by $\Delta$ the Laplace operator on $\R^n$. The
operator $A(1+\Delta)^{n/2}$ is then a pseudodifferential operator
with almost periodic coefficients and nonpositive order. Hence,
the operator
$[A(1+\Delta)^{n/2}]^\sharp=A^\sharp(1+\Delta^\sharp)^{n/2}$
belongs to the von Neumann algebra $\cnn$. Now the operator
$(1+\Delta^\sharp)^{-n/2}$ is a Fourier multiplier defined by the
function $\zeta\mapsto (1+\zeta^2)^{-n/2}$. Hence if, for $\lambda
 >0$, $E_\lambda$ is the spectral projection of the operator
$(1+\Delta)^{-n/2}$ corresponding to the interval $(0,\lambda)$
then the operator $1-E_\lambda$ is a Fourier multiplier 
as well. Let us suppose it is defined
by the function $\zeta \mapsto
f_\lambda((\zeta^2+1)^{-n/2})$. It follows that the
trace $\tau$ of the operator $1-E_\lambda$ is given by
$$
\int_{\R^n} f_\lambda(\frac{1}{(\zeta^2+1)^{n/2}}) d\zeta.
$$
It is easy to compute this integral and to show that it is
proportional to $\frac{1}{\lambda}$. So the infimum of those
$\lambda$ for which $\tau (1-E_\lambda) \leq t$ is precisely
proportional to $\frac{1}{t}$. Hence the operator
$(1+\Delta^\sharp)^{-n/2}$, and hence $A$,  belongs to the Dixmier
ideal $\clL^{1,\infty}(\cnn, \tau)$.

In order to compute the Dixmier trace of the operator $A$, we
apply \cite[Theorem 10.1]{Shu} to deduce that the spectral
$\tau-$density $N_A(\lambda)$ of $A$ has the asymptotic expansion
$$
N_A(\lambda) = \frac{\gamma_0(A)}{\lambda} (1+ o(1)), \quad \lambda \to +\infty,
$$
where $\gamma_0(A)$ is given by:
$$
\gamma_0(A) = \frac{1}{n} \int_{\R^n_B\times \IS^{n-1}} a_{-n} (x, \zeta)
d\alpha_B(x) d\zeta.
$$
Now, if $A$ is positive then by \cite[Proposition 1]{BeF}:
$$
\tau_\omega (A) = \lim_{\lambda\to +\infty} \lambda N_A(\lambda) = \gamma_0(A).
$$
This proves the theorem for positive $A$. Since the principal symbol map is
a homomorphism, we deduce the result for general $A$.
\end{proof}

The normalisation we have chosen for the trace in the von Neumann
setting of this Section eliminates a factor of
$\frac{1}{(2\pi)^n}$ which occurs on the Wodzicki residue in the
type $I$ theory.

\subsection{Almost periodic spectral triple}

  We denote by $\cA$ the $*$-subalgebra
of $\CAP^{\infty}(\R^n)$ of smooth almost periodic functions on
$\R^n$. We take the Hilbert space on which the algebra $\cnn$ acts
to be $B^2(\R^n)\otimes L^2(\R^n)$ where $B^2(\R^n)$ is the
Hilbert space of almost periodic functions on $\R^n$ where the
norm and inner product are given by the restriction of the Haar
trace on $\CAP^{\infty}(\R^n)$ to $\cA$ (note that $B^2(\R^n)\cong
\ell^2(\R^n_d)$). This type $II_{\infty}$ von Neumann algebra is
endowed with a faithful normal semifinite trace that we denote by
$\tau$. (We note that the explicit formula for $\tau$ is as given
in the previous subsection.)

The usual Dirac operator on $\R^n$ is denoted by $\clD _0$. So, if
$\cS$ is the spin representation of $\R^n$ then $\clD _0$ acts on
smooth $\cS-$valued functions on $\R^n$. The operator $\clD _0$ is
$\Z^n-$periodic and it is affiliated with the von Neumann algebra
$\cnn_\cS = \cnn \otimes \End(\cS)$. This latter is also a type
$II_\infty$ von Neumann algebra with the trace $\tau\otimes \tr$.
More generally, for any $N\geq 1$, we shall denote by $\cnn_{\cS,
N}$ the von Neumann algebra $\cnn\otimes \End(\cS\otimes \C^N)$
with the trace $\tau\otimes \tr$.

The algebra $\cA$ and its closure are faithfully represented
as $*-$subalgebras of the von Neumann algebra $\cnn_\cS$.
In the same way the algebra $\cA\otimes M_N(\C)$ can be viewed
as a $*-$subalgebra of $\cnn_{\cS, N}$. More precisely, if $a\in \cA$ then
the operator $a^\sharp$ defined by
$$
(a^\sharp f)(x,y) := a(x+y) f(x,y), \quad \forall f\in
B^2(\R^n)\otimes L^2(\R^n),
$$
belongs to $\cnn_{S,N}$. The operator $a^\sharp$ is just the one associated
with
the zero-th order differential operator corresponding to multiplication by $a$.
The same formula allows to represent $\cA$ in $\cnn_\cS$.

\begin{prop}\
The triple $(\cA, \cnn_\cS, \clD _0^\sharp)$ is a semifinite
  spectral triple of finite dimension equal to $n$.
\end{prop}

\begin{proof}\
Note that the algebra $\cA$ is unital. The differential operator
$\clD _0$ is known to be elliptic periodic and self-adjoint on
$\R^n$. Therefore, the operator $\clD _0^\sharp$ is affiliated
with the von Neumann algebra $\cnn_\cS$ and it is self-adjoint as
a densely defined  unbounded operator on the Hilbert space
$B^2(\R^n)\otimes L^2(\R^n)$ with  $\clD _0^2 = \Delta\otimes Id$
with $\Delta$ the usual Laplacian. For any smooth almost periodic
function $f$ on $\R^n$, the commutator $[\clD _0, f]$ is a $0-$th
order almost periodic differential operator and so $[\clD
_0^\sharp, f]$ belongs to the von Neumann algebra $\cnn$.

On the other hand, the pseudodifferential operator $T=(\Delta + I)^{-1/2}$
is the Fourier multiplier associated with the function
$k\mapsto \frac{1}{(\|k\|^2 +1)^{1/2}}$. Therefore, its singular numbers
$\mu_t(T)$ can be computed explicitly as in the proof of Theorem \ref{Dixmier}
and shown to be proportional to $t^{-1/n}$.
\end{proof}

\section{Lesch's Index Theorem}

In this Section we describe
(following \cite{CPS2}) a proof of an index theorem due to M. Lesch
\cite{Le91JOT}
(see also \cite{PR}) for Toeplitz operators with noncommutative symbol
that relies on the zeta function approach to the Dixmier trace.
We begin with a  unital $C^*$-algebra $\mathcal A$
with a faithful finite trace,
$\tau_\A$ satisfying $\tau_\A(1)=1$ and a continuous action
  $\alpha$ of $\mbR$ on $\mathcal A$
leaving $\tau_\A$ invariant.

We let $H_{\tau_\A}$ denote the Hilbert space completion of $\mathcal A$
  in the inner
product $(a|b)=\tau(b^*a)$. Then
  $\mathcal A$ is a Hilbert algebra and the left regular
representation of $\mathcal A$
  on itself extends by continuity to a representation,
$a \mapsto \pi_{\tau_\A}(a)$ of $\mathcal A$
  on $H_{\tau_\A}$ \cite{Dix}. In what follows, we will drop
the notation $\pi_{\tau_\A}$ and just denote the action of $\mathcal A$
  on $H_{\tau_\A}$
by juxtaposition.

We now look at the induced representation,
$\tilde{\pi}$, of the crossed product
$C^*$-algebra
$\clA \times_{\alpha}\mbR$ on $L^2(\mbR,H_{\tau_\A})$.
That is, $\tilde{\pi}$ is the
representation $\pi\times\lambda$ obtained from the covariant pair,
$(\pi,\lambda)$ of
representations of the system $(\clA ,\mbR,\alpha)$ defined
  for $a \in \clA $,
$t,s \in \mbR$ and $\xi \in L^2(\mbR,H_{\tau_\A})$ by:
$$(\pi(a)\xi)(s) = \alpha_s^{-1}(a)\xi(s)$$
and
$$\lambda_t(\xi)(s) = \xi(s-t).$$

Then, for a function $x\in L^1(\mbR,\clA )
\subset \clA \times_{\alpha}\mbR$,
  the action of $\tilde{\pi}(x)$ on a vector $\xi$
in $L^2(\mbR,H_{\tau_\A})$ is defined as follows:
$$(\tilde{\pi}(x)\xi)(s) =
\int_{-\infty}^{\infty}\alpha_s^{-1}(x(t))\xi(s-t)dt.$$

Now the twisted convolution algebra
$L^1(\mbR,\clA ) \cap L^2(\mbR,H_{\tau})$ is
a dense subspace of $L^2(\mbR,H_{\tau})$ and also a Hilbert algebra in the
given inner product. As such, there is a canonical faithful, normal,
semifinite trace, $\tau$,
on the von Neumann algebra
$$\mathcal N=
(\tilde{\pi}(\clA \times_{\alpha}\mbR))^{\prime\prime}.$$
For functions $x,y: \mbR\to \clA \subset H_{\tau_\A}$ which are in
$L^2(\mbR,H_{\tau_\A})$ and whose twisted left convolutions $\tilde{\pi}(x),
\tilde{\pi}(y)$
define bounded operators on $L^2(\mbR,H_{\tau_\A})$, this trace is given by:
$$\tau(\tilde{\pi}(y)^*\tilde{\pi}(x))= \langle x|y\rangle =
\int_{-\infty}^{\infty} \tau_\A(x(t)y(t)^*)dt.$$

In particular, if we identify
$L^2(\mbR)=L^2(\mbR)\otimes 1_\clA \subset L^2(\mbR,H_{\tau_\A})$
  then any scalar-valued
function $x$ on $\mbR$ which is the Fourier transform $x=\widehat{f}$ of a
bounded
$L^2$ function, $f$ will have the properties that $x\in L^2(\mbR,H_{\tau_\A})$
and $\tilde{\pi}(x)$ is a bounded operator. For such scalar functions $x$, the
operator $\tilde{\pi}(x)$ is just the usual convolution by the function $x$
and is
usually denoted by $\lambda(x)$ since it is just the integrated form of
$\lambda$. The next Lemma follows easily from these considerations.

\begin{lemma}
With the hypotheses and notation discussed above

(i) if $h\in L^2(\mbR)$ with $\lambda(h)$ bounded and $a\in \mathcal A$,
then defining $f:\mbR\to H_{\tau}$ via $f(t)=ah(t)$ we see that
$f\in L^2(\mbR,H_{\tau_\A})$ and $\tilde{\pi}(f)=\pi(a)\lambda(h)$ is bounded,

(ii) if $g\in L^1(\mbR)\cap L^{\infty}(\mbR)$ and $a\in \mathcal A$ then
$\pi(a)\lambda(\hat{g})$ is trace-class in $\mathcal N$ and
$$Tr(\pi(a)\lambda(\hat{g}))=\tau_\A(a)\int_{-\infty}^{\infty}g(t)dt.$$
\end{lemma}

\begin{proof}
To see part (i), let $\xi\in C_c(\mbR,H_{\tau_\A})\subseteq
L^2(\mbR,H_{\tau_\A}).$
Then
\begin{eqnarray*}
(\tilde{\pi}(f)\xi)(s)&=&\int_{-\infty}^{\infty}\alpha_s^{-1}(f(t))\xi(s-t)dt\\
&=&\int_{-\infty}^{\infty}\alpha_s^{-1}(a)h(t)\xi(s-t)dt\\
&=&\alpha_s^{-1}(a)\int_{-\infty}^{\infty} h(t)\xi(s-t)dt\\
&=&\alpha_s^{-1}(a)(\lambda(h)\xi)(s)\\
&=&(\pi(a)\lambda(h)\xi)(s).
\end{eqnarray*}

To see part (ii) we can assume that $g$ is nonnegative and $a$ is
self-adjoint. Then let
$g=g^{1/2}g^{1/2}$ so that $g^{1/2}\in L^2\cap L^{\infty}$ and so
$\lambda(\widehat{g^{1/2}})$ is bounded. Now,
$$\pi(a)\lambda(\widehat{g}
)=\pi(a)\lambda(\widehat{g^{1/2}})\pi(1_\clA )\lambda
(\widehat{g^{1/2}}).$$
Then, $\pi(a)\lambda(\widehat{g^{1/2}})=\tilde{\pi}(x)$ where $x(t)=
a\widehat{g^{1/2}}(t)$
and $\pi(1_\clA )\lambda(\widehat{g^{1/2}})=\tilde{\pi}(y)$
  where $y(t)=
1_\clA \widehat{g^{1/2}}(t).$
So, $\tilde{\pi}(x)$ and $\tilde{\pi}(y)$ are in $\clN _{sa}$ and
$\pi(a)\lambda(\widehat{g})=\tilde{\pi}(x)\tilde{\pi}(y).$

Hence,
\begin{eqnarray*}
\tau(\pi(a)\lambda(\widehat{g}))&=&\tau(\tilde{\pi}(x)\tilde{\pi}(y))\\
&=&\int_{-\infty}^{\infty}\tau(x(t)y(t))dt\\
&=&\tau_\A(a)\int_{-\infty}^{\infty}\left|\widehat{g^{1/2}}(t)\right|^2dt
=\tau_\A(a)\int_{-\infty}^{\infty} g(s)ds.
\end{eqnarray*}\end{proof}

By construction
$\mathcal N$ is a semifinite von Neumann algebra with faithful, normal,
semifinite
trace, $\tau$, and a faithful representation $\pi:\clA \to \mathcal N$
\cite{Dix}.
For each $t\in \mbR$, $\lambda_t$ is a unitary in
$U(\mathcal N)$. In fact the one-parameter unitary group
$\{\lambda_t\; \vert\; t\in\mbR\}$ can be written
$\lambda_t=e^{it\D}$ where $\D$  is
the unbounded
self-adjoint operator $$\D=\frac{1}{2\pi i}\frac{d}{ds}$$ which is affiliated
with $\mathcal N$. In the Fourier transform
picture (i.e., the spectral picture for $\D$)
of the previous proposition, $\D$
becomes multiplication by the independent variable and so $f(\D)$ becomes
pointwise multiplication by the function $f$. That is,
$$\tilde{\pi}(\hat{f}) = \lambda(\hat{f}) = f(\D).$$
And, hence, if $f$ is a bounded $L^1$ function, then:
$$\tau(f(\D)) = \int_{-\infty}^{\infty}f(t)dt.$$
By this discussion and the previous lemma, we have the following result

\begin{lemma}
If $f\in L^1(\mbR)\cap L^{\infty}(\mbR)$ and $a\in \mathcal A$
  then $\pi(a)f(\D)$ is
trace-class in $\mathcal N$ and
$$\tau(\pi(a)f(\D))=\tau(a)\int_{-\infty}^{\infty}f(t)dt.$$
\end{lemma}

  We let $\delta$ be the densely defined (unbounded) $*$-derivation
on $\mathcal A$ which is the infinitesimal generator of the representation
$\alpha:\mbR\to Aut(\clA )$
  and let $\hat{\delta}$ be the unbounded $*$-derivation
on $\mathcal N$ which is the infinitesimal generator of the
representation $Ad\circ\lambda:\mbR\to Aut(\mathcal N)$
(here $Ad(\lambda_t)$ denotes conjugation by $\lambda_t^*\cdot\lambda_t$).
Now if
$a\in dom(\delta)$ then clearly $\pi(a)\in dom(\hat{\delta})$ and
$\pi(\delta(a)) = \hat{\delta}(\pi(a))$. By \cite{BR} Proposition 3.2.55
(and its proof) we have that $\pi(\delta(a))$ leaves the domain of $\D$
invariant and  $$\pi(\delta(a))=2\pi i[\D,\pi(a)].$$
We are now in a position to state and prove Lesch's index theorem.

\begin{thm}[\cite{CPS2}]
Let $\tau_\A$
be a faithful finite trace on the unital $C^*$-algebra, $\mathcal A$,
which is invariant for an action $\alpha$ of $\mbR$. Let $\mathcal N$ be the
semifinite von Neumann algebra
$(\tilde{\pi}(\clA \times_{\alpha}\mbR))^{\prime\prime}$,
and let $\D$ be the
infinitesimal generator of the canonical representation $\lambda$ of $\mbR$
in $U(\mathcal N)$. Then, the representation
$\pi:\clA \to \mathcal N$ defines a
$\clL ^{(1,\infty)}$ summable Breuer-Fredholm module
$(\clN ,\D)$ for
  $\mathcal A$. Moreover,
if $P$ is the nonnegative spectral projection for $\D$ and
  $u\in U(\clA )$ is also
in the domain of $\delta$, then $T_u:=P\pi(u)P$ is Breuer-Fredholm in
$P\mathcal NP$ and
$$ind (T_u) = \frac{1}{2\pi i}\tau(u\delta(u^*)).$$
\end{thm}

\begin{proof}
It is easy to see that $\D$ satisfies
$(1+\D^2)^{-1/2}\in \clL ^{(1,\infty)}$. By the previous discussion,
for any $a\in dom(\delta)$ we have $\pi(\delta(a))=2\pi i[\D,\pi(a)].$
Since the
domain of $\delta$ is dense in $\mathcal A$ we see that $\pi$ defines a
$\clL ^{(1,\infty)}$ summable Breuer-Fredholm module for $\mathcal A$.

By Theorem~\ref{GTOind},
$$ind (T_u) =\lim_{p\to 1^+}\frac{1}{2}(p-1)
Tr(\pi(u)[\D,\pi(u^*)](1+\D^2)^{-p/2})$$
and hence
\begin{eqnarray*}
ind (T_u)&=&\lim_{p\to 1^+}\frac{1}{2}(p-1)\frac{1}{2\pi i}
Tr(\pi(u\delta(u^*))(1+\D^2)^{-p/2})\\
&=&\lim_{p\to 1^+}\frac{1}{2}(p-1)\frac{1}{2\pi i}
\tau(u\delta(u^*))\int_{-\infty}^{\infty}(1+t^2)^{-p/2}dt\\
&=&\lim_{p\to 1^+}\frac{1}{2\pi i}\tau(u\delta(u^*))\frac{1}{2}(p-1)
{C}_{p/2} = \frac{1}{2\pi i}\tau(u\delta(u^*))
\end{eqnarray*}
\end{proof}

\section{The Hochschild Class of the Chern Character}\label{MR}

The theorem we discuss in this Section, for $\clN=\clL(\HH)$ and
$1<p<\infty$ ($p$ integral), was proved in lectures by Alain
Connes at the Coll\`{e}ge de France in 1990. A version of this
argument appeared in \cite{GVF}. The extension of this argument to
general semifinite von Neumann algebras, with the additional
hypothesis that the unbounded self-adjoint operator $\D$ (which
will form part of a spectral triple) have bounded inverse, is
presented by Benameur and Fack, \cite{BeF}. A simpler strategy
using the pseudodifferential calculus of Connes-Moscovici,
\cite{CM}, was communicated to us by Nigel Higson. In conjunction
with the results in \cite{CPS2}, Higson's argument appears to
generalise to the semifinite case, however, we will not describe
the details here focusing instead on another approach.

In \cite{CPRS1} the theorems of Connes and Benameur-Fack were extended.
First a proof was given for the case $p=1$
and second the hypothesis, in the type $II_\infty$ case, that $\D$
has bounded inverse was removed. This is crucial due to the
`zero-in-the-spectrum' phenomenon for $\D$. That is, for type $II$
$\clN$, zero is generically in the point and/or continuous
spectrum, \cite{FW} and hence we are not dealing with just the
simple problem that arises in the type $I$ case posed by a finite
dimensional kernel. One feature of the approach in \cite{CPRS1} is
that the strategy of the proof is the same for all $p\geq 1$, and
is independent of the type of the von Neumann algebra $\clN$.

We explain the general semifinite version of the type $I$ result
in \cite[IV.2.$\gamma$]{Co4} which identifies the Hochschild class
of the Chern character of a $(p,\infty)$-summable spectral triple.

\begin{thm}\label{thm8} Let $(\A,\HH,\D)$ be a $QC^k$
$(p,\infty)$-summable spectral triple with $p\geq 1$ integral and
$k=\max\{2,p-2\}$. Then

1) A Hochschild cocycle on $\A$ is defined by
\ben \phi_\omega(a_0,...,a_p)=\lambda_p\tau_\omega(\Gamma
a_0[\D,a_1]\cdots[\D,a_p](1+\D^2)^{-p/2}),\een

2) For all Hochschild $p$-cycles $c\in C_p(\A)$ (i.e., $bc=0$),
\ben\la \phi_\omega,c\ra=\la Ch_{F_\D},c\ra,\een
where $Ch_{F_\D}$ is the Chern character in cyclic cohomology of the
pre-Fredholm module
over $\A$ with $F_\D=\D(1+\D^2)^{-1/2}$. \end{thm}

\begin{rems}
Here $\tau_\omega$ is the Dixmier trace associated to any state
  $\omega\in D(\mbR_+^*)$.
\end{rems}
The two most important corollaries of Theorem \ref{thm8} are the following.
\begin{cor}\label{ismeas} Let $(\A,\HH,\D)$ be as in Theorem \ref{thm8}.
If $c=\sum_i a_0^i\otimes a_1^i\otimes\cdots \otimes
a_p^i$ is a Hochschild $p$-cycle, then
\ben \Gamma\sum_ia_0^i[\D,a_1^i]\cdots[\D,a_p^i](1+\D^2)^{-p/2}\een
is (Dixmier)-measurable.
\end{cor}
This corollary is relevant to the axioms of noncommutative spin geometry
since it tells us that when we use the noncommutative integration
theory provided by the Dixmier trace we do not need to worry
which functional $\omega$ we use to define the integration when we apply it
to Hochschild cycles.
\begin{cor} With $(\A,\HH,\D)$ as in Theorem \ref{thm8}, and supposing
that $Ch_{F_\D}$ pairs nontrivially with $HH_p(\A)$, then
\ben \tau_\omega((1+\D^2)^{-p/2})\neq 0.\een
\end{cor}

\begin{rems}
The hypothesis of the Corollary is that there exists some
Hochschild $p$-cycle such that $\la ICh_{F_\D},c\ra\neq 0$
(where $I$ is the map defined by Connes, \cite{Co4} subsection III.1.$\gamma$
pp 199-200). Computing this
pairing using Theorem \ref{thm8} above, we see that  $(1+\D^2)^{-p/2}$ can
not have zero Dixmier trace for any  choice of Dixmier functional $\omega$.
For if $(1+\D^2)^{-p/2}$ did have vanishing Dixmier trace,
and $c=\sum_ia^i_0\otimes\cdots\otimes a^i_p$ is any Hochschild cycle
\bean |\la ICh_{F_\D},c\ra|&=&\left|\sum_i\tau_\omega
\left(\Gamma a^i_0[\D,a^i_1]\cdots[\D,a^i_p](1+\D^2)^{-p/2}\right)
   \right|\nno
&\leq&\sum_i\n\Gamma a^i_0[\D,a^i_1]\cdots[\D,a^i_p]\n\tau_\omega
\left((1+\D^2)^{-p/2}\right)=0.\eean
Hence if the pairing is nontrivial, the Dixmier trace can not vanish
on $(1+\D^2)^{-p/2}$.
\end{rems}
There are subtle points in the proof. One first assumes that the triple
   triple $(\A,\HH,\D)$ has $\D$ invertible (by replacing $(\A,\HH,\D)$ by
$(\A,\HH^2,\D_m)$ with $\D_m=\bma \D & m\\ m & -\D\ema$
if necessary).
Then one needs to verify that
$\phi_\omega$ in Theorem \ref{thm8} does
not depend on this replacement.
One also requires the verification that the functional
$\phi_\omega$ is indeed a Hochschild cocycle but this quite simple.

\begin{lemma}\label{iscocycle} Let $p\geq 1$
and suppose that $(\A,\HH,\D)$ is a $QC^1$
$(p,\infty)$-summable spectral triple. Then the multilinear functional
\ben \phi_\omega(a_0,...,a_p)=
\lambda_p\tau_\omega(\Gamma a_0[\D,a_1]\cdots[\D,a_p](1+\D^2)^{-p/2})\een
is a Hochschild cocycle.
\end{lemma}

\begin{proof} {By Lemma 3 of \cite{CPRS1}
and the trace property of the Dixmier trace, we
have
\bean (b\phi_\omega)(a_0,...,a_p)&=&(-1)^{p-1}\lambda_p
\tau_\omega(\Gamma a_0[\D,a_1]
\cdots[\D,a_{p-1}]a_p(1+\D^2)^{-p/2})\nno
&&\qquad-(-1)^{p-1}\lambda_p\tau_\omega(\Gamma a_0[\D,a_1]
\cdots[\D,a_{p-1}](1+\D^2)^{-p/2}a_p).\eean
As $(\A,\HH,\D)$ is $QC^1$,
\bean [(1+\D^2)^{-p/2},a_p]
&=&-\sum_{k=0}^{p-1}(1+\D^2)^{-(p-k)/2}[(1+\D^2)^{1/2},a_p]
(1+\D^2)^{-(1+k)/2},\eean
and this is trace class. So  $a_p(1+\D^2)^{-p/2}=(1+\D^2)^{-p/2}a_p$
modulo trace class operators, and so
the two terms above cancel.}
\end{proof}
}

\section{Lidskii type formula for Dixmier traces}

A semifinite analogue of the classical Lidskii theorem stated in terms
of the (so-called) Brown spectral measure $\mu_T$ of $T\in \clN$
asserts~\cite{Brown} that
$$
\tau(T)=\int_{\sigma(T)\setminus\{0\}}\lambda d\mu_T(\lambda).
$$
In the case, when $\clN=\clL(\clH)$ and $\tau$ is the standard
trace $Tr$ the equality above reduces to the classical case
asserting that the trace $Tr(T)$ of an arbitrary  trace class
operator $T$ is given by the sum $\sums{n\geq 1} \lambda(T),$
where $\{\lambda_n(T)\}_{n\geq 1}$ is the sequence of eigenvalues
of $T,$ arranged in decreasing order of absolute values of
$\lambda_n$ and counting multiplicities. Note, that in the case
$T\geq 0,$ the equality $\Tr(T) = \sums{n\geq 1} \lambda_n(T)$
follows immediately from the spectral theorem for compact
operators. If $T^*=T,$ then again, by the spectral theorem we can
select the orthonormal basis of $\clH$ consisting of eigenvalues
of $T$ and still infer
Lidskii's theorem without any difficulty. Here
it is worth observing that the assumption $T^*=T$ belongs to the ideal
$\clL^1(\clH)$ 
of all trace class operators on $\clH$
 implies the absolute convergence of the series
$\sums{n\geq 1} |\lambda_n(T)|.$ The latter fact guarantees the
convergence of the series $\sums{n\geq 1}\lambda_n(T)$ in whatever
ordering of the set of all eigenvalues for $T$ is chosen (in
particular, for the decreasing ordering of absolute values of $T$).

The core difference of this situation with the setting of Dixmier
traces consists in the fact that the series $\sums{n\geq 1}
|\lambda_n(T)|$ diverges for every normal $T\in\DixIdeal{1}\setminus
\clL^1(\clH).$ Therefore, even though for a given
$T=T^*\in\DixIdeal{1}$ we define $\tau_\omega(T)$ as the difference
$\tau_\omega(T_+)-\tau_\omega(T_-)$ where each number is computed
according to the definition $\limo{N} \ILogN
\sums{n=1}^N\lambda_n(T_\pm),$ it is by no means clear that we have
$$
   \tau_\omega(T)=\limo{N} \ILogN \sums{n=1}^N\lambda_n(T)
$$
for the special enumeration of the set $\{\lambda_n(T)\}_{n\geq
1}$ given by the decreasing order of absolute values of
$|\lambda_n(T)|;$ or for that matter for \emph{any} enumeration of
this set. This difficulty becomes even more pronounced in the case
of a general semifinite von Neumann albegra.


Note that the restriction $\mu_t(T)\le\frac{C}{t}$, $t\ge1$
imposed on $T\in \clL^{(1,\infty)}\nt$ in the theorem below, implies
  that $T$ belongs to the ideal $N(\psi_1)\nt$ (see the definition
of $\psi_1$ in Section~\ref{MFandSS} and the definition of r.i.
ideal $N(\psi)$ in Section~\ref{SSFonMS}). Everywhere below, we
assume that $\omega \in CD(\mbR_+^*)$ (see
Section~\ref{Connes-Dixmier traces}). When we apply $\omega$ to a
sequence $x\in\ell_\infty$, we identify $x$ with its image $i(x)$
in  $L_\infty(\mbR_+)$ (see Remark~\ref{IsometricEmbeddingRem}).
\begin{thm}[\cite{AzSu}]\label{LidskiiForDixmierTr}
\begin{thlist}
\item If $\nt$ is a semifinite von Numann albegra and $T\in
   \clL^{(1,\infty)}\nt$ satisfies $\mu_t(T)\le\frac{C}{t}$,
   $t\ge1$ for some $C>0,$
   then
$$
\tau_\omega
(T)=\omega\text{-}\lim_{t\to\infty}\frac{1}{\log(1+t)}\int_{\lambda\notin\frac{1}{t}G}\lambda
d\mu_T(\lambda).
$$
\item
   Let $T$ be a compact operator on a Hilbert space $\clH,$ such that
   $\mu_n(T)\leq \frac{C}{n}, \ n\geq 1$ for some $C>0.$ Let
   $\lambda_1,\,\lambda_2,\,\dots$ be the list of eigenvalues of the
   operator $T$ counted with multiplicities such that
   $|\lambda_1|\geq|\lambda_2|\geq\dots\ .$ Then
\begin{equation}\label{TRofIdealeq}
\tau_\omega(T)= \limo{t} \ILog \sums{\lambda \in \sigma(T), \
   \lambda\notin \frac 1t G} \lambda\mu_T(\lambda) = \limo{N}\frac
1{\log (1+N)} \sums{i=1}^N \lambda_i,
\end{equation}
where $\mu_T(\lambda)$ is the algebraic multiplicity of the eigenvalue
$\lambda$ and $G$ is an arbitrary bounded neighborhood of $0\in\mbC.$

\end{thlist}

\end{thm}

Here, we emphasize the fact that $G$ is an \emph{arbitrary}
bounded neighborhood of $0.$ In the case $\clN=\clL(\clH)$,
formula~\eqref{TRofIdealeq} says that we need to compute the sum
of all eigenvalues $\lambda_j(T)$ which do not belong to the
"squeezed" neighborhood $\frac 1t G$ as $t\to \infty$ and then to
apply $\omega$-limit.  It is interesting to emphasize the special
case of formula~\eqref{TRofIdealeq} for measurable operators
belonging to the ideal $N(\psi_1)$. The result below should be
compared with the results of Theorems~\ref{th6.5}
and~\ref{thTisCDmeasurable}.
\begin{cor}
If $T$ is a (Connes-Dixmier) measurable operator satisfying the
assumption of Theorem~\ref{LidskiiForDixmierTr}~(i) (respectively, (ii)), then
$$
\tau_\omega(T)= \lim_{t\to\infty} \ILog \int_{\lambda
   \not\in\frac{1}{t}G} \lambda d\mu_T(\lambda) \quad\left(\text{resp.}=
\lim_{N\to\infty}\frac 1{\log (1+N)} \sums{i=1}^N \lambda_i\right),
$$
for every $\omega\in CD(\mbR_+^*)$.
\end{cor}

Fix an orthonormal basis in $\clH$ and identify every element
$x\in \clL(\clH)$ with its matrix $(x_{ij})_{i,j=1}^\infty$. It is
well-known~\cite{GK} that the triangular truncation operator
$\frT$ given by
$$\frT(x)=\left\{\begin{array}{cc}
  x_{ij}, & i\geq j\\ 0, & i<j \end{array}\right.$$
acts boundedly from the trace class $\clL^1(\clH)$ into
$\clL^{(1,\infty)}(\clH)$. Noting that $\frT(x)-\diag(x)$ is quasinilpotent for
every $x\in \clL^1(\clH)$, we obtain the following


\begin{cor} The operator $\frT(x)$ is Connes-Dixmier measurable for every
$x\in\clL^1(\clH)$,
moreover $\tau_\omega(\frT(x))=0$, $\omega \in CD(\mbR_+^*)$.
\end{cor}
We conclude with the following result immediately following from
Theorem~\ref{LidskiiForDixmierTr} (see
also~\cite[Proposition~1]{Fa04FA}).
\begin{cor}
   Let $M$ be a compact $n$-dimensional Riemannian manifold and let $T$
   be a pseudodifferential operator of order $-n$ on $M.$ Then
  $$
    \tau_\omega(T) = \liminfty{N} \frac 1{\log N} \sums{k=1}^N \lambda_k.
  $$
\end{cor}

\mathsurround 0pt
\def\itakdalee{$\dots$}


\end{document}